\documentclass[%
 aip,
 amsmath,amssymb,
 reprint,%
]{revtex4-1}

\usepackage{graphicx}
\usepackage{dcolumn}
\usepackage{bm}
\usepackage{color}
\usepackage{hyperref}

\newtheorem{theorem}{Theorem}[section]

\newtheorem{remark}[theorem]{Remark}

\newcommand{\gsyn}{g_{\rm syn}}
\newcommand{\ms}{M_S}
\newcommand{\ls}{L_s}
\newcommand{\mss}{M_{SS}}

\newcommand{\fsn}{\mathrm{FSN}}

\renewcommand{\epsilon}{\varepsilon}

\newcommand{\RED}[1]{{\color{black}{#1}}}

\newcommand{\subfigimg}[3][,]{%
	\setbox1=\hbox{\includegraphics[#1]{#3}}
	\leavevmode\rlap{\usebox1}
	\rlap{\hspace*{-5pt}\raisebox{\dimexpr\ht1-2\baselineskip}{#2}}
	\phantom{\usebox1}
}
\begin{document}


\title{Mixed Mode Oscillations in a Three-Timescale Coupled Morris-Lecar System} 



\author{Ngoc Anh Phan}
\email[]{ngocanh-phan@uiowa.edu}
\affiliation{Department of Mathematics, University of Iowa, Iowa City, IA 52242}

\author{Yangyang Wang}
\email[]{yangyangwang@brandeis.edu}
\affiliation{Department of Mathematics, Brandeis University, Waltham, MA 02453}


\begin{abstract}
Mixed mode oscillations (MMOs) are complex oscillatory behaviors of multiple-timescale dynamical systems in which there is an alternation of large-amplitude and small-amplitude oscillations. It is well known that MMOs in two-timescale systems can arise either from a canard mechanism associated with folded node singularities or a delayed Andronov-Hopf bifurcation (DHB) of the fast subsystem. While MMOs in two-timescale systems have been extensively studied, less is known regarding MMOs emerging in three-timescale systems. In this work, we examine the mechanisms of MMOs in coupled Morris-Lecar neurons with three distinct timescales. We investigate two kinds of MMOs occurring in the presence of a singularity known as canard-delayed-Hopf (CDH) and in cases where CDH is absent. In both cases, we examine how features and mechanisms of MMOs vary with respect to variations in timescales. Our analysis reveals that MMOs supported by CDH demonstrate significantly stronger robustness than those in its absence. Moreover, we show that the mere presence of CDH does not guarantee the occurrence of MMOs. This work yields important insights into conditions under which the two separate mechanisms in two-timescale context, canard and DHB, can interact in a three-timescale setting and produce more robust MMOs, particularly against timescale variations. 
\end{abstract}


\maketitle 

\begin{quotation}
One of the most common types of complex oscillatory dynamics observed in systems with multiple timescales is mixed mode oscillations (MMOs). MMOs are characterized by patterns that involve the interspersion of small-amplitude and large-amplitude oscillations. Over the years, the theory of MMOs in fast-slow systems has been well-developed. Recently, there has been more progress on the analysis of MMOs in three-timescale systems. Nonetheless, MMOs in the latter case are still much less understood. In this work, we contribute to the investigation of MMOs in the three-timescale settings by considering coupled Morris-Lecar neurons. We uncover the properties and geometric mechanisms underlying two different MMO patterns in our three-timescale system. One of them involves the interaction of the two distinct MMO mechanisms, showing a high degree of robustness to timescale perturbations, whereas the other lacks such mechanism and is thus vulnerable to timescale variations. Based on our analysis, we establish conditions that lead to more robust generation of MMOs in three-timescale problems, particularly against perturbations in timescales.
\end{quotation}

\section{Introduction}\label{sec:intro}
Mixed mode oscillations (MMOs) are frequently perceived in the dynamical systems involving multiple timescales \citep{Desroches2012}; these are complex oscillatory dynamics characterized by the concatenation of small-amplitude oscillations (SAOs) and large-amplitude excursions in each periodic cycle. Such phenomena have been recognized in many branches of sciences including physics, chemistry \citep{Hudson1979,Awal2023}, and particularly life sciences such as \citep{Krupa2008a,Yu2008,Krupa2012,Teka2012, Vo2010, Vo2014,Curtu2010,CR2011,Harvey2011,Kugler2018,Kimrey2020,2ndKimrey2020,Pavlidis2022,Bat2021}. 

Theoretical analysis of MMOs in systems with two distinct timescales has been well developed with the implementation of the geometric singular perturbation theory (GSPT) \cite{Fenichel1979}; see \cite{Desroches2012} for review. Two common mechanisms leading to the occurrence of MMOs in multiple timescale problems are canard dynamics associated with the twisting of slow manifolds due to folded singularities \citep{SW2001,Wechselberger2005} and a slow passage \RED{through} the delayed Andronov-Hopf bifurcation (DHB) of the fast subsystem \citep{Baer1989, Neishtadt1987, Neishtadt1988, Hayes2016}. While in two-timescale settings, these two mechanisms remain separated, they can coexist and interact in three-timescale regime \citep{Teka2012, Vo2013, Maess2014,Letson2017}. 

Compared with the extremely well-studied MMOs in two-timescale problems, the theory of MMOs in the three-timescale settings has been less well-developed. Traditionally, three-timescale problems are simplified to two-timescale problems which is the natural setting for geometric singular perturbation theory \citep{Baldemir2020}. However, many real-world systems have more than two timescales \citep{WR2016, WR2017, WR2020, Maess2014,Jalics2010, Vo2013, Chumakov2015, PW2014}. It has also been established that a two-timescale decomposition fails to capture certain aspects of the system's dynamics \citep{Nan2015}. Therefore, classifying three timescales into two groups is not a sufficient approach for modelling and analysis. 

MMOs in three-timescale systems have been studied before (see, e.g., \cite{Krupa2008a, Krupa2008b, Krupa2012, Vo2013, Maess2014, Maess2016, Letson2017, DK2018, Kak2022, Jalics2010, Kak2023a, Kak2023b, Sadhu2022}, for examples and references). 
Initial approaches were to consider three-dimensional systems
\begin{equation}
    \begin{array}{rcl}
       \varepsilon \frac{dx}{dt}  &=& f(x,y,z), \vspace{0.1in}  \\
        \frac{dy}{dt} &=& g(x,y,z), \vspace{0.1in}\\
        \frac{dz}{dt} &=& \delta h(x,y,z),
    \end{array}
\end{equation}
with special cases $\delta=\varepsilon$ or $\delta=\sqrt{\varepsilon}$ \citep{Krupa2008a,Krupa2008b,Maess2014,Jalics2010,Maess2016}. MMOs were shown to emerge through an effect analogous to a slow passage through a canard explosion \citep{Krupa2008b,Jalics2010,Maess2014}.
More recently, there has been a growing interest in MMOs with independent singular perturbation parameters $\varepsilon$ and $\delta$, as explored in various three-dimensional models \citep{Letson2017,Kak2022,Kak2023a, Kak2023b}. In particular, Ref.~\onlinecite{Letson2017} centered on a novel singularity type denoted as \textit{canard-delayed-Hopf} (CDH) singularity, which naturally arises in three-timescale settings when the two mechanisms for MMOs (the fast subsystem Hopf and a folded node) coexist and interact. The authors investigated the existence and properties of MMOs near the CDH singularity. 

In this paper, we contribute to the investigation of MMOs in three-timescale settings by considering a model of 4-dimensional coupled Morris-Lecar neurons \citep{ML1981,RE1998} that was introduced by \cite{Nan2015}. 
The model equations are given by
\begin{widetext}
\begin{equation}\label{eq:main}
    \begin{array}{rcl}
        \frac{dV_1}{dt} &=&(I_1-g_{\rm Ca}m_\infty(V_1)(V_1-V_{\rm Ca})-g_{\rm K}w_1(V_1-V_{\rm K})-g_{\rm L}(V_1-V_{\rm L})-g_{\rm syn}S(V_2)(V_1-V_{\rm syn}))/C_1, \vspace{0.1in} \\
        \frac{dw_1}{dt}&=&\phi_1(w_\infty(V_1)-w_1)/\tau_{w}(V_1), \vspace{0.1in}\\
        \frac{dV_2}{dt}& =&(I_2-g_{\rm Ca}m_\infty(V_2)(V_2-V_{\rm Ca})-g_{\rm K}w_2(V_2-V_{\rm K})-g_{\rm L}(V_2-V_{\rm L}))/C_2, \vspace{0.1in}\\
        \frac{dw_2}{dt}&=&\phi_2(w_\infty(V_2)-w_2)/\tau_{w}(V_2),
    \end{array}
\end{equation}
\end{widetext}
with
\begin{equation}\label{eq:term}
    \begin{array}{rcl}
         S(V_i)&=&\alpha(V_i)/(\alpha(V_i)+\beta), \vspace{0.1in}\\
         \alpha(V_i)&=&1/(1+{\rm exp}(-(V_i-\theta_s)/\sigma_s)), \vspace{0.1in}\\
         m_{\infty}(V_i)&=&0.5(1+{\rm tanh}((V_i-K_1)/K_2)), \vspace{0.1in}\\
         w_{\infty}(V_i)&=&0.5(1+{\rm tanh}((V_i-K_{3})/K_{4})), \vspace{0.1in}\\
         \tau_{w}(V_i)&=&1/{\rm cosh}((V_i-K_{3})/2K_{4}).
    \end{array}
\end{equation}

\begin{table}[!t]
    \footnotesize
    \caption{\label{tab:par}The values of the parameters in the model given by  (\ref{eq:main}) and (\ref{eq:term}).}
    \begin{center}
        \begin{tabular}{|c c|c c|c c|}
            \hline
            \multicolumn{6}{|c|}{Parameter values} \\
            \hline
            $C_1$\hspace{3mm} & $8\hspace{1mm} {\rm \mu F/cm^2}$ & $I_1$\hspace{3mm} & $0\hspace{1mm} {\rm \mu A/cm^2}$ & $\phi_1$\hspace{3mm} & $0.01$\\

            $C_2$\hspace{3mm} & $100\hspace{1mm} {\rm \mu F/cm^2}$ & $I_2$\hspace{3mm} & $60\hspace{1mm} {\rm \mu A/cm^2}$ & $\phi_2$\hspace{3mm} & $0.001$\\

            $V_{\rm Ca}$\hspace{3mm} &$120\hspace{1mm} {\rm mV}$ & $g_{\rm Ca}$\hspace{3mm} & $4\hspace{1mm} {\rm mS/cm^2}$ & $K_1$\hspace{3mm} & $-1.2\hspace{1mm} \rm mV$\\

            $V_{\rm K}$\hspace{3mm} & $-84\hspace{1mm} \rm mV$ & $g_{\rm K}$\hspace{3mm} & $8\hspace{1mm} \rm mS/cm^2$ & $K_2$\hspace{3mm} & $18\hspace{1mm} \rm mV$\\

            $V_{\rm L}$\hspace{3mm} & $-60\hspace{1mm} \rm mV$ & $g_{\rm L}$\hspace{3mm} & $2\hspace{1mm} \rm mS/cm^2$ & $K_{3}$\hspace{3mm} & $12\hspace{1mm} \rm mV$\\

            $V_{\rm syn}$\hspace{3mm} & $30\hspace{1mm} \rm mV$& $\theta_s$\hspace{3mm} & $-20\hspace{1mm} \rm mV$ & $K_{4}$\hspace{3mm} & $17.4\hspace{1mm} \rm mV$\\

            $\beta$\hspace{3mm} & $0.5\hspace{1mm} \rm ms^{-1}$ & $\sigma_s$\hspace{3mm} & $10\hspace{1mm} \rm mV$ & &   \\
            \hline
        \end{tabular}
    \end{center}
\end{table}

Table \ref{tab:par} lists the parameter values for the model chosen to ensure that \eqref{eq:main} exhibits three distinct timescales where $V_1$ is fast, $w_1, V_2$ are slow and $w_2$ is superslow. In a more biologically realistic model for calcium and voltage interactions, $V_1$ might
represent membrane potential, while $V_2$ might represent intracellular calcium concentration with appropriate adjustments to parameter units and functional terms (see, e.g., \citep{WR2016, WR2017}). For the physiological description of functions in \eqref{eq:main} and \eqref{eq:term}, we refer readers to \cite{Nan2015} for details. 
\begin{figure*}[!htp]
\centering
    \begin{tabular}         
        {@{}p{0.48\linewidth}@{\quad}p{0.48\linewidth}@{}}
        \subfigimg[width=\linewidth]{\bf{\small{}}}{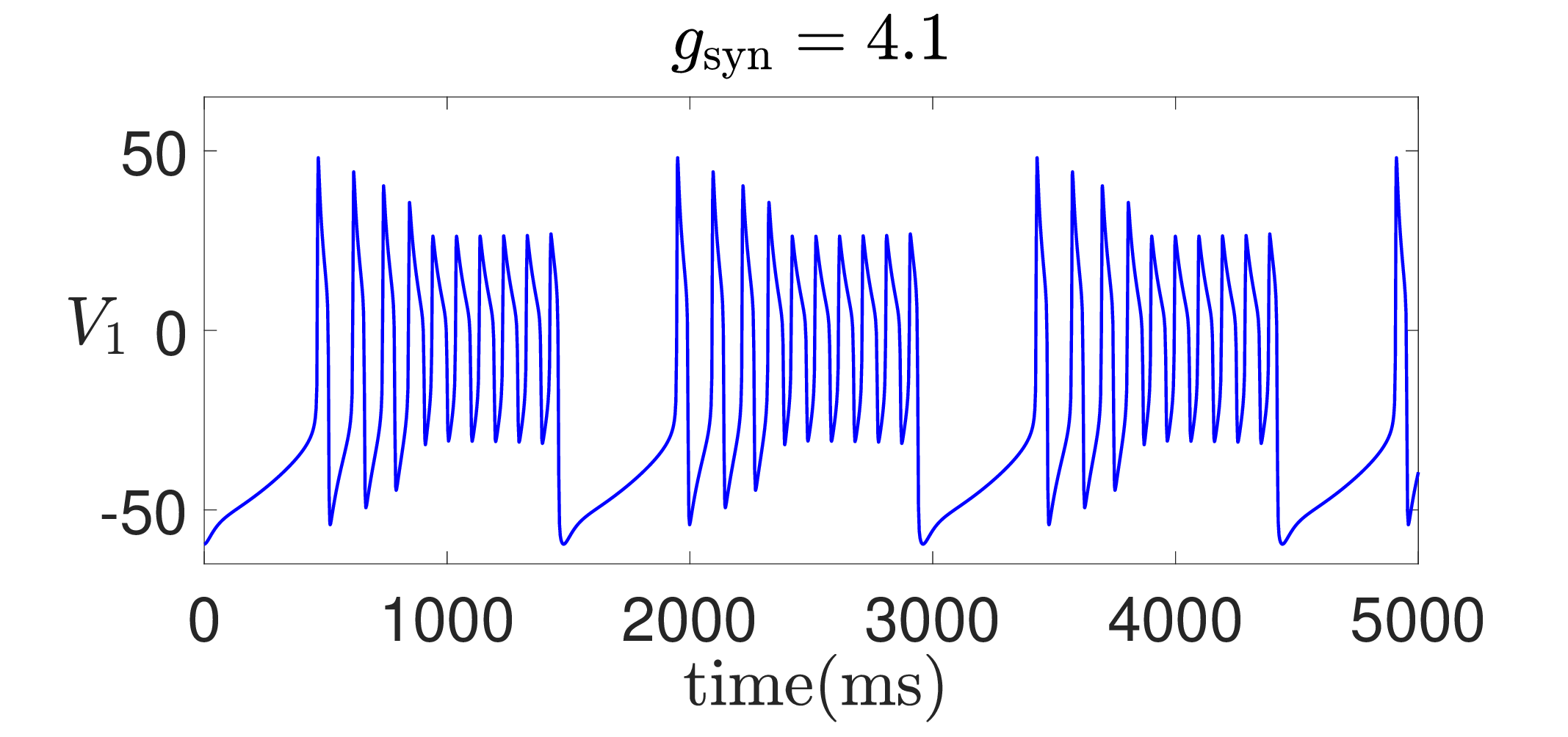}&
       \subfigimg[width=\linewidth]{\bf{\small{}}}{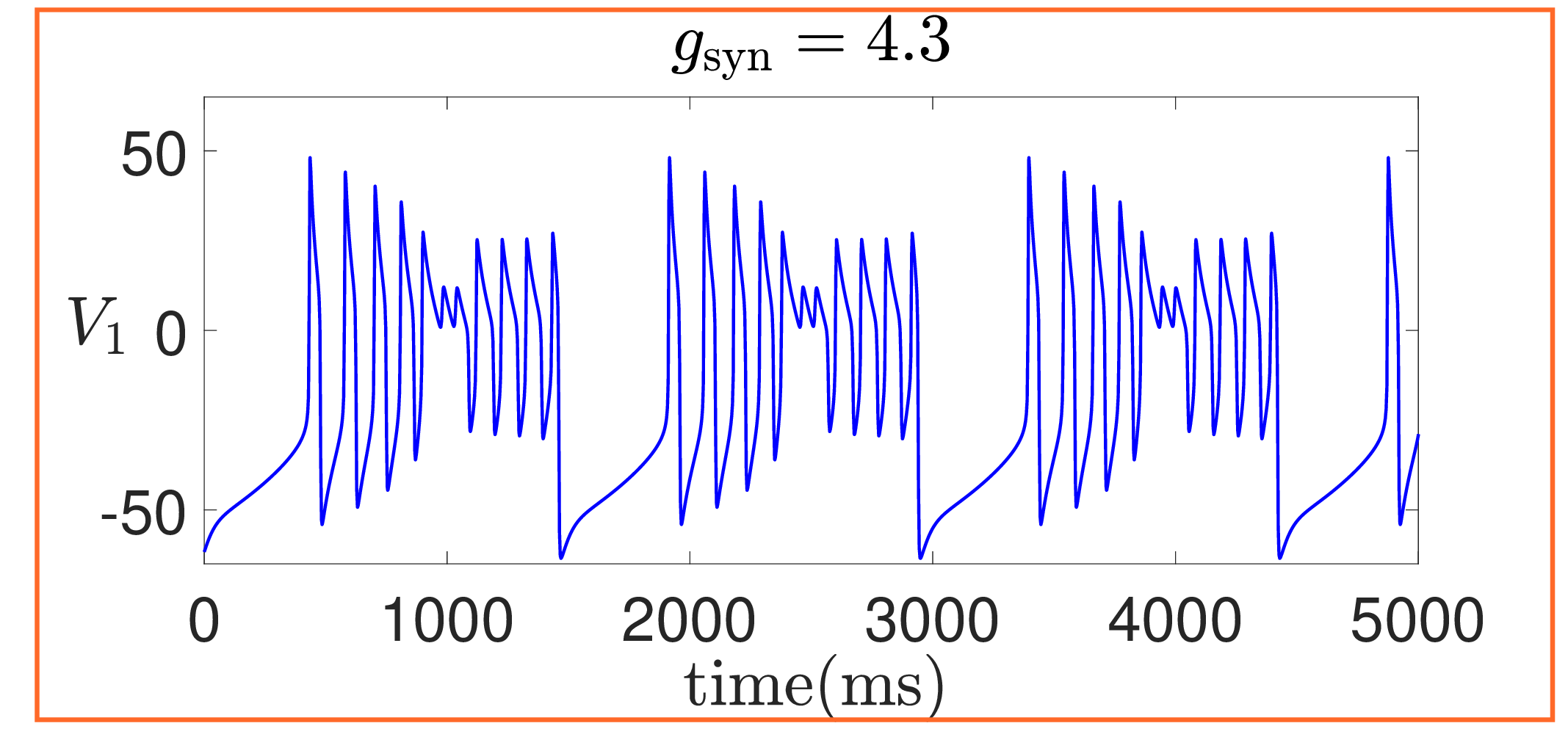}\\
        \subfigimg[width=\linewidth]{\bf{\small{}}}{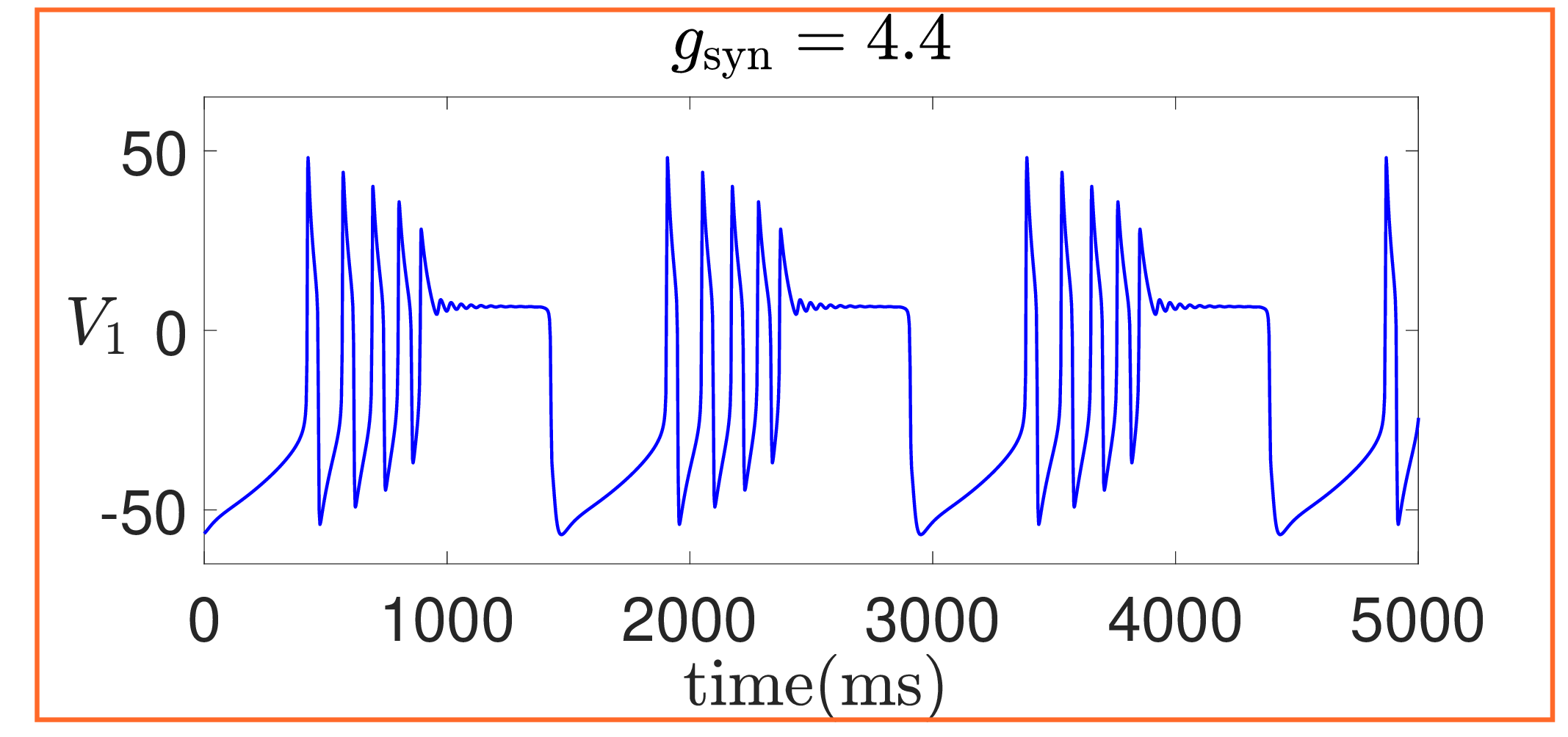}&
       \subfigimg[width=\linewidth]{\bf{\small{}}}{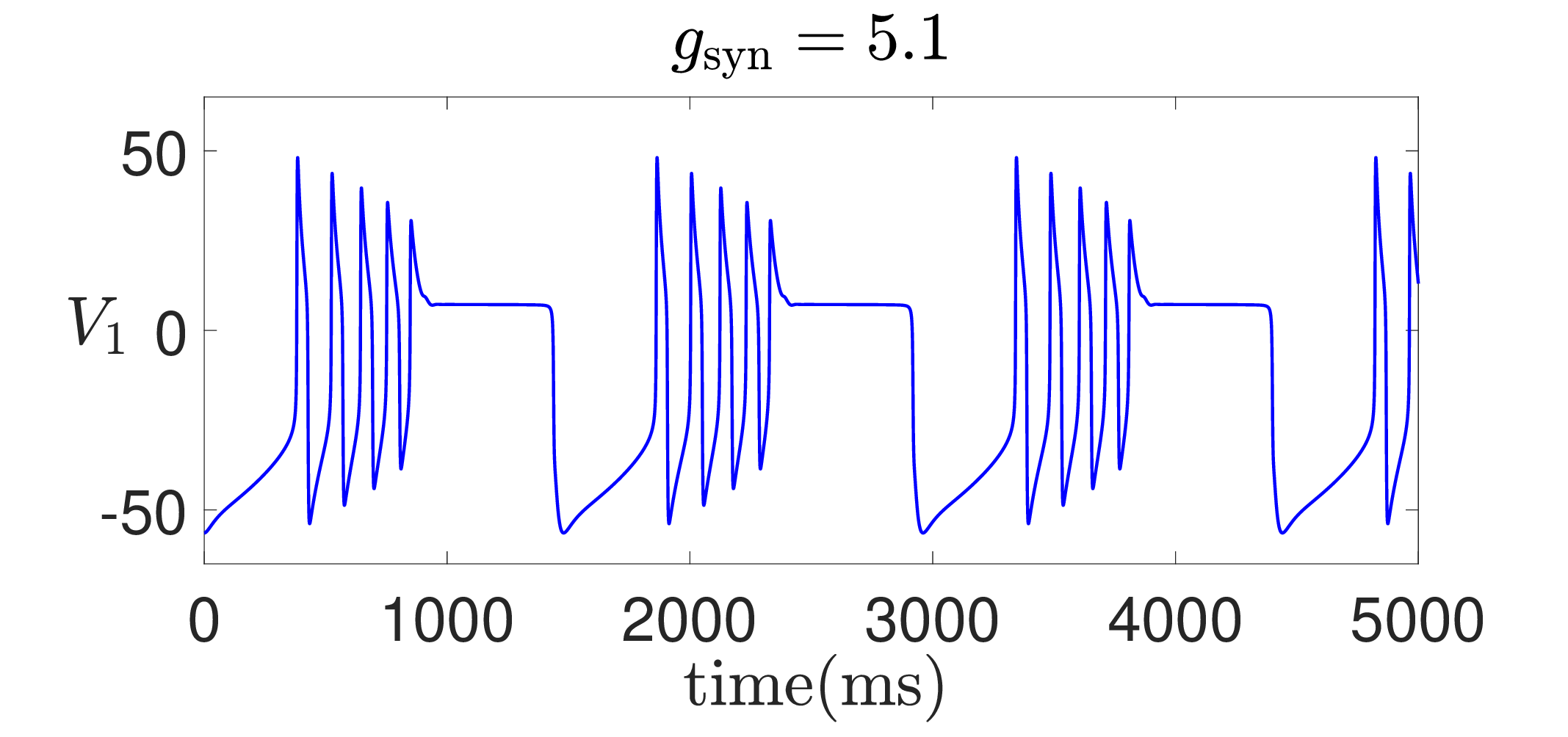}
   \end{tabular}
    \caption{Time traces of the model \eqref{eq:main} for \RED{different values of $g_{\rm syn}$. MMOs are observed for $\gsyn=4.3$ and $\gsyn=4.4$.}}
    \label{fig:simu}
\end{figure*}

In the absence of coupling $g_{\rm syn}=0$, $(V_1, w_1)$ is excitable with an attracting critical point at relatively low $V_1$ value, whereas $(V_2, w_2)$ is oscillatory with an attracting limit cycle independent of the value of $g_{syn}$ and the dynamics of $(V_1,w_1$). 
To analyze the three-timescale coupled Morris-Lecar neurons, Ref.~\onlinecite{Nan2015} extended two approaches previously developed in the context of GSPT for the analysis of two-timescale systems to the three-timescale setting and showed these two approaches complemented each other nicely. By varying $g_{\rm syn}$ in system \eqref{eq:main}, the authors identified various solution features that truly require three timescales, thus demonstrating the functional relevance of three timescales in the model. While system \eqref{eq:main} exhibits both the fast subsystem Hopf and folded nodes that can support MMOs, MMOs were not observed within the parameter regime examined by  Ref.~\onlinecite{Nan2015} \RED{(e.g., $\gsyn=4.1$ and $5.1$ in Figure \ref{fig:simu})}. 

\begin{figure}[!htp]
\begin{center}
   \includegraphics[width=\linewidth]{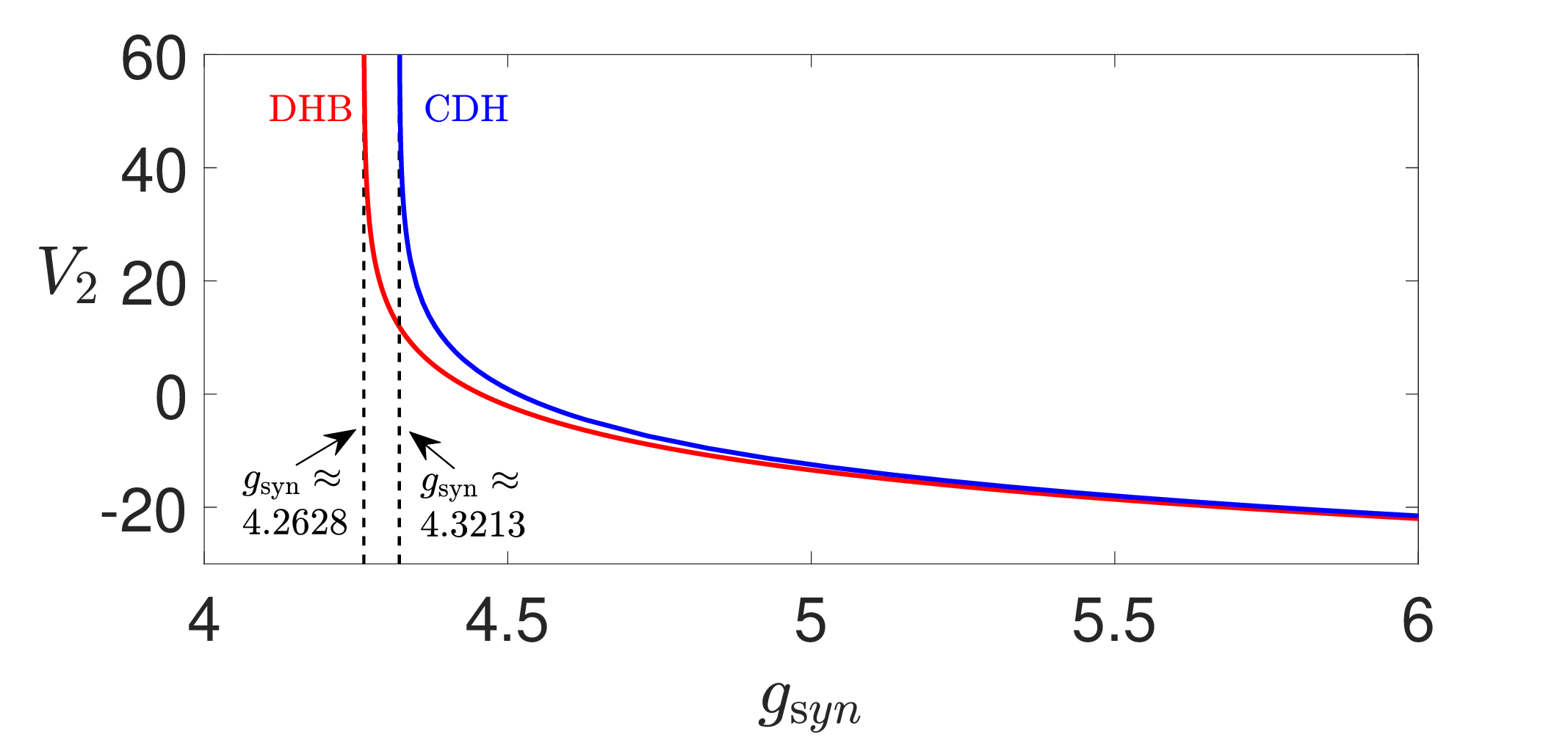}
\end{center}
    \caption{\RED{Bifurcation curves of DHB (red) and CDH singularities (blue) for \eqref{eq:main} with respect to $g_{\rm syn}$. Specifically, these curves represent the upper DHB and upper CDH, corresponding to larger $V_1$ and $V_2$ values. The lower CDH and DHB, associated with smaller $V_i$ values, are not presented here. 
    The two vertical asymptotes are given by $g_{\rm syn} \approx 4.2628$ and $g_{\rm syn} \approx 4.3213$. }}
    \label{fig:bd-gsyns}
\end{figure}

The goal of this work is the analysis of MMOs and their robustness in three-timescale systems by focusing on a coupled Morris-Lecar system \eqref{eq:main}. 
\RED{Based on our simulations, we have selected two MMO solutions on which to focus our analysis. Specifically, we consider $g_{\rm syn}=4.3$ and $g_{\rm syn}=4.4$ (with the unit of $\rm mS/cm^2$), as highlighted in Figure~\ref{fig:simu}. To provide further insight into our choice of $\gsyn$ values, we perform a bifurcation analysis to explore the effect of $\gsyn$ on a singularity called \textit{canard-delayed-Hopf} (CDH) that was firstly introduced by \citep{Letson2017} (see Figure \ref{fig:bd-gsyns}, blue). As noted above, this singularity plays a crucial role in organizing MMOs within the three-timescale setting. We show in Sections \ref{sec:gspt} that the full system \eqref{eq:main} may exhibit two CDH points - one at larger $V_i$ values, $i\in\{1,2\}$ (denoted as upper CDH) and the other at smaller $V_i$ (denoted as lower CDH, see Figure \ref{fig:weak-eigen-direc}). Similarly, \eqref{eq:main} may exhibit an upper DHB and a lower DHB. However, we demonstrate in Sections \ref{4p3} and \ref{4p4} that only the upper CDH or DHB can support MMOs. 

In Figure \ref{fig:bd-gsyns}, we plot the bifurcation curves of the upper CDH and the upper DHB in the $(\gsyn, V_2)$ plane, both of which approach vertical asymptotes as $\gsyn$ decreases. Our first choice of $\gsyn=4.3$ represents $\gsyn$ values between the two asymptotes, at which the upper DHB exists but there is no upper CDH. In contrast, $\gsyn=4.4$ lies to the right of the CDH asymptote and serves as a representative scenario where both the upper CDH and DHB exist. When $\gsyn<4.2628$ (e.g., $\gsyn=1$ and $4.1$ as considered in \citep{Nan2015}), the system \eqref{eq:main} does not produce MMOs. Although both CDH and DHB are present for $\gsyn=5.1$, the absence of MMOs as shown in Figure \ref{fig:simu} is due to the real nature of eigenvalues in a relevant subsystem during the transition from $V_1$ spikes to a $V_1$ plateau \citep{Nan2015}. 
}

\begin{figure*}[!htp]
\centering
    \begin{tabular}         
        {@{}p{0.48\linewidth}@{\quad}p{0.48\linewidth}@{}}
        \subfigimg[width=\linewidth]{\bf{\small{(A)}}}{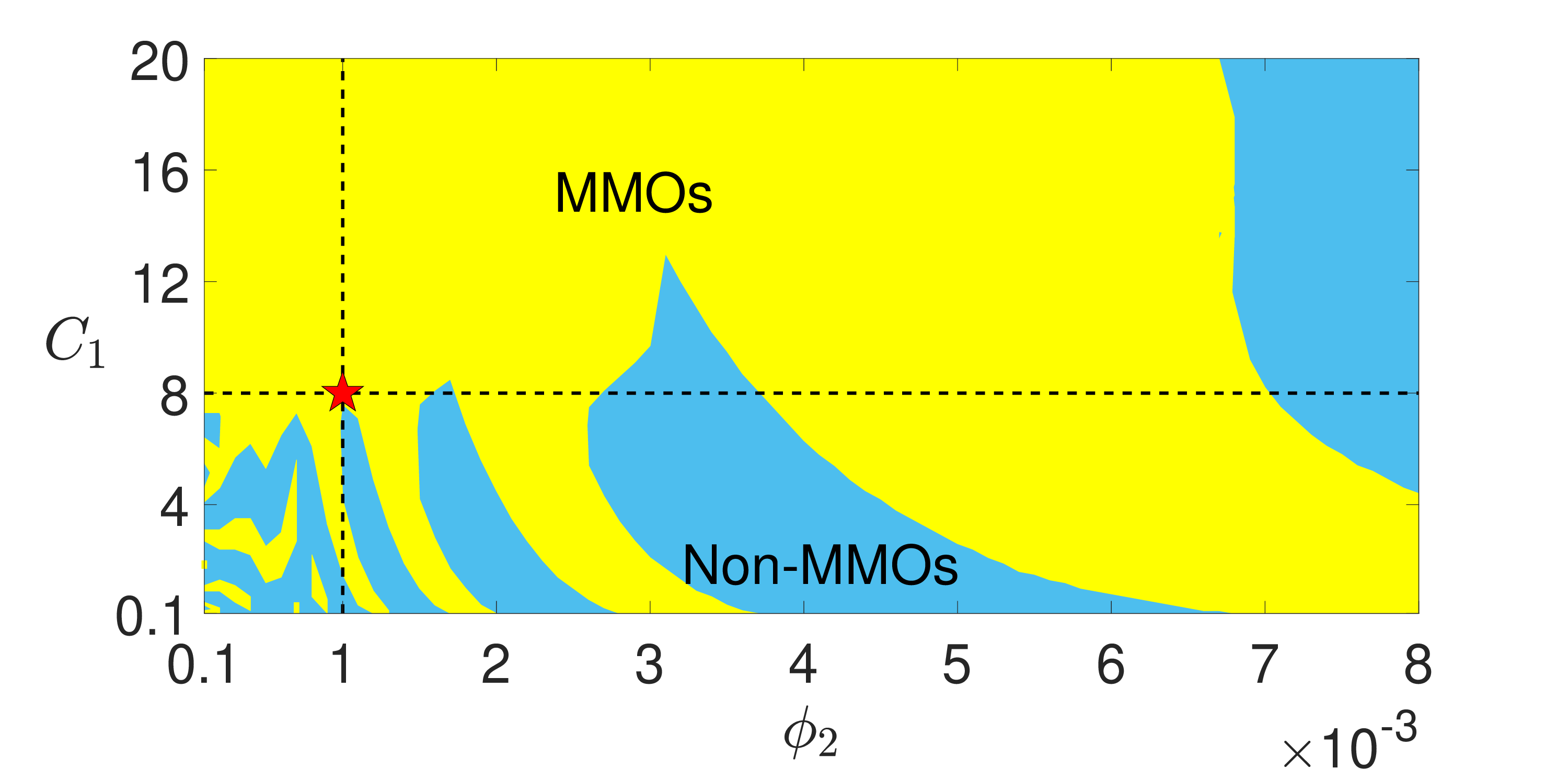}&
        \subfigimg[width=\linewidth]{\bf{\small{(B)}}}{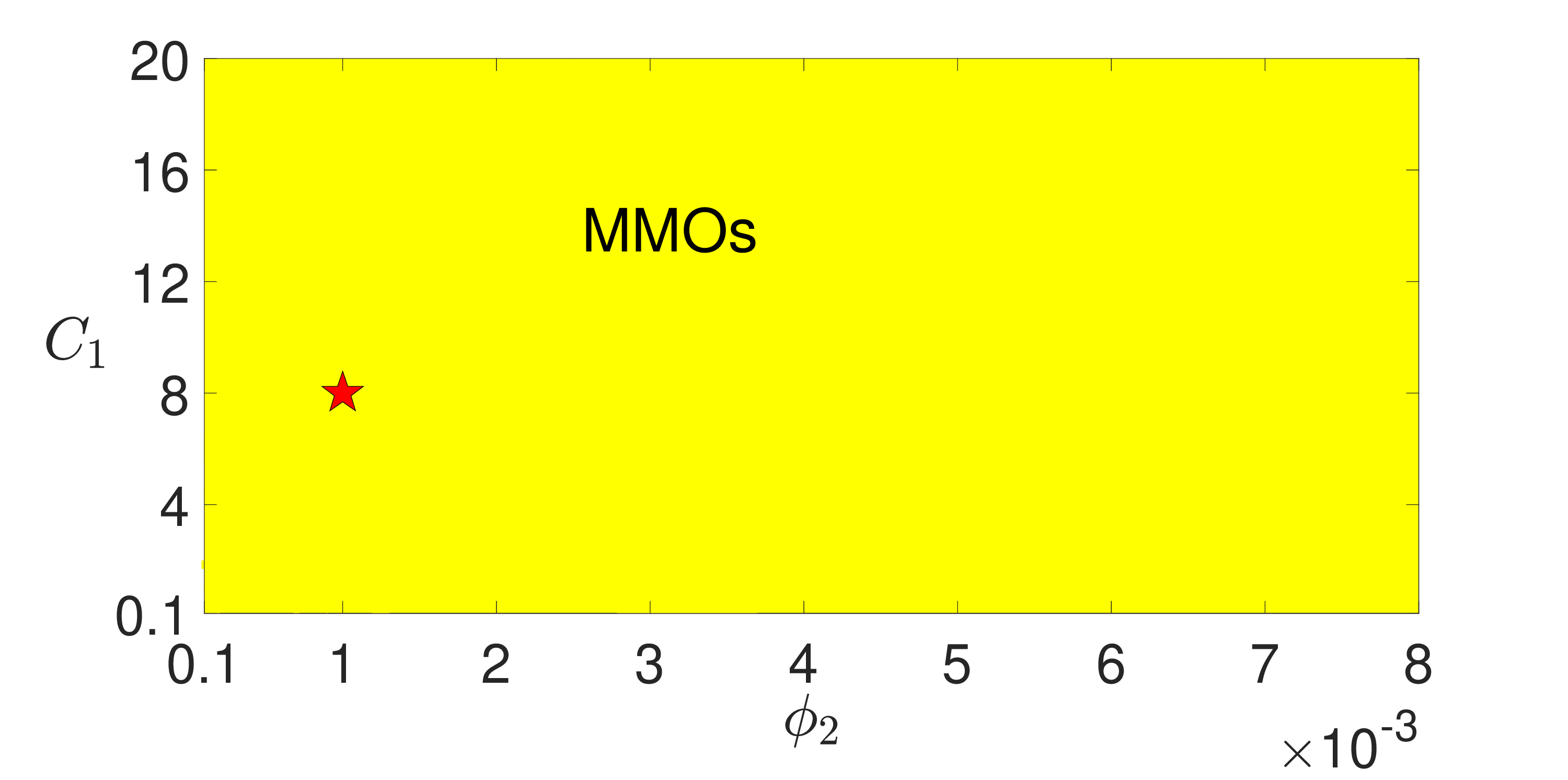}
    \end{tabular}
    \caption{Regions of MMOs (yellow) and non-MMOs (blue) of the full system \eqref{eq:main} in $(\phi_2, C_1)$-space for (A) $\gsyn=4.3$ and (B) $\gsyn=4.4$. Increasing $C_1$ slows down the fast variable $V_1$, whereas increasing $\phi_2$ speeds up the superslow variable $w_2$. \RED{The timescales of $w_1$ and $V_2$ remain unaffected.}
    The red star marks the default parameter values of $C_1$ and $\phi_2$ as given in Table \ref{tab:par}. (A) $g_{\rm syn}=4.3$. While MMOs are robust to increasing $C_1$ and decreasing $\phi_2$, decreasing $C_1$ or increasing $\phi_2$ leads to multiple transitions between MMOs and non-MMOs (crossings between the dashed lines with the yellow/blue boundaries). (B) $g_{\rm syn}=4.4$. MMOs are robust to changes of both $C_1$ and $\phi_2$ over the ranges of $0.1\leq C_1\leq 20$ and $0.1\leq \phi_2\leq 8e-3$. Note that the MMOs at $\gsyn=4.4$ will eventually vanish for $C_1$ and $\phi_2$ large enough at which there is no more timescale separation (data not shown).}
    \label{fig:c1-phi2}
\end{figure*}

Our analysis suggests that the two MMO solutions at \RED{$\gsyn=4.3$ and $\gsyn=4.4$} arise from distinct mechanisms, resulting in remarkably different sensitivities to variations in timescales (i.e., varying $C_1$ and $\phi_2$), as illustrated in Figure~\ref{fig:c1-phi2}. Utilizing the extended GSPT \citep{Fenichel1979,Nan2015}, \RED{we show that there is no interaction of different MMO mechanisms due to the lack of a nearby CDH singularity when $\gsyn=4.3$ (see Figure \ref{fig:bd-gsyns}). Instead, the MMOs at $\gsyn=4.3$ solely depend on the delayed Hopf mechanism and are sensitive to variations in timescales (Figure \ref{fig:c1-phi2}A). 
In contrast, there exists a CDH in the middle of the SAOs when $\gsyn=4.4$ as discussed above. We demonstrate that this CDH allows the fast subsystem Hopf and a canard point to coexist and interact to co-modulate properties of the local oscillatory behavior, resulting in MMOs with significantly stronger robustness than $\gsyn=4.3$ (Figure \ref{fig:c1-phi2}B). 
In summary, our findings reveal that MMOs near a CDH exhibit stronger robustness compared to those governed by a single mechanism, and that the CDH singularity is a determining factor in whether the two distinct MMO mechanisms can interact or not. However, it is essential to note that not all CDH singularities can support local MMOs. Specifically, we demonstrate that  the lower CDH does not support MMOs. }

 

Our work is novel in two main aspects. First, to the best of our knowledge, our study is the first to investigate the geometric conditions that lead to robust occurrences of MMOs in three-timescale systems. It is worth noting that while Ref.~\onlinecite{Kak2022} also considered the robustness of MMOs in a three-timescale system, their focus was specifically on MMOs with double epochs of SAOs. 
Second, we discovered that the CDH singularities do not always enable the two MMO mechanisms to interact and produce MMOs. This is different from past studies \citep{Vo2013,Letson2017} where the CDH always leads to occurrence of MMOs.  From analyzing system \eqref{eq:main}, we found that CDH singularities that lie close to the actual fold of the superslow manifold (defined later by \eqref{superslow}) do not support MMOs regardless of perturbation sizes $\varepsilon$ and $\delta$. 

As the first step of our timescale decomposition approach, we perform a dimensional analysis of \eqref{eq:main} to reveal the important timescales. This transforms \eqref{eq:main} to the following three-timescale problem
\begin{equation}\label{eq:slow}
    \begin{array}{rcl}
        \varepsilon\frac{dV_1}{dt_s} &=& f_1(V_1,w_1,V_2), \vspace{0.1in} \\
        \frac{dw_1}{dt_s}&=&g_1(V_1,w_1), \vspace{0.1in} \\
        \frac{dV_2}{dt_s}& =&f_2(V_2,w_2), \vspace{0.1in}\\
        \frac{dw_2}{dt_s}&=&\delta g_2(V_2,w_2),
    \end{array}
\end{equation}
where $\varepsilon=0.1$,  $\delta=0.053$, $t_s$ is the slow dimensionless time variable, $f_1$, $f_2$, $g_1$ and $g_2$ are $O(1)$ functions specified in \eqref{eq:slow-appendix} in the Appendix \ref{ap:nondim} which include details of the nondimensionalization procedure. For simplicity, we did not rescale $V_1$ and $V_2$ in \eqref{eq:slow} as the scalings of voltage have no influence on the timescales. 

We call system \eqref{eq:slow} that is described over the \textit{slow timescale} the \textit{slow system} in which $V_1$ evolves on a timescale of $O(\varepsilon^{-1})$, $(w_1,V_2)$ on $O(1)$ and $w_2$ on $O(\delta)$. 
Introducing a superslow time $t_{ss}=\delta t_s$ yields an equivalent description of dynamics:
\begin{equation}\label{eq:ss}
    \begin{array}{rcl}
        \varepsilon\delta\frac{dV_1}{dt_{ss}}&=&f_1(V_1,\,V_2,\,w_1), \vspace{0.1in} \\
        \delta\frac{dw_1}{dt_{ss}}&=&g_1(V_1,\,w_1), \vspace{0.1in} \\
        \delta\frac{dV_2}{dt_{ss}}&=&f_2(V_2,\,w_2) ,\vspace{0.1in} \\
        \frac{dw_2}{dt_{ss}}&=&g_2(V_2,\,w_2),
    \end{array}
\end{equation}
which evolves on the {\em superslow timescale} and is called the \textit{superslow system}.
Alternatively, defining a fast time $t_f = t_s/ \varepsilon$,
we obtain the following \textit{fast system}:
\begin{equation}\label{eq:fast}
    \begin{array}{rcl}
        \frac{dV_1}{dt_f}&=&f_1(V_1,\,V_2,\,w_1), \vspace{0.1in} \\
        \frac{dw_1}{dt_f}&=&\varepsilon g_1(V_1,\,w_1), \vspace{0.1in} \\
        \frac{dV_2}{dt_f}&=&\varepsilon f_2(V_2,\,w_2), \vspace{0.1in} \\
        \frac{dw_2}{dt_f}&=&\varepsilon\delta g_2(V_2,\,w_2),
    \end{array}
\end{equation}
which evolves on the {\em fast timescale}.

The paper is organized as follows. In Section \ref{sec:gspt}, we perform a geometric singular perturbation analysis of the 3-timescale problem \eqref{eq:main} by treating $\varepsilon$ as the principal perturbation parameter while keeping $\delta$ fixed, treating $\delta$ as the principal perturbation parameter while keeping $\varepsilon$ fixed, and by treating $\varepsilon$ and $\delta$ as two independent perturbation parameters. We review both mechanisms for MMOs and discuss their interaction at the double singular limit $(\varepsilon,\delta)\to (0,0)$. Notation, subsystems, \RED{construction of singular orbits at various singular limits} and other preliminaries relating to the method of GSPT are all presented in Section \ref{sec:gspt}. 
\RED{In Section \ref{4p3}, we investigate the mechanism and sensitivity of MMOs when $g_{\rm syn}=4.3$ to variations in perturbation sizes $\varepsilon$ or $\delta$ (i.e. varying $C_1$ and $\phi_2$ in \eqref{eq:main}). We explain the transitions between MMO and non-MMO dynamics as we vary one perturbation parameter while keeping the other fixed at its default value, as illustrated by the two lines in Figure \ref{fig:c1-phi2}A. Our analysis indicates that MMOs persist for $\delta \leq \mathcal{O}(\varepsilon)$. Increasing $\delta$ via increasing $\phi_2$ or decreasing $\varepsilon$ via decreasing $C_1$ to a degree where $\delta > \mathcal{O}(\varepsilon)$ leads to multiple MMOs/non-MMOs transitions. MMOs are completely lost when $\delta \geq \mathcal{O}(\varepsilon^{\frac{1}{3}})$.
While Figure \ref{fig:c1-phi2}A also shows transitions occurring when both $\varepsilon$ and $\delta$ are relatively large, the analysis of these transitions is beyond the standard GSPT and falls outside the scope of this paper.}
In Section \ref{4p4}, we uncover the dynamic mechanism underlying MMOs from \eqref{eq:main} when $g_{\rm syn}=4.4$. \RED{In this case, the existence of a CDH in the middle of the SAOs enables two different mechanisms to co-modulate properties of MMOs}. We explain why MMOs organized by a CDH singularity as seen in the case of $\gsyn=4.4$ exhibit remarkable robustness against variations in timescales (see Figure \ref{fig:c1-phi2}B). \RED{Our analysis demonstrate that when $\delta=\mathcal{O}(\varepsilon)$, MMOs exhibit both canard and DHB features. Upon tuning $\delta\geq O(\sqrt{\varepsilon})$, DHB-like features disappear and the canard mechanism dominates.}
Finally, we conclude in Section \ref{sec:discussion} with a discussion.

\section{Geometric Singular Perturbation Analysis}\label{sec:gspt}

In this section, we apply the extended geometric singularity perturbation analysis \citep{Nan2015,Fenichel1979} to the three-timesale coupled Morris-Lecar system \eqref{eq:slow} by treating $\varepsilon$ as the only singular perturbation parameter \citep{SW2001,Wechselberger2005}, treating $\delta$ as the only singular perturbation parameter \citep{Baer1989, Neishtadt1987, Neishtadt1988, Hayes2016}, and finally treating $\varepsilon$ and $\delta$ as two independent singular perturbation parameters \citep{Vo2013,Nan2015,Letson2017}.

Although the detailed GSPT analysis and derivation of subsystems have been previously presented in \citep{Nan2015}, we provide a brief overview in this paper for the sake of completeness. However, the focus of our current work is on the investigation of MMOs, which is distinct from the emphasis of \citep{Nan2015}. Specifically, we concentrate on reviewing and discussing the canard mechanism in subsection \ref{subsec:canard}, delayed Hopf bifurcation in subsection \ref{sec:delayed-HB}, and their interactions in subsection \ref{subsec:interaction}.

\subsection{Singular Limits}\label{subsec:singular-limits} 

\subsubsection{$\varepsilon\to 0$ singular limit.} Fixing $\delta > 0$ and taking $\varepsilon\to0$ in the fast system (\ref{eq:fast}) yields the one-dimensional (1D) {\em fast layer problem}, a system that describes the dynamics of the fast variable, $V_1$, for fixed values of the other variables, 
\begin{equation}\label{eq:fastlayer}
    \begin{array}{rcl}
        \frac{dV_1}{dt_f}=f_1(V_1,\,w_1,\,V_2).
    \end{array}
\end{equation} 
The set of equilibrium points of the fast layer problem is called the \textit{critical manifold} and is denoted as $\ms$:
\begin{equation}\label{eq:ms}
    \ms:=\{(V_1,\,w_1,\,V_2,\,w_2):\, f_1(V_1,\, w_1, \,V_2)=0\}\,.
\end{equation}
Although $\ms$ is a three-dimensional (3D) manifold in $\mathbb{R}^4$ space, it does not  depend on $w_2$.
We can solve $f_1(V_1,\, w_1, \, V_2)=0$ for $w_1$ as a function of $V_1$ and $V_2$ and can therefore represent $\ms$ as
\begin{equation}\label{eq:critmfld}
    w_1=F_1(V_1,\,V_2)
\end{equation}
for a function $F_1$.
It is well known that, for sufficiently small $\varepsilon>0$, normally hyperbolic parts of $\ms$ each perturb to a locally invariant manifold called a {\em slow manifold}, on which $w_1$ is given by an $O(\varepsilon)$-perturbation of $F_1$ \citep{Fenichel1979}; we simply use $\ms$ as a convenient numerical approximation of these slow manifolds.

$\ms$ is a 3D folded manifold with two-dimensional (2D) fold surface, $\ls$, given by 
\begin{equation}\label{eq:L-F1}
    \ls:=\{(V_1,\,w_1,\,V_2,\,w_2)\in \ms: F_{1V_1}=0\},
\end{equation}
or equivalently
\begin{equation}\label{eq:L-f1}
    \ls:=\left\{(V_1,\,w_1,\,V_2,\,w_2) \in \ms :\, 
    \ f_{1V_1}=0 \right\},
\end{equation}
\RED{where $F_{1V_1}$ and $f_{1V_1}$ denote the partial derivatives of $F_1$ and $f_1$ with respect to $V_1$}. The fold surface divides the critical manifold $M_S$ into attracting upper $(M_S^U)$ and lower $(M_S^L)$ branches where $F_{1V_1}<0$ and repelling middle branch $(M_S^M)$ where $ F_{1V_1}>0$.

Taking the same limit, i.e.,  $\varepsilon\to0$ with $\delta > 0$, in the slow system (\ref{eq:slow}) 
yields the 3D \textit{slow reduced problem}, a system that describes the dynamics of $w_1, V_2, w_2$ along $\ms$, 
\begin{equation}\label{eq:slowreduced}
    \begin{array}{rcl}
        \frac{dw_1}{dt_s}&=&g_1(V_1,\,w_1), \vspace{0.1in} \\
        \frac{dV_2}{dt_s}&=&f_2(V_2,\,w_2), \vspace{0.1in}\\
        \frac{dw_2}{dt_s}&=&\delta g_2(V_2,\,w_2),
    \end{array}
\end{equation}
where $f_1=0$.

\subsubsection{$\delta\to 0$ singular limit.} Alternatively, fixing $\varepsilon>0$ and taking $\delta\to0$ in the slow system (\ref{eq:slow}) yields the 3D {\em slow layer problem} in the form
\begin{equation}\label{eq:slowlayer}
    \begin{array}{rcl}
        \varepsilon\frac{dV_1}{dt_s}&=&f_1(V_1,\,w_1,\,V_2), \vspace{0.1in}\\
        \frac{dw_1}{dt_s}&=&g_1(V_1,\,w_1), \vspace{0.1in}\\
        \frac{dV_2}{dt_s}&=&f_2(V_2,\,w_2),
    \end{array}
\end{equation}
where the superslow variable $w_2$ is a parameter. 

The set of equilibrium points of the slow layer problem \eqref{eq:slowlayer} is defined to be the {\em superslow manifold} and is denoted as $\mss$
\begin{equation}\label{superslow}
    \begin{array}{rcl}
    \mss&:=&\{(V_1,\,w_1,\,V_2,\,w_2):\, \vspace{0.1in} \\ &&f_1(V_1,\,w_1,\,V_2)=
    g_1(V_1,\,w_1)=f_2(V_2,\,w_2)=0\}.
    \end{array}
\end{equation}
$\mss$ is a 1D subset of 
$\ms$. Similarly to $\ms$, the normally hyperbolic parts of $\mss$ perturb to nearly locally invariant manifolds for $\delta$ sufficiently small. Later in subsection \ref{sec:delayed-HB}, we will discuss the bifurcations of the slow layer problem \eqref{eq:slowlayer}, i.e., nonhyperbolic regions on $\mss$ where Fenichel's theory (GSPT) breaks down.

Taking the same singular limit in the superslow system \eqref{eq:ss} leads to the 1D \textit{superslow reduced problem}
\begin{equation}\label{eq:ssreduced}
    \begin{array}{rcl}
        \frac{dw_2}{dt_{ss}}&=&g_2(V_2,\,w_2),
    \end{array}
\end{equation}
where $f_1=g_1=f_2=0$. 
The superslow motions of trajectories of \eqref{eq:ssreduced} are slaved to $\mss$ until nonhyperbolic points are reached.

\subsubsection{$\varepsilon\to 0, \,\delta\to 0$ double singular limits.} Both the slow reduced problem \eqref{eq:slowreduced} and the slow layer problem \eqref{eq:slowlayer} still include two distinct timescales. Further taking the limit $\delta \to 0$ in (\ref{eq:slowreduced}) or taking the limit $\varepsilon\to 0$ in \eqref{eq:slowlayer}
yields the same {\em slow reduced layer problem}, 
\begin{equation}\label{eq:slowreducedlayer}
    \begin{array}{rcl}
        0&=&f_1(V_1,\,V_2,\,w_1), \vspace{0.1in}\\
        \frac{dw_1}{dt_s}&=&g_1(V_1,\,w_1), \vspace{0.1in}\\
        \frac{dV_2}{dt_s}&=&f_2(V_2,\,w_2),
    \end{array}
\end{equation}
which describes the slow motion along $\ms$ and the superslow variable $w_2$ is fixed as a constant. 

It follows that the double singular limits lead to three subsystems: the fast layer problem \eqref{eq:fastlayer}, the slow reduced layer problem \eqref{eq:slowreducedlayer} and the superslow reduced problem \eqref{eq:ssreduced}. In addition to the naturally expected fast/slow transitions and slow/superslow transitions, transitions directly from superslow to fast dynamics and from superslow to fast-slow relaxation oscillations have also been observed in \citep{Nan2015}. 

\subsection{Slow Reduced Problem and Canard Dynamics}\label{subsec:canard}

To investigate canard dynamics, 
we project the slow reduced problem \eqref{eq:slowreduced} onto $(V_1,\,V_2,\,w_2)$ to obtain a complete description of the dynamics along $\ms$. To this end, we differentiate the graph representation of $\ms$ given by $w_1=F_1(V_1,\,V_2)$ to obtain 
\begin{equation} \label{proj_slowreduced}
    \begin{array}{rcl}
         F_{1V_1}\frac{dV_1}{dt_{s}} & = & g_1(V_1,\,w_1)- F_{1V_2}f_2(V_2,\,w_2), \vspace{.1in}\\
        \frac{dV_2}{dt_{s}} & = & f_2(V_2,\,w_2),\vspace{.1in}\\
        \frac{dw_2}{dt_{s}} & = &  \delta g_2(V_2,\,w_2),
    \end{array}
\end{equation}
\RED{where $F_{1V_2} := \frac{\partial F_1}{\partial V_2}$.} 
Note that the reduced system \eqref{proj_slowreduced} is singular at the fold surfaces $\ls$ \eqref{eq:L-F1}. Nonetheless, this singular term can be removed by a time rescaling $t_{s}=- F_{1V_1} t_{d}$ and we obtain the following desingularized system 
\begin{equation} \label{desingu_slowreduced}
    \begin{array}{rcl}
         \frac{dV_1}{dt_{d}} & = & F_{1V_2}f_2(V_2,\,w_2)-g_1(V_1,\,w_1):=F(V_1,\,w_1,\,V_2,\,w_2), \vspace{.1in}\\
        \frac{dV_2}{t_{d}} & = &- F_{1V_1}f_2(V_2,\,w_2):=G(V_1,\,V_2,\,w_2),\vspace{.1in}\\
        \frac{dw_2}{dt_{d}} & = & -\delta F_{1V_1}g_2(V_2,\,w_2):=H(V_1,\,V_2,\,w_2,\,\delta).
    \end{array}
\end{equation}
We observe that the desingularized system (\ref{desingu_slowreduced}) is equivalent to (\ref{proj_slowreduced}) on the attracting branch, i.e, for $F_{1V_1}<0$, but has the opposite orientation on the repelling branch, i.e, for $F_{1V_1}>0$.

The desingularized system \eqref{desingu_slowreduced} has two kinds of singularities: ordinary and folded singularities. The ordinary singularities are the true equilibria of the full system (\ref{eq:main}), which is defined by  
\begin{equation} \label{ordinary}
    \begin{array}{rcl}
        E&:=&\{ (V_1,\,w_1,\,V_2,\,w_2)\in M_S: \vspace{0.1in} \\
        &&g_1(V_1,\,w_1)=f_2(V_2,\,w_2)=g_2(V_2,\,w_2)=0 \}.
    \end{array}
\end{equation}
For the chosen parameter set in Table \ref{tab:par}, $E$ always lies on the repelling branch of $\ms$ and hence is unstable.
In contrast to the ordinary singularities, the folded singularities are not equilibria of the full system. They lie on one-dimensional curves along the fold surface $\ls$ defined by
\begin{equation} \label{folded}
    \RED{\mathcal{M}}:= \{ (V_1,\,w_1,\,V_2,\,w_2)\in \ls: F(V_1,\,w_1,\,V_2,\,w_2)=0 \}.
\end{equation}
Folded singularities are special points that allow trajectories of \eqref{proj_slowreduced} to cross the fold $\ls$ with nonzero speed. Such solutions are called singular canards \citep{SW2001}. Note that when projecting to $(V_1,\,w_1,\,V_2)$-space, the condition $F=0$ is redundant and the fold surfaces $\ls$ become curves that overlap with the folded singularity curves $\RED{\mathcal{M}}$.

Since there is a curve of folded singularities, the Jacobian of \eqref{desingu_slowreduced} evaluated along $\RED{\mathcal{M}}$ (denoted as $J_D$, see Appendix \ref{app:cond-FSN}) always has a zero eigenvalue and the eigenvector corresponding to this zero eigenvalue is tangent to \RED{the curve of folded singularities}. Generically, the other two eigenvalues ($\lambda_w,\lambda_s$) where $|\lambda_w|<|\lambda_s|$ have nonzero real part and are used to classify the folded singularities. Folded singularities with two real eigenvalues with the same sign (resp., with opposite signs) are called \textit{folded nodes} (resp., \textit{folded saddles}). Those with complex eigenvalues are called \textit{folded foci}, which does not produce canard dynamics. 

\begin{figure*}[!t]
\begin{center}
\begin{tabular}{@{}p{0.48\linewidth}@{\quad}p{0.48\linewidth}@{}}
\subfigimg[width=\linewidth]{\bf{\small{(A)}}}{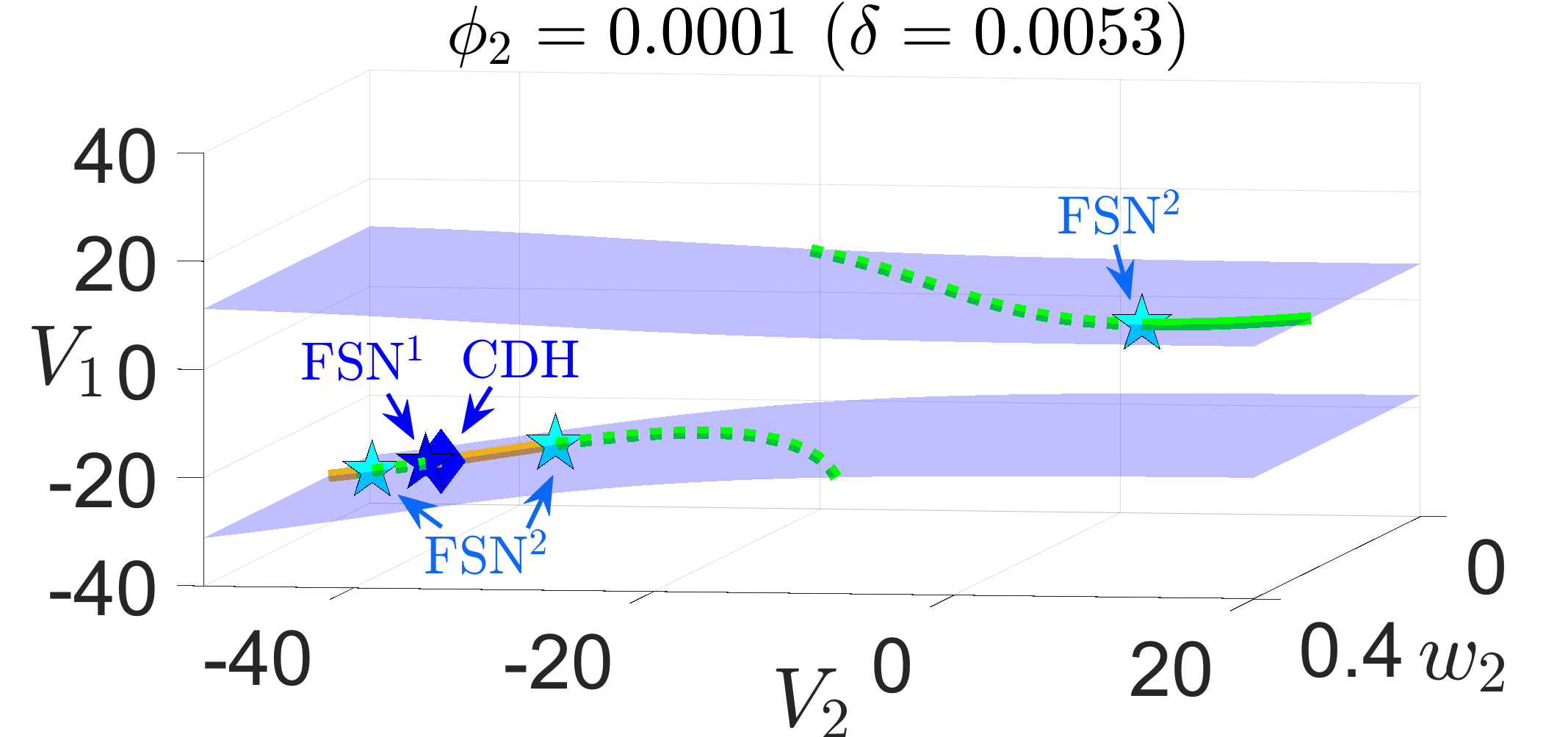}
&\subfigimg[width=\linewidth]{\bf{\small{(B)}}}{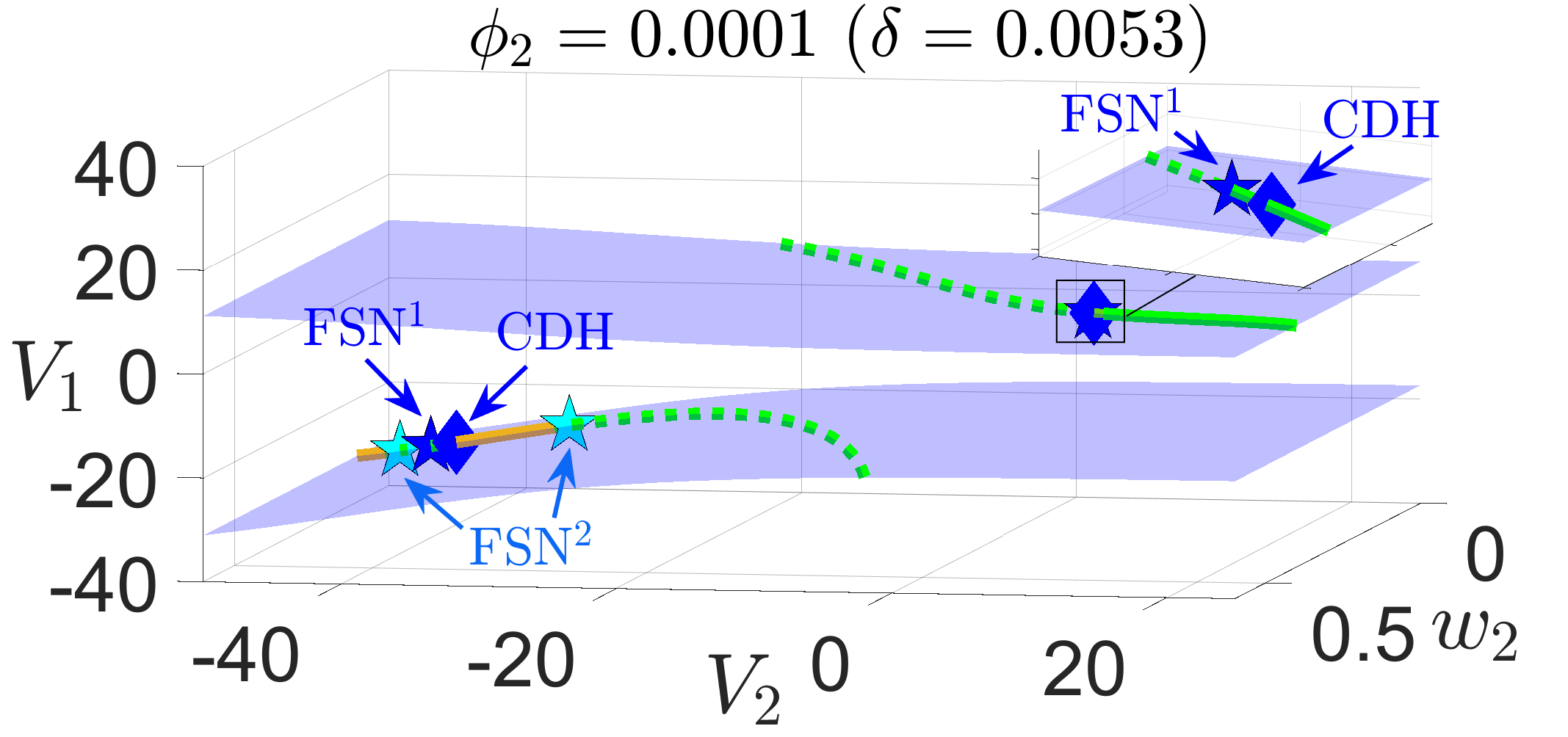}\\
\subfigimg[width=\linewidth]{\bf{\small{(C)}}}{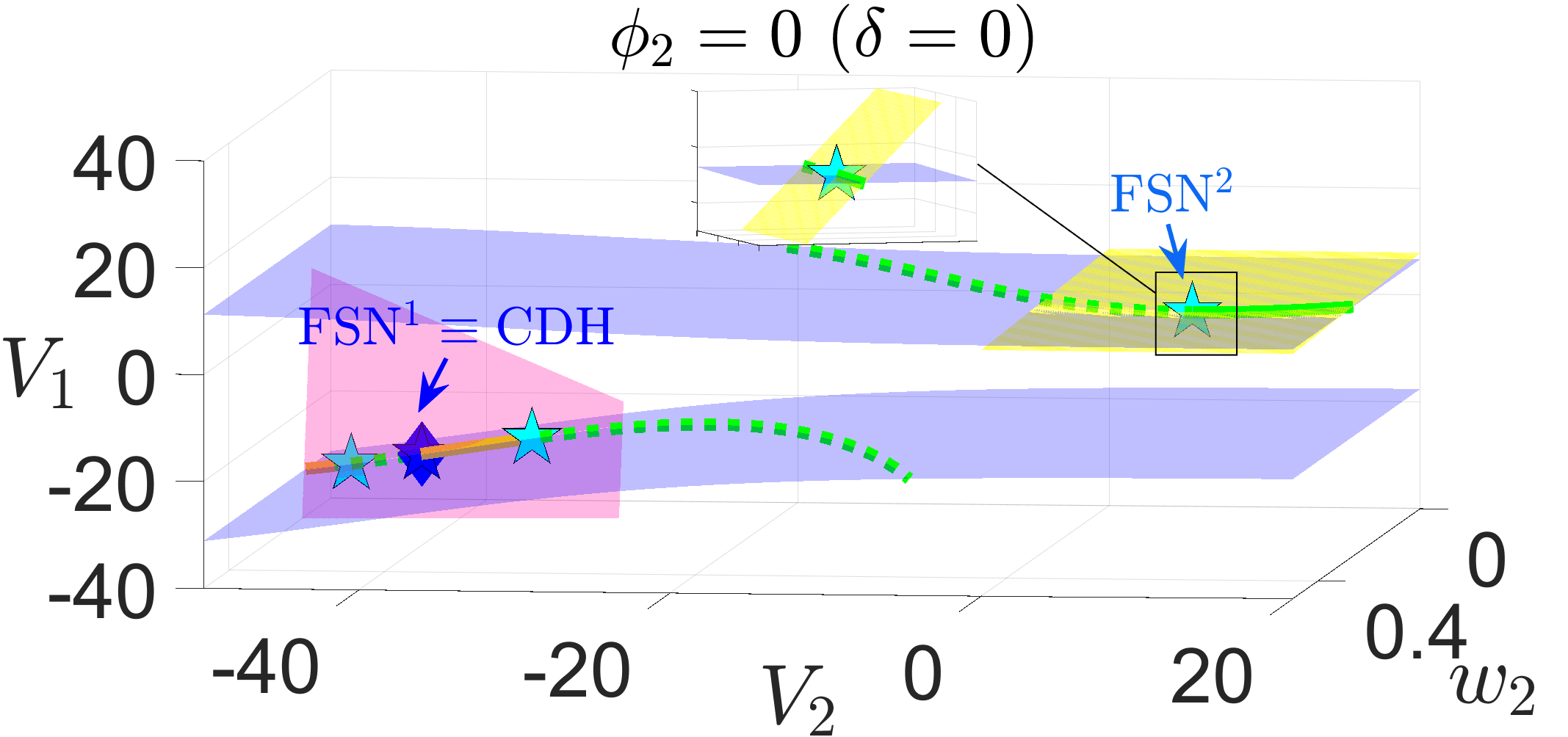}
&\subfigimg[width=\linewidth]{\bf{\small{(D)}}}{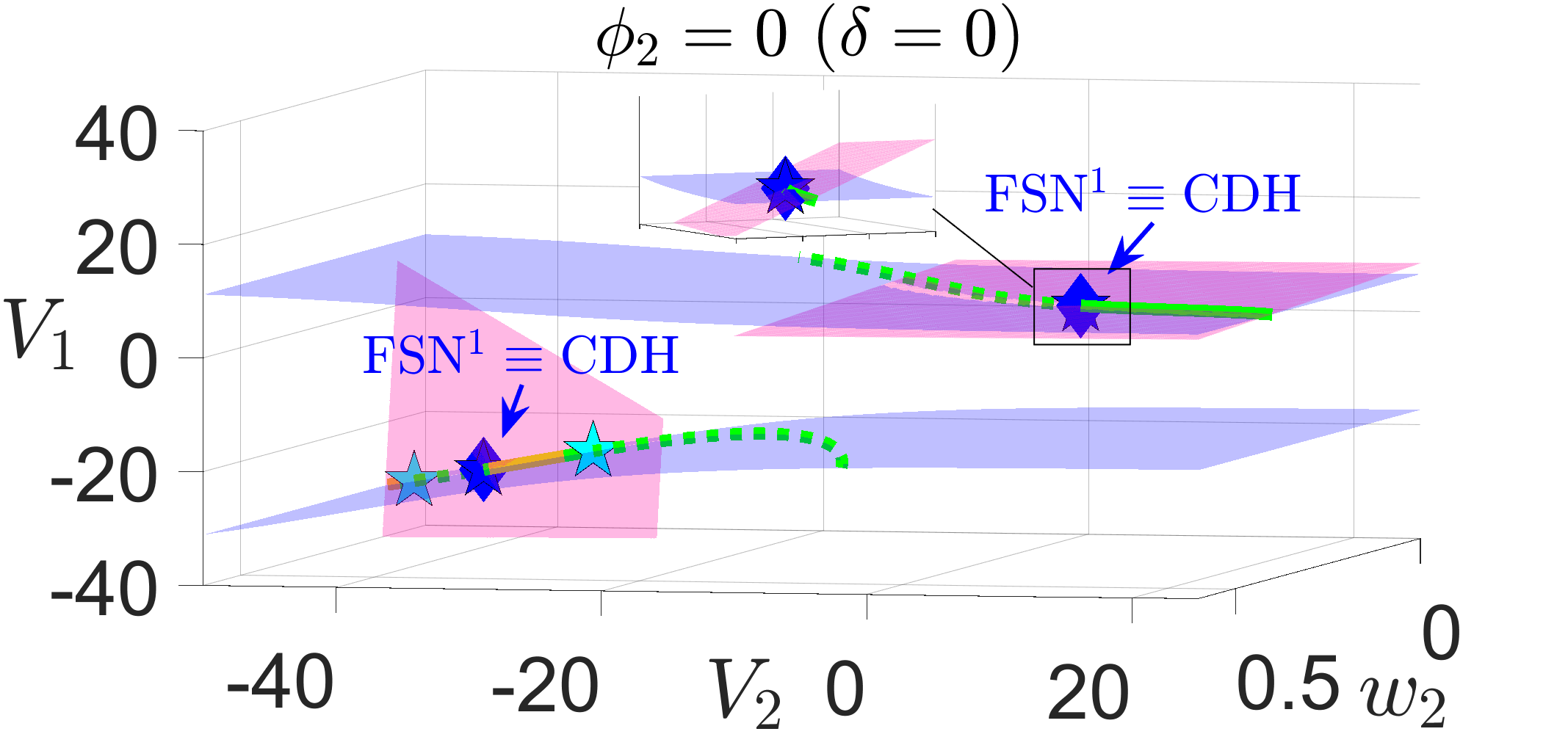}
\end{tabular}
\end{center}
\caption{Projections to $(V_1,V_2,w_2)$-space of the critical manifold fold surfaces $\ls$ (blue surface) for (A, C) $\rm g_{syn}=4.3$ and (B, D) $\rm g_{syn}=4.4$. Also shown are the curves of folded singularities $\RED{\mathcal{M}}$ including folded node (solid green), folded saddle (dashed green), and two types of folded saddle-nodes $\fsn$ (blue star: $\fsn^1$; cyan star: $\fsn^2$). The yellow curve consists of mostly folded foci points and small segments of other singularities (e.g., folded node, folded saddle) that are barely visible and hence are not displayed here. In the top two panels (A, B) when $\delta\neq 0$, an $\fsn^1$ point (blue star) is $O(\delta)$ close to a CDH singularity (blue diamond), whereas an $\fsn^2$ (cyan star) is far away from any CDH. 
In the lower two panels (C, D) at the singular limit $\delta=0$, the $\fsn^1$ point becomes a CDH singularity (blue star overlapping with blue diamond). \RED{The center subspace of an $\fsn^1$ (resp., an $\fsn^2$) is denoted by a pink plane (resp., a yellow plane). It follows that the center manifolds of both $\fsn^1$ and $\fsn^2$ are transverse to $\ls$.}}
\label{fig:weak-eigen-direc}
\end{figure*}

In the stable folded node case, we have strong and weak eigenvalues $0>\lambda_w>\lambda_s$. The singular \textit{strong canard} is the unique solution corresponding to the strong stable manifold tangent to the strong eigendirection. For each folded node, the corresponding strong canard and the fold surface $\ls$ form a two dimensional trapping region (the funnel) on the attracting branch of $\ms$ such that all solutions in the funnel converge to that folded node. The funnel family of all folded nodes of $\RED{\mathcal{M}}$ and the fold surface $\ls$ then form a three dimensional funnel volume. Trajectories that land inside the funnel volume will be drawn into one of the folded nodes, passing through the fold surface from an attracting $\ms$ to a repelling $\ms$ due to a cancellation of a simple zero in \eqref{proj_slowreduced}, and such solutions are so-called singular canards. 

\subsubsection{Folded saddle node (FSN)}
In \eqref{desingu_slowreduced}, a degenerate singularity arises when a second eigenvalue, $\lambda_w=0$, becomes zero. This folded singularity is referred to as a ``folded saddle node" ($\fsn$), and is characterized by the condition
\begin{equation} \label{FSNI}
    \begin{array}{rcl}
        \fsn&:=& \{ (V_1,\,w_1,\,V_2,\,w_2) \in \RED{\mathcal{M}}: \vspace{0.1in} \\
        &&f_2 Q(V_1,\,V_2,\,w_2)=\delta P(V_1,\,V_2,\,w_2)\},
    \end{array}
\end{equation}
where \RED{$Q$ and $P$ are defined in \eqref{eq:app-FSN} in Appendix \ref{app:cond-FSN}, which contains a detailed derivation of the FSN condition.} Similar to \citep{Vo2013}, we demonstrate in the appendix that our system can exhibit an $\fsn$ in two different ways: either 
\begin{equation}\label{FSN-1way}
    \begin{array}{rcl}
        \mathrm{FSN}^1 &:=& \{ (V_1,\,w_1,\,V_2,\,w_2) \in \ls: \vspace{0.1in}\\
        &&f_2=\delta K_1(V_1,\,V_2,\,w_2), \ g_1= \delta K_2(V_1,\,V_2,\,w_2)\},
    \end{array}
\end{equation}
or
\begin{equation}\label{FSN-2way}
    \begin{array}{rcl}
        \mathrm{FSN}^2 &:=& \{ (V_1,\,w_1,\,V_2,\,w_2) \in \RED{\mathcal{M}}: \vspace{0.1in} \\
        &&Q(V_1,\,V_2,\,w_2)= \delta K_3(V_1,\,V_2,\,w_2)\},
    \end{array}
\end{equation}
where $K_1, K_2$ and $K_3$ are defined in Appendix \ref{app:cond-FSN}. \RED{We discuss below in Remark \ref{rem:FSN-types} that both $\fsn^1$ and $\fsn^2$ are novel types of $\fsn$ \citep{Letson2017}}.

\RED{
\begin{remark}\label{rem:FSN-types}
    In our parameter regime, the ordinary singularity point lies in the middle branch of critical manifold and is not involved in any bifurcations of the folded singularities. Hence, the $\fsn$ singularities \eqref{FSNI} are neither of type II nor type III \citep{VW2015,KW2010,Roberts2015,Letson2017}. While both $\fsn^1$ and $\fsn^2$ are saddle-node bifurcations of folded singularities, the corresponding center manifolds of these $\fsn$ singularities are transverse to the fold surface at the singularties (see Figure~\ref{fig:weak-eigen-direc}C and D, yellow and pink planes). We follow \citep{Letson2017} to denote both $\fsn^1$ and $\fsn^2$ as novel types of $\fsn$. 
\end{remark}
}

\begin{remark}\label{rm:fsnI}
The condition \eqref{FSN-1way} suggests that an $\fsn^1$ is $O(\delta)$ close to the intersection point of the superslow manifold $M_{ss}$ and the fold surface $\ls$, which was defined as the canard-delayed-Hopf (CDH) singularity in \citep{Letson2017}. In contrast, $\fsn^2$ is always far away from a CDH. Figure \ref{fig:weak-eigen-direc} shows the positions of $\fsn^1$ (blue star), $\fsn^2$ (cyan star) and CDH (blue diamond) in $(V_1,\,V_2,\,w_2)$-space, for $\delta\neq 0$ (top panels) and the singular limit $\delta=0$ (bottom panels). 
\RED{Recall the bifurcation of the upper CDH with respect to $\gsyn$ is shown in Figure~\ref{fig:bd-gsyns}, which explains why there is no upper CDH for $\gsyn=4.3$. }
It is worth noting that a CDH point of \eqref{eq:main} is always a folded singularity because the critical manifold $\ms$ does not depend on the superslow variable $w_2$. 
\end{remark}

\subsection{Slow Layer Problem and Delayed Hopf Bifurcations}\label{sec:delayed-HB}
In this subsection, we turn to the slow layer problem \eqref{eq:slowlayer} resulting from the $\delta\to 0$ singular limit, which exhibits delayed Hopf bifurcations that allow for interesting dynamics. 

Let $J_{\rm SL}$ denote the Jacobian matrix of \eqref{eq:slowlayer} evaluated along the superslow manifold $\mss$, which is given by
\begin{equation}\label{eq:JSL} 
    J_{\rm SL} = \begin{pmatrix}
    \frac{1}{\varepsilon}f_{1V_1} & \frac{1}{\varepsilon}f_{1w_1} & \frac{1}{\varepsilon} f_{1V_2}\\ 
    g_{1V_1} & g_{1w_1} & 0\\
    0 & 0 & f_{2V_2}
    \end{pmatrix},\
\end{equation}
\RED{where the nonzero entries denote partial derivatives.}

The eigenvalues of $J_{\rm SL}$ are given by $f_{2V_2}$ and the eigenvalues of 
\begin{equation*} 
    J = \begin{pmatrix}
   \frac{1}{\varepsilon}f_{1V_1} & \frac{1}{\varepsilon}f_{1w_1}\\ g_{1V_1} & g_{1w_1}
    \end{pmatrix}.\
\end{equation*}
Thus, the Hopf bifurcation points on $\mss$ are given by ${\rm tr}(J)=\frac{1}{\varepsilon}f_{1V_1} + g_{1w_1} = 0$ and $\det J>0$. The former defining condition can be rewritten as
\begin{equation}\label{cond:H_ss}
    \mss^H:=\left\{(V_1,\,w_1,\,V_2,\,w_2) \in \mss :\, f_{1V_1}=-\varepsilon g_{1w_1} \right\}.
\end{equation}
\begin{remark}\label{rem:dh}
It follows from \eqref{cond:H_ss} that an $\mss^H$ is $O(\varepsilon)$ close to the intersection of $\mss$ and $\ls$, i.e., a CDH singularity. The subsystem HB bifurcation $\mss^H$ is also known as delayed Hopf bifurcation (DHB).
\end{remark}

The isolated fold bifurcation points on $\mss$ are located by letting $\det J_{\rm SL}=0$.  That is, 
\begin{equation}\label{cond:L_ss}
    \begin{array}{rcl}
        L_{ss} &:=& \{(V_1,\,w_1,\,V_2,\,w_2) \in \mss :\, \\
        && f_{2 V_2} = 0 \text{ or } f_{1 V_1} g_{1w_1}- f_{1w_1} g_{1 V_1}=0 \}.
    \end{array}
\end{equation}
The fold points $L_{ss}$ that satisfy the former condition (denoted as $L_{ss}^1$) are the folds of the $V_2$-nullcline (see Figure \ref{fig:transition}, green circle and green triangle), which correspond to the transition between superslow dynamics along $\mss$ and slow jumps. $L_{ss}$ given by the latter condition (denoted as $L_{ss}^2$) corresponds to the actual fold of $\mss$ when projected to $(V_1,w_1,V_2)$-space. Since $g_{1w_1}<0$, $f_{1w_1}<0$ and $g_{1V_1}>0$, it follows that $L_{ss}^2$ lies on the middle branch of $\ms$ ($f_{1V_1}=f_{1w_1} g_{1 V_1}/g_{1w_1}>0$) and hence will not play a role in dynamics. At the double singular limit $(\varepsilon, \delta)\to (0,0)$, the $L_{ss}^2$ fold point of $\mss$ will occur at the fold curve of the critical manifold and become a CDH. This can be shown by analyzing the slow reduced layer problem \eqref{eq:slowreducedlayer} obtained from the double singular limits. 

\subsection{Interaction between canard and delayed hopf} \label{subsec:interaction}

To investigate the interaction between the canard and the delayed Hopf mechanisms in the double limit case ($\varepsilon\to 0, \,\delta\to 0$), we need to examine the slow reduced layer problem \eqref{eq:slowreducedlayer}. The corresponding desingularized system is given by  
\begin{equation} \label{desingu-2}
    \begin{array}{rcl}
        \frac{dV_1}{dt_{d}} & = & F(V_1,\,w_1,\,V_2,\,w_2), \vspace{.1in}\\
        \frac{dV_2}{dt_{d}} & = & G(V_1,\,V_2,\,w_2), \vspace{.1in}\\
        \frac{dw_2}{dt_{d}} & = & 0,
    \end{array}
\end{equation}
where $F$ and $G$ are defined in (\ref{desingu_slowreduced}). Note that \eqref{desingu-2} is the $\delta\to 0 $ limit of the desingularized system \eqref{desingu_slowreduced} from subsection \ref{subsec:canard}. The folded singularities of \eqref{desingu-2} are exactly the same as $\RED{\mathcal{M}}$ given by \eqref{folded},
whereas the ordinary singularities of \eqref{desingu-2} are relaxed to be $\mss$. 
The FSN condition at the double singular limit can be obtained from letting $\delta\to 0$ in the $\fsn^1$ condition \eqref{FSN-1way} or the $\fsn^2$ condition \eqref{FSN-2way}.  
This implies that, at the double singular limit, an FSN$^1$ becomes a CDH singularity (see Figure \ref{fig:weak-eigen-direc}C and D)
\begin{equation}\label{eq:FSN1-double-limit}  
    \begin{array}{rcl}
        \mathrm{FSN}^1_{(\varepsilon, \delta)\to (0, 0)} &:=& 
        \{ (V_1,\,w_1,\,V_2,\,w_2) \in \ls: f_2= g_1= 0\} \\
        &=& \mss\cap \ls,
    \end{array}
\end{equation}
whereas an $\fsn^2$ singularity is characterized by
\begin{equation}\label{eq:FSN2-double-limit}
    \begin{array}{rcl}
        \mathrm{FSN}^2_{(\varepsilon, \delta)\to (0, 0)} &:=& \{ (V_1,\,w_1,\,V_2,\,w_2) \in \RED{\mathcal{M}}: \\ &&Q(V_1,\,V_2,\,w_2)=0\}.
    \end{array}
\end{equation}

According to Remarks \ref{rm:fsnI} and \ref{rem:dh}, an $\fsn^1$ singularity from the $\varepsilon$ viewpoint converges to a CDH as $\delta\to 0$ and a DHB $\mss^H$ from the $\delta$ viewpoint converges to a CDH as $\varepsilon\to 0$. It is natural to expect that a CDH singularity point should serve as the interplay between the canard dynamics and the delayed Hopf bifurcation to produce MMOs, as seen in \citep{Letson2017}. 
However, this is not always the case. Specifically, we find that while the CDH on the upper fold surface $\ls$ (see Figure \ref{fig:weak-eigen-direc}B and D) supports MMOs with a high level of robustness due to the coexistence and interaction of two distinct MMO mechanisms, no MMO dynamcis were observed near the lower CDH. It is worth highlighting that both the upper and lower CDH points in our system are of the same type as the CDH investigated in \citep{Letson2017}. 

The CDH points in system \eqref{eq:main} are $\fsn^1$ singularities at the double singular limit. We prove in Appendix \ref{app:CDH-vc} that the CDH singularity in \eqref{eq:main} is a novel type of saddle-node bifurcation of folded singularities as described in \citep{Letson2017}, with the center manifold of the CDH transverse to the fold $\ls$ of the critical manifold. 
This is further confirmed in Figure~\ref{fig:weak-eigen-direc}C and D, \RED{as discussed in Remark \ref{rem:FSN-types}.}

In the case of $g_{syn}=4.3$ (see the left panels of Figure \ref{fig:weak-eigen-direc}), there is no upper CDH. As a result, there is no coexistence and interaction of canard and delayed Hopf mechanisms, leading to MMOs that are sensitive to variation of timescales (see Section \ref{4p3}). For $g_{\rm syn}=4.4$, an upper CDH exists. We show in Section \ref{4p4} that this CDH serves as an organizing center for the local small-amplitude oscillatory dynamics, which results in robust MMOs through the interplay of the DHB and canard mechanisms. In both cases, there exist CDH points on lower $\ls$. However, MMOs are not observed in the neighborhood of any lower CDH (see Section \ref{4p3} for more discussions). 

\RED{\subsection{Singular orbit construction}\label{sec:singular-orbit}

To understand the dynamics of \eqref{eq:main} using GSPT \citep{Fenichel1979}, we need to construct singular periodic orbits  by concatenating solution segments of singular limit systems. According to GSPT, such a singular oscillation will generically perturb to a periodic solution of the full problem as we move away from the singular limit. 

We now use our analysis from previous subsections to construct singular approximations of \eqref{eq:main}
in order to understand the full system trajectory $\Gamma_{(\varepsilon,\delta)}$. We let $\Gamma^F_{(0,\delta)}$ and $\Gamma^S_{(0,\delta)}$ denote trajectories of the fast layer problem \eqref{eq:fastlayer} and the slow reduced problem \eqref{eq:slowreduced} obtained from the $\varepsilon\to 0$ singular limit. Let $\Gamma^S_{(\varepsilon,0)}$ and $\Gamma^{SS}_{(\varepsilon,0)}$ denote trajectories of the slow layer problem \eqref{eq:slowlayer} and the superslow reduced problem \eqref{eq:ssreduced} obtained from the $\delta\to 0$ singular limit. 
The process of constructing singular orbits at the double singular limits is more complicated since there are more than two singular limit systems. We let $\Gamma^x_{(0,0)}$, $x\in\{\rm F, S, SS\}$ denote the fast, slow and superslow flows at the double singular limits.   

\begin{figure}[!htp]
    \centering
   \includegraphics[width=\linewidth]{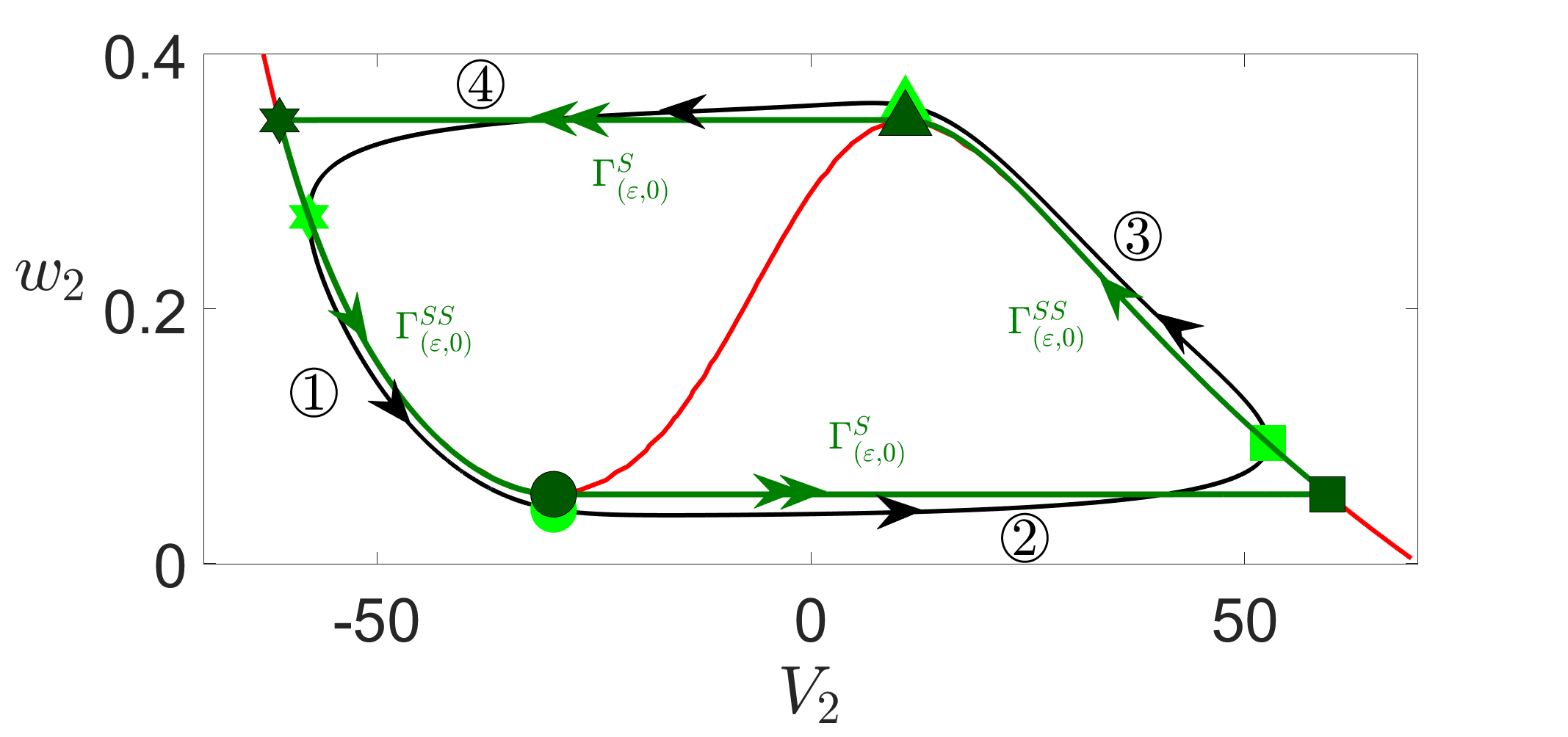}
    \caption{
    Projection of \RED{the singular orbit (green) and} the solution trajectory $\Gamma_{(\varepsilon,\delta)}$ (black) of the full system  \eqref{eq:main} onto the phase plane of ($V_2,w_2$) system with parameters given in Table \ref{tab:par}. The red curve is the $V_2$-nullcline\RED{, which is the projection of the superslow manifold $\mss$.} \RED{The dark and light green symbols mark the key transitions between the slow and superslow sections of the singular orbit and the perturbed oscillation, respectively:} the star and square indicate the transitions from the slow to the superslow motions, and the circle and triangle mark the transitions from superslow to slow sections at the fold of the $V_2$-nullcline. The circled numbers indicate four phases of the oscillations: superslow excursions along $\mss$ during \textcircled{1} and \textcircled{3} and slow jumps at the fold of $\mss$ during \textcircled{2} and \textcircled{4}. }
    \label{fig:transition}
\end{figure}

We start by showing the singular orbit construction for the $(V_2,w_2)$-subsystem (Figure \ref{fig:transition}), which is a relaxation oscillation that is independent of $\gsyn$, $(V_1,w_1)$, and $\varepsilon$. The singular orbit (Figure \ref{fig:transition}, green trajectory) consists of a superslow excursion along the left branch of $V_2$-nullcline ($\Gamma^{SS}_{(\varepsilon,0)}$, phase \textcircled{1}), a slow jump at the lower fold of $V_2$-nullcline up to its right branch ($\Gamma^S_{(\varepsilon,0)}$, phase \textcircled{2}), a superslow excursion through the active phase ($\Gamma^{SS}_{(\varepsilon,0)}$, phase \textcircled{3}) and a slow jump back to its left branch ($\Gamma^S_{(\varepsilon,0)}$, phase \textcircled{4}). Green symbols mark points at the key transition between the four different sections of the oscillation and will be used for later analysis. When projected onto $(V_2,w_2)$-space, the $V_2$-nullcline (red) corresponds to the superslow manifold $\mss$ and the singular orbit $\Gamma^S_{(\varepsilon,0)}\cup \Gamma^{SS}_{(\varepsilon,0)}$ overlaps with the singular orbit from the double singular limits due to its independence on $\varepsilon$. For sufficiently small perturbation $\delta$, $\Gamma^S_{(\varepsilon,0)}\cup \Gamma^{SS}_{(\varepsilon,0)}$ perturbs to the full system trajectory $\Gamma_{(\varepsilon,\delta)}$ shown by the black curve.

\begin{figure*}[!htp]
    \begin{center}
        \begin{tabular}{@{}p{0.48\linewidth}@{\quad}p{0.48\linewidth}@{}}
            \subfigimg[width=\linewidth]{\bf{\small{(A)}}}{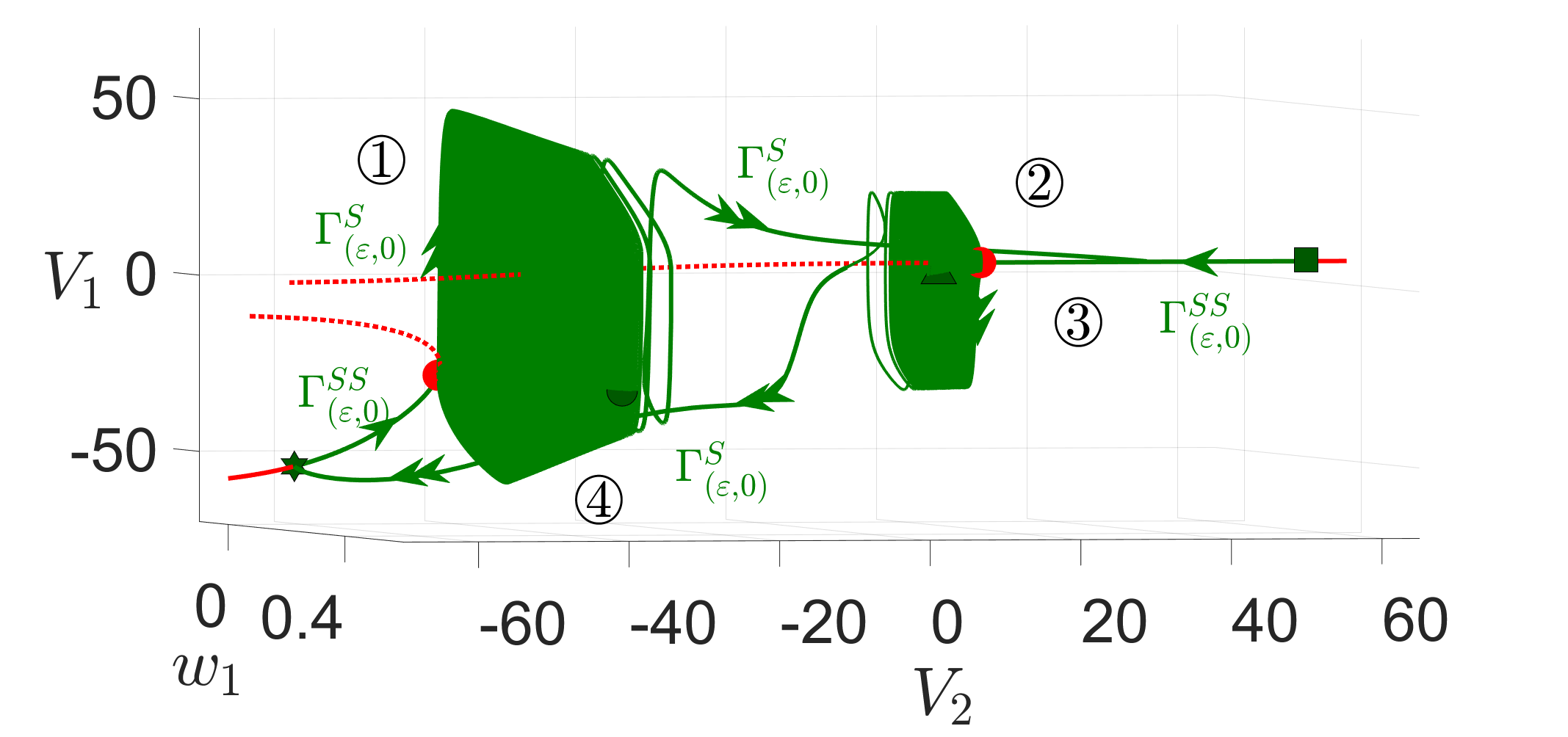} &
            \subfigimg[width=\linewidth]{\bf{\small{(B)}}}{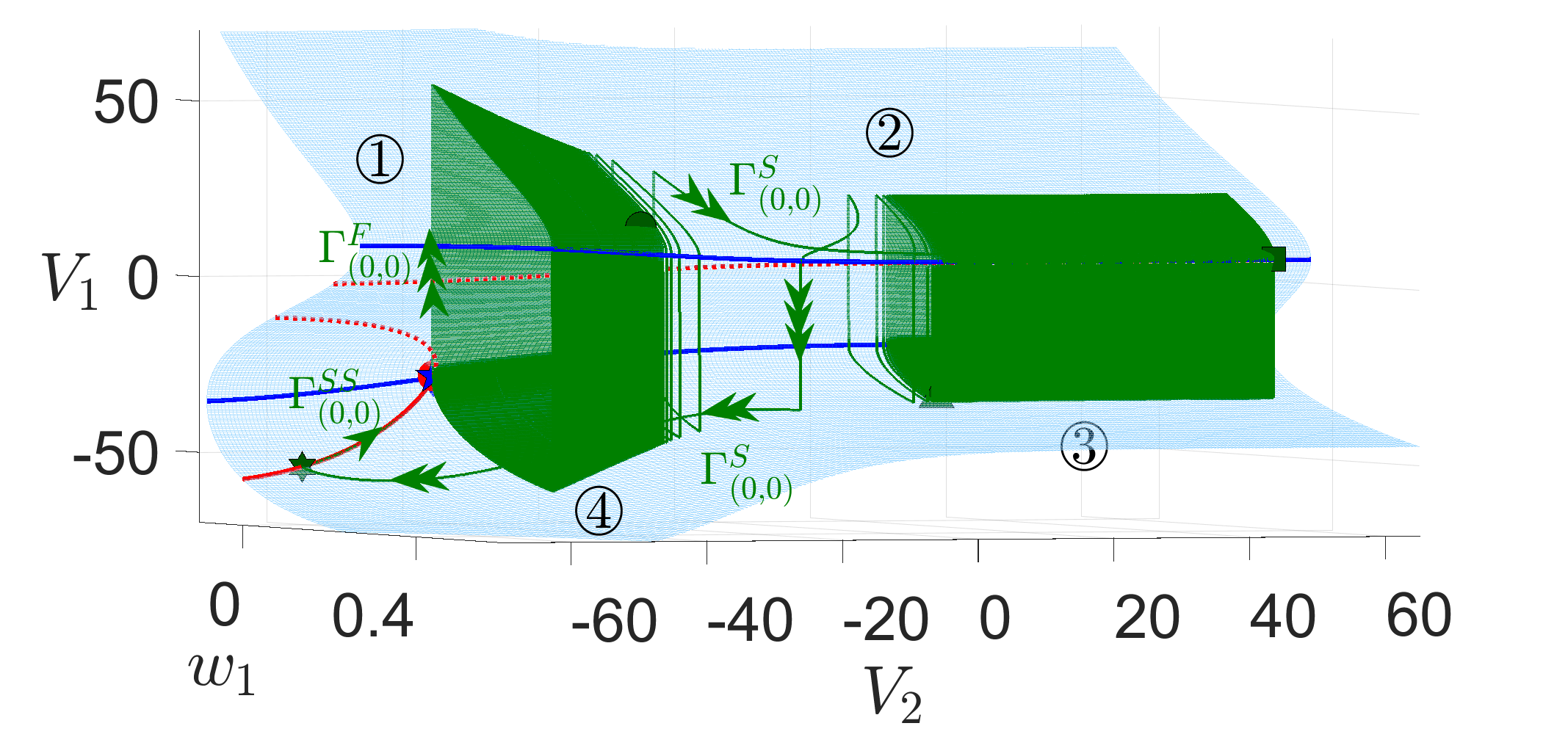}\\
            \subfigimg[width=\linewidth]{\bf{\small{(C)}}}{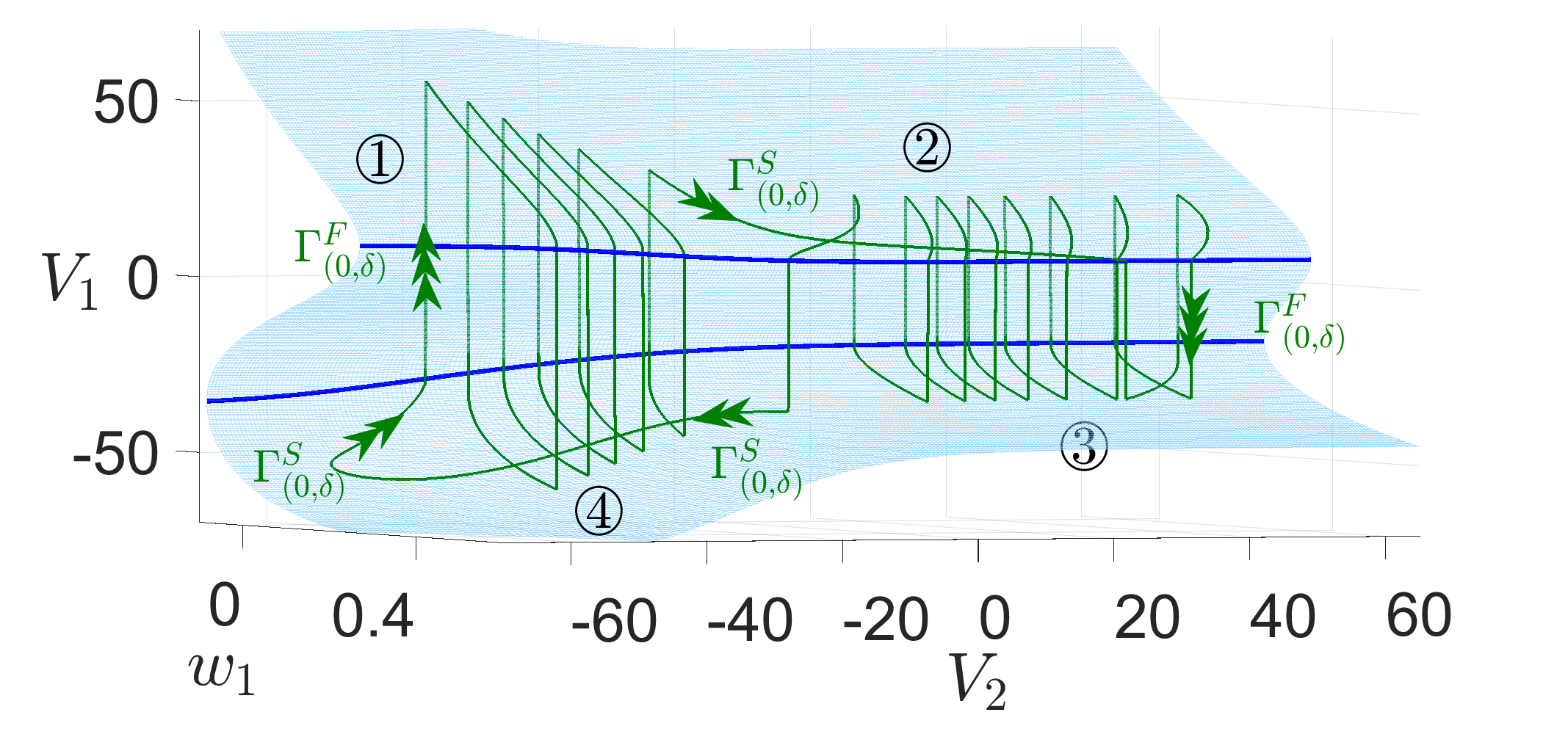} &
            \subfigimg[width=\linewidth]{\bf{\small{(D)}}}{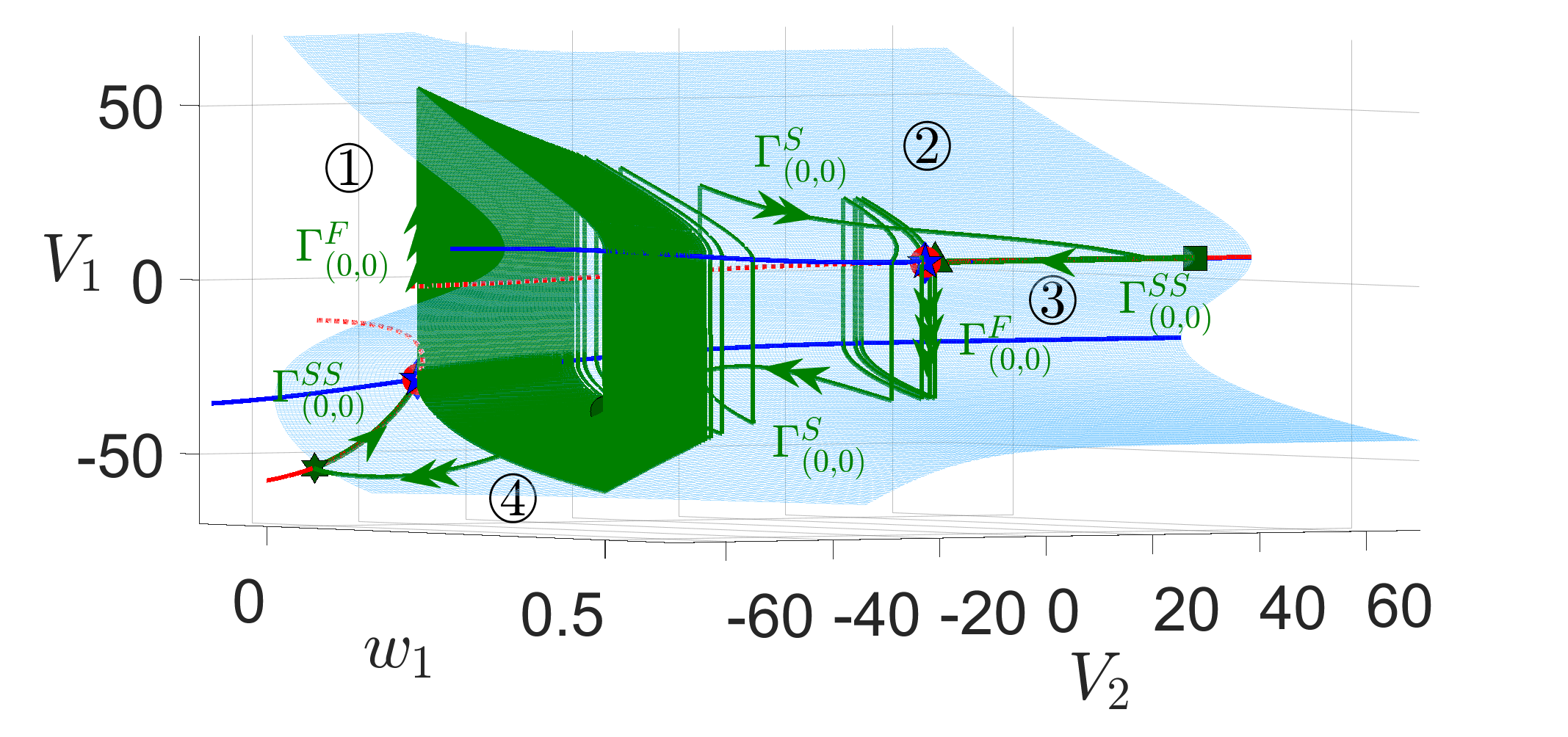}
        \end{tabular}
    \end{center}
\caption{\RED{Projection of singular periodic orbit (green curve) of system \eqref{eq:main} for $g_{\rm syn}=4.3$ (A, B, C) and $g_{\rm syn}=4.4$ (D) to $(V_1, w_1, V_2)$-space. Specifically, panels (A) and (C) show the singular orbits for $(\varepsilon, \delta)=(0.1,0)$, and $(\varepsilon, \delta)=(0,0.053)$, respectively. Panels (B) and (D) show the singular orbits at the double singular limit $(\varepsilon,\delta) =(0, 0)$.
Also shown are the superslow manifold $M_{SS}$ (red curves), the critical manifold $M_S$ (blue surface) and folds of $M_S$ (blue curves). The solid (resp., dashed) red curves are the branches of $M_{SS}$ that are attracting (resp., repelling) under the superslow flow. Green symbols represent transition points between slow and superslow flow for the $(V_2,w_2)$ oscillation as shown in Figure \ref{fig:transition}. The red circle denotes the fast subsystem DHB $M_{SS}^H$. $\fsn^1$ and CDH are denoted by the blue star and diamond, respectively.  }}
\label{fig:singular-orbit}
\end{figure*}

Figure \ref{fig:singular-orbit}A illustrates the singular orbit $\Gamma^S_{(\varepsilon,0)}\cup\Gamma^{SS}_{(\varepsilon,0)}$ for $\gsyn=4.3$ in $(V_1, V_2, w_1)$-space, together with the superslow manifold $\mss$ (red solid curve $\mss^a$: attracting; red dashed curve $\mss^r$: repelling). In other panels, the critical manifold $M_S$ (blue surface) and its folds $\ls$ (blue curves) are also plotted. $\ms$ is separated into three sheets by the folds $\ls$, in which the upper ($\ms^U$) and lower ($M_S^L$) branches are stable and the middle branch ($M_S^M$) is unstable. 

Starting from the green star in panel (A), the singular orbit is in phase \textcircled{1} and evolves along the lower branch of $\mss$ under the superslow reduced problem \eqref{eq:ssreduced}. After hitting the DHB where the stability of $\mss$ changes, the slow layer problem \eqref{eq:slowlayer} takes over but with $V_2$ still evolving on a superslow timescale (see Figure \ref{fig:transition}, phase \textcircled{1}). Thus, the singular orbit during the rest of this phase is a continuum of $(V_1, w_1)$ relaxation oscillations. 
As the evolution speed of $V_2$ changes from superslow to slow at the green circle (beginning of phase \textcircled{2}), a few more spikes occur before the slow flow $\Gamma^S_{(\varepsilon,0)}$ travels to the green square on $\mss$, at which phase \textcircled{3} begins. After that, the superslow reduced problem takes over until the singular orbit $\Gamma^{SS}_{(\varepsilon,0)}$ reaches the upper DHB at which $\mss$ becomes unstable. As such, a family of $(V_1,w_1)$ oscillations emerges as we observed during phase \textcircled{1}. As phase \textcircled{4} begins at the green triangle, several additional $V_1$ spikes occur before the slow flow $\Gamma^S_{(\varepsilon,0)}$ travels back to the green star, thus completing a full cycle. 

Figure \ref{fig:singular-orbit}B shows the singular orbit at the double singular limits. It closely resembles the orbit in panel (A), with the notable exception that there is no longer a superslow segment along the upper $\mss$. Instead, we observe a continuum of $V_1$ spikes throughout phase \textcircled{3}. This is because the upper $\mss$ becomes unstable for all $V_2$ values as $\varepsilon\to 0$, which will be discussed further in subsection \ref{sec:effect_on_DHB}. In (B), the singular orbit consists of $\Gamma^F_{(0,0)}$ (triple arrow) that are fast $V_1$ jumps from $\ls$, $\Gamma^S_{(0,0)}$ (double arrow) which travels along stable branches of $\ms$ or the intersection of $\ms$ and the $V_2$-nullcline, and $\Gamma^{SS}_{(0,0)}$ (single arrow) that follows the stable branch of $\mss$. 

While the $V_2$ value at the slow/superslow transition (green circle or green triangle) is uniquely determined, the corresponding $(V_1,w_1)$ values can assume arbitrary positions on the upper or lower sheet of $\ms$ for that fixed $V_2$. Thus, infinitely many singular orbit segments can be constructed during phases \textcircled{2} and \textcircled{4}, although we only plot one in panel (B) for clarity. For the same reason, the singular orbits in (A) are also not unique since there exist infinitely many ways to choose a starting position for $\Gamma^S_{(\varepsilon,0)}$ during phases \textcircled{2} and \textcircled{4}. 

Comparing (A) and (B) indicates that to obtain MMO dynamics in the perturbed system, $\varepsilon$ cannot be too small. To explain the MMO dynamics for $\gsyn=4.3$ in section \ref{4p3}, we mainly make use of the singular orbit in panel (A) with $0<\varepsilon\ll 1$. 
This orbit, aside from the upper superslow segment where the stability of $\mss$ differs between the two panels, can be viewed as an $\mathcal{O}(\varepsilon)$ perturbation of the singular orbit in panel (B). Hence, when discussing the fast-slow oscillations in the full system, we still refer to different segments of them as 
being governed by \eqref{eq:fastlayer} which describes fast $V_1$ jumps and \eqref{eq:slowreducedlayer} which describes the slow motion along $\ms$.  

For completeness, we also plot the singular orbit for $\varepsilon=0$ but $\delta\neq 0$ in Figure \ref{fig:singular-orbit}C. Instead of a continuum of $V_1$ spikes, we observe a finite number of $V_1$ spikes during phases \textcircled{1} and \textcircled{3} since $\delta\neq 0$. In this limit, the fast segment ($\Gamma^F_{(0,\delta)}$, triple arrow) is the same as $\Gamma^F_{(0,0)}$, whereas the slow segments $\Gamma^S_{(0,\delta)}$ are $\mathcal{O}(\delta)$ perturbations of $\Gamma^S_{(0,0)}$ or $\Gamma^{SS}_{(0,0)}$ from (B). The latter has also been illustrated in Figure \ref{fig:transition}. 

Lastly, the singular orbit for $\gsyn=4.4$ at the double singular limits is shown in Figure \ref{fig:singular-orbit}D. There are two major differences between $\gsyn=4.3$ and $\gsyn=4.4$: Firstly, in contrast to the $\gsyn=4.3$ case where the upper DHB vanishes, it now converges to the upper CDH at the double singular limits in (D). Secondly, $\Gamma^{SS}_{(0,0)}$ in panel (D) follows the upper stable branch of $\mss$ until reaching a saddle-node bifurcation at the green triangle where $\mss$ becomes unstable. 
This results in a unique construction of singular orbit segment during phase \textcircled{4} due to its unique starting position. Nonetheless, the construction of singular orbits during phase \textcircled{2} remains non-unique as discussed in the $\gsyn=4.3$ case.   
Later in section \ref{4p4}, we show how this singular orbit perturbs to MMOs at $\gsyn=4.4$ for various perturbations. 
}

\begin{figure*}[!t]
    \begin{center}
        \includegraphics[width=0.8\linewidth]{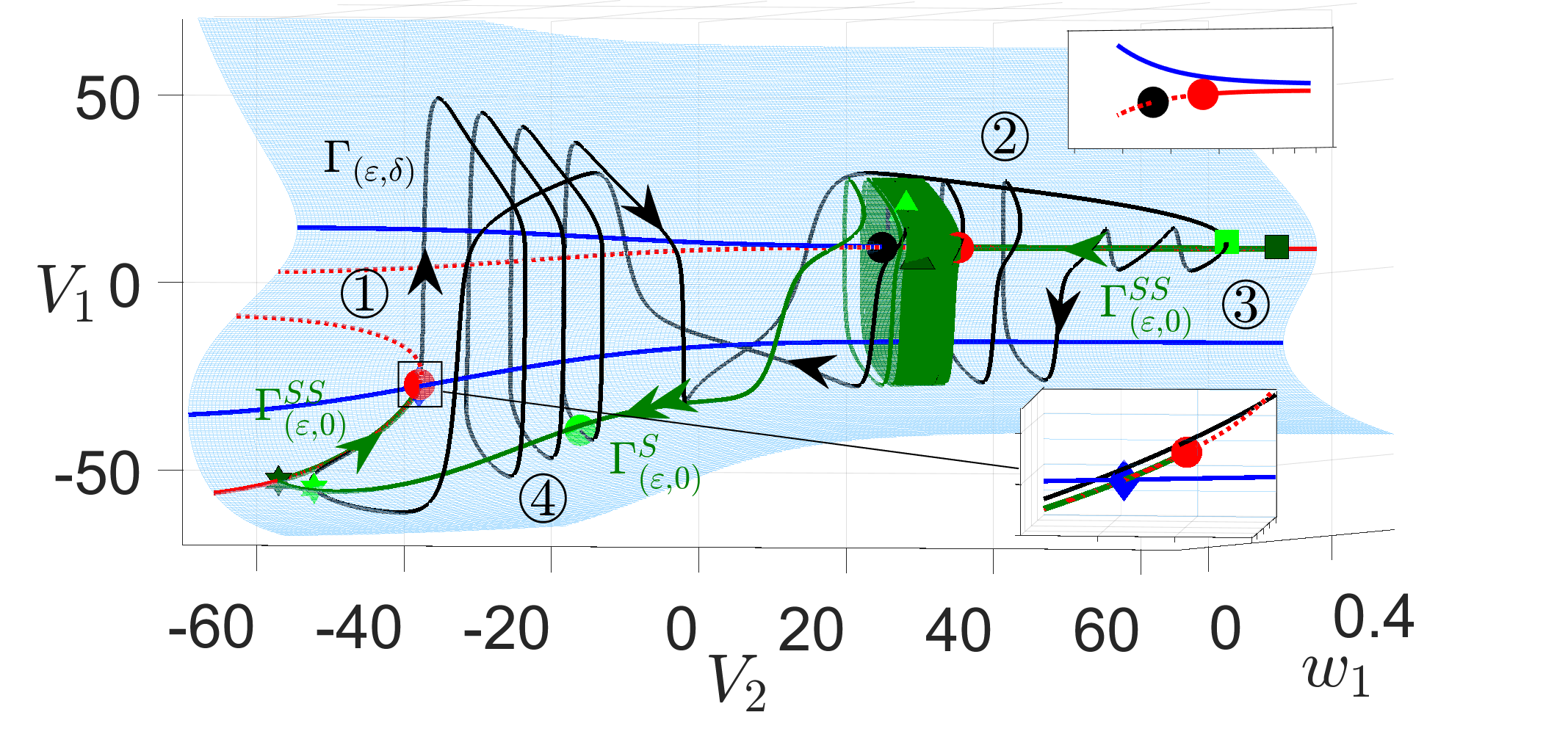}
    \end{center}
\caption{Projection of an attracting MMO solution trajectory (black curve) of system \eqref{eq:main} for $\gsyn=4.3$ from Figure \ref{fig:simu} to $(V_1,w_1,V_2)$-space. \RED{Also shown are parts of the singular orbit from Figure \ref{fig:singular-orbit}A (green curves)},
the critical manifold $M_S$ (blue surface), folds of $M_S$ (blue curves) and the superslow manifold $M_{SS}$ (red curves). The upper inset shows the upper branch of $M_{SS}$ is always below the upper fold of $M_S$, indicating the absence of an upper CDH. The black circle near the upper fold is the true equilibrium of the full system \eqref{eq:main}, which is unstable. The lower inset shows a magnified view around the lower CDH (blue diamond) at which the lower branch of $M_{SS}$ intersects the lower fold of $M_S$. The lower HB bifurcation (red circle) is $O(\varepsilon)$ close to the CDH singularity. \RED{Other color coding and symbols have the same meaning as in Figures \ref{fig:transition} and \ref{fig:singular-orbit}.}}
\label{fig:Ms-gsyn4p3}
\end{figure*}

\section{Analysis of MMOs when $g_{\rm syn}=4.3\ \rm{mS/cm^2} $}\label{4p3}

In this section, we study MMOs that arise when there is no CDH singularity in the middle of the small-amplitude oscillations (SAOs). We show that the only existing mechanism for MMOs at $\gsyn=4.3$ is the delayed-hopf bifurcation (DHB) (see subsection \ref{subsec:MMOs4.3}) and explain why the absence of an upper CDH leads to the sensitivity of MMOs to timescale variations (see subsections \ref{sec:vary-eps-delta-4p3} and \ref{sec:vary_eps_delta_mmo_4p3}). In particular, we explain the complex transitions between MMOs and non-MMOs due to changes of $\varepsilon$ or $\delta$ in subsection \ref{sec:vary_eps_delta_mmo_4p3} (see Figure \ref{fig:c1-phi2}A). 
We also discuss why there is no MMOs near the lower CDH in subsection \ref{sec:noSTOs-CDH}. 

\subsection{Relation of the trajectory to $M_S$ and $\mss$}\label{subsec:MMOs4.3}

The MMO solution of \eqref{eq:main} for $\gsyn=4.3$ from \RED{Figure \ref{fig:simu}} is projected onto the $(V_1,w_1,V_2)-$space in Figure \ref{fig:Ms-gsyn4p3}. The full system equilibrium (black circle) lies on $M_S^M$ and is unstable. The stability of the superslow manifold $\mss$ changes at the two DHBs that are subcritical (red circles). In particular, as $V_2$ decreases, the upper branch of $\mss$ changes from stable-focus (with one negative real eigenvalue and a pair of complex-conjugate eigenvalues whose real parts are negative) to saddle-focus (one negative real eigenvalue and a pair of complex-conjugate eigenvalues whose real parts are positive). Note that the upper fold $\ls$ always lies above the upper branch of the superslow manifold $M_{SS}$ so there is no upper CDH, whereas the lower fold intersects the lower branch of $M_{SS}$ at a CDH singularity (blue diamond in the lower inset). This CDH is a folded focus and will become an FSN$^1$ as $\delta\to 0$ as discussed in \S\ref{subsec:interaction}. Moreover, the nearby HB bifurcation (red circle in the lower inset) is $O(\varepsilon)$ close to this CDH. 

\RED{Figure \ref{fig:Ms-gsyn4p3} shows that the singular orbit from Figure \ref{fig:singular-orbit}A (green curve) is a suitable predictor of the full trajectory (black curve). For the sake of clarity, we choose to not display the entire singular orbit but instead focus on regions where small amplitude oscillations (SAOs) emerge. 
During phase \textcircled{3}, the upper superslow segment $\Gamma^{SS}_{(\varepsilon,0)}$ perturbs to $\Gamma_{(\varepsilon,\delta)}$ which displays two SAOs around the stable branch of $\mss$. 
These SAOs soon transition to large-amplitude oscillations before crossing the DHB at the red circle to reach the unstable branch of $\mss$. 
As phase \textcircled{4} begins at the light green triangle, the slow jump down of $V_2$ brings the trajectory to a region of $M_S^L$ where there is a nearby stable $M_{SS}$ that attracts the trajectory. From the green star, the trajectory follows $M_{SS}$ on the superslow timescale to the lower fold of $M_S$ (see Figure \ref{fig:Ms-gsyn4p3}, lower inset), where it jumps up to $M_S^U$ on the fast timescale under \eqref{eq:fastlayer}, which corresponds to the onset of spikes in $V_1$. The spikes persist until reaching the green circle, at which phase \textcircled{2} begins and the slow jump up of $V_2$ brings the trajectory to the green square, completing a full cycle.}

\RED{We claim that the only mechanism underlying the MMOs at $\gsyn=4.3$ is the DHB mechanism. This is not surprising given that the upper $\Gamma^{SS}_{(\varepsilon,0)}$ switches to a continuum of big spikes when the upper DHB vanishes at the singular limit $\varepsilon= 0$, as shown in Figure \ref{fig:singular-orbit}B.}
To understand why canard dynamics are not involved, we view the trajectory in $(V_1, V_2, w_2)$-space (see Figure \ref{fig:proj-v1v2w2-gsyn4p3}). Unlike Figure \ref{fig:Ms-gsyn4p3} where $\ls$ and curves of folded singularity (\RED{$\mathcal{M}$}) overlap, the $(V_1,V_2,w_2)$-projection captures the structure of the fold surfaces $\ls$ and lets us examine the position of the solution trajectory relative to \RED{$\mathcal{M}$}. To illustrate how SAOs arise from the DHB mechanism, we look at the projection of the trajectory onto $(V_2, V_1)$-plane, which includes the periodic orbits born at the upper Hopf bifurcation (see Figure \ref{fig:Ms-V1V2-gsyn4p3}). 

\subsubsection{Canard mechanism does not contribute to MMOs}\label{sec:no-canard}

\begin{figure}[!htp]
    \begin{center}       \includegraphics[width=0.9\linewidth]{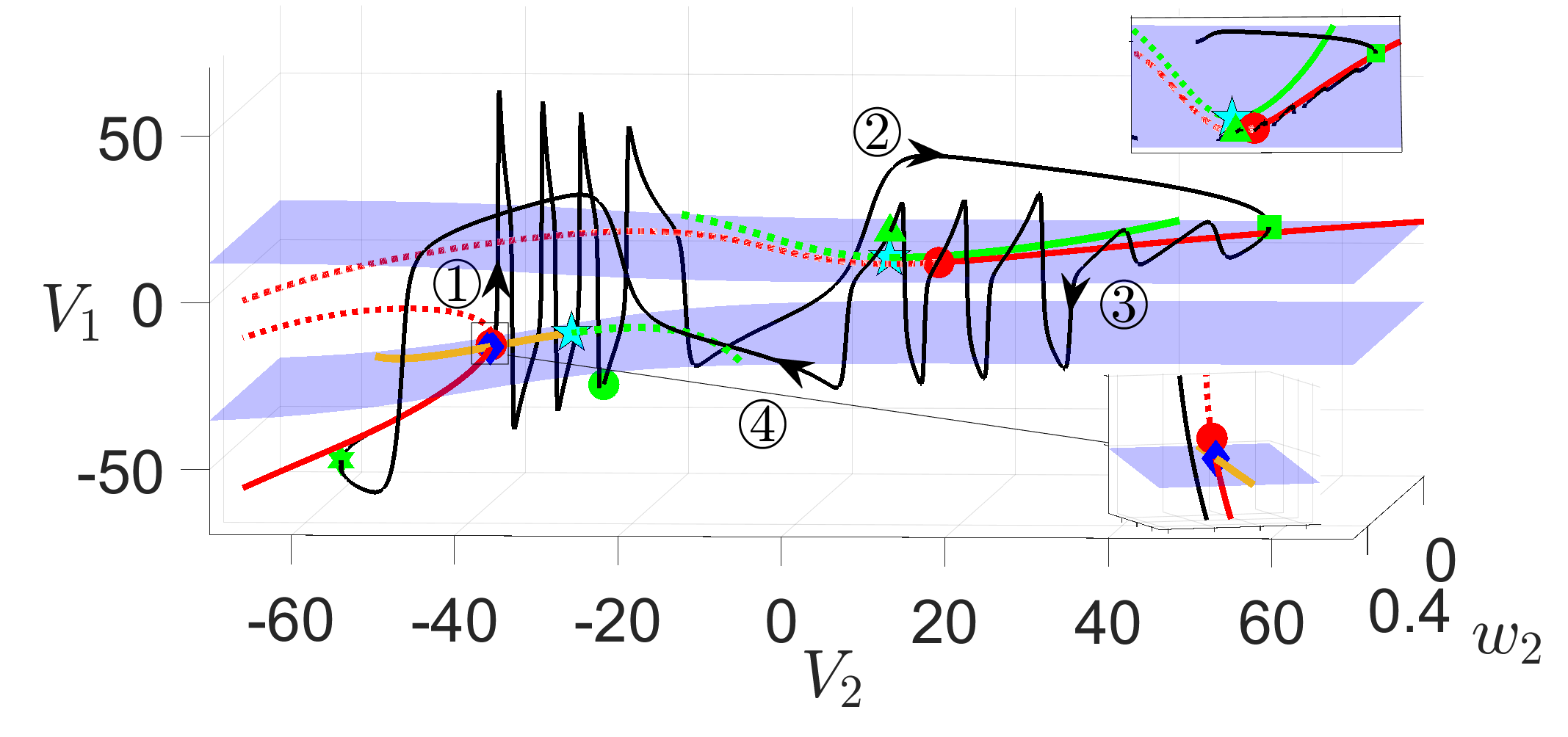}
    \end{center}
\caption{The solution from Figure \ref{fig:Ms-gsyn4p3} projected onto ($V_1,V_2, w_2$)-space. Also shown are the projections of the fold surface $\ls$ (blue surface), superslow manifold $M_{SS}$ (red curve) and folded singularities \RED{$\mathcal{M}$} (green and yellow curves). The cyan star denotes $\fsn^2$ defined in \eqref{FSN-2way}. Color codings of the \RED{folded singularity} curves are the same as in Figure \ref{fig:weak-eigen-direc}. Other symbols have the same meanings as in Figure \ref{fig:Ms-gsyn4p3}.}
\label{fig:proj-v1v2w2-gsyn4p3}
\end{figure}

In Figure \ref{fig:proj-v1v2w2-gsyn4p3}, the solution of \eqref{eq:main} for $\gsyn=4.3$ is projected onto the space $(V_1,V_2,w_2)$ with two separate fold surfaces $\ls$ (blue) and two branches of folded singularities. In the lower \RED{folded singularity} curve, there is an $\fsn^2$ (cyan star) separating the folded singularities that are mostly folded foci (yellow curve) and folded saddle (dashed green). The inset around the lower CDH (blue diamond) shows that the trajectory crosses the lower fold at a regular jump point and hence immediately jumps up to $\ms^U$. The upper \RED{folded singularity} curve also has an $\fsn^2$ (cyan star) marking the boundary between folded saddles (green dashed curve) and stable folded nodes (green solid curve). However, as shown in the upper inset, the trajectory crosses the upper $\ls$ at the normal jump points that are distant from folded nodes. 
Hence, the emergence of SAOs in the MMOs is not due to the canard mechanism. 

\subsubsection{MMOs arise from the delayed Hopf mechanism}\label{4p3-Hopf-MMO}

\begin{figure*}[!t]
    \begin{center}
        \begin{tabular}{@{}p{0.48\linewidth}@{\quad}p{0.48\linewidth}@{}}
            \subfigimg[width=\linewidth]{\bf{\small{(A)}}}{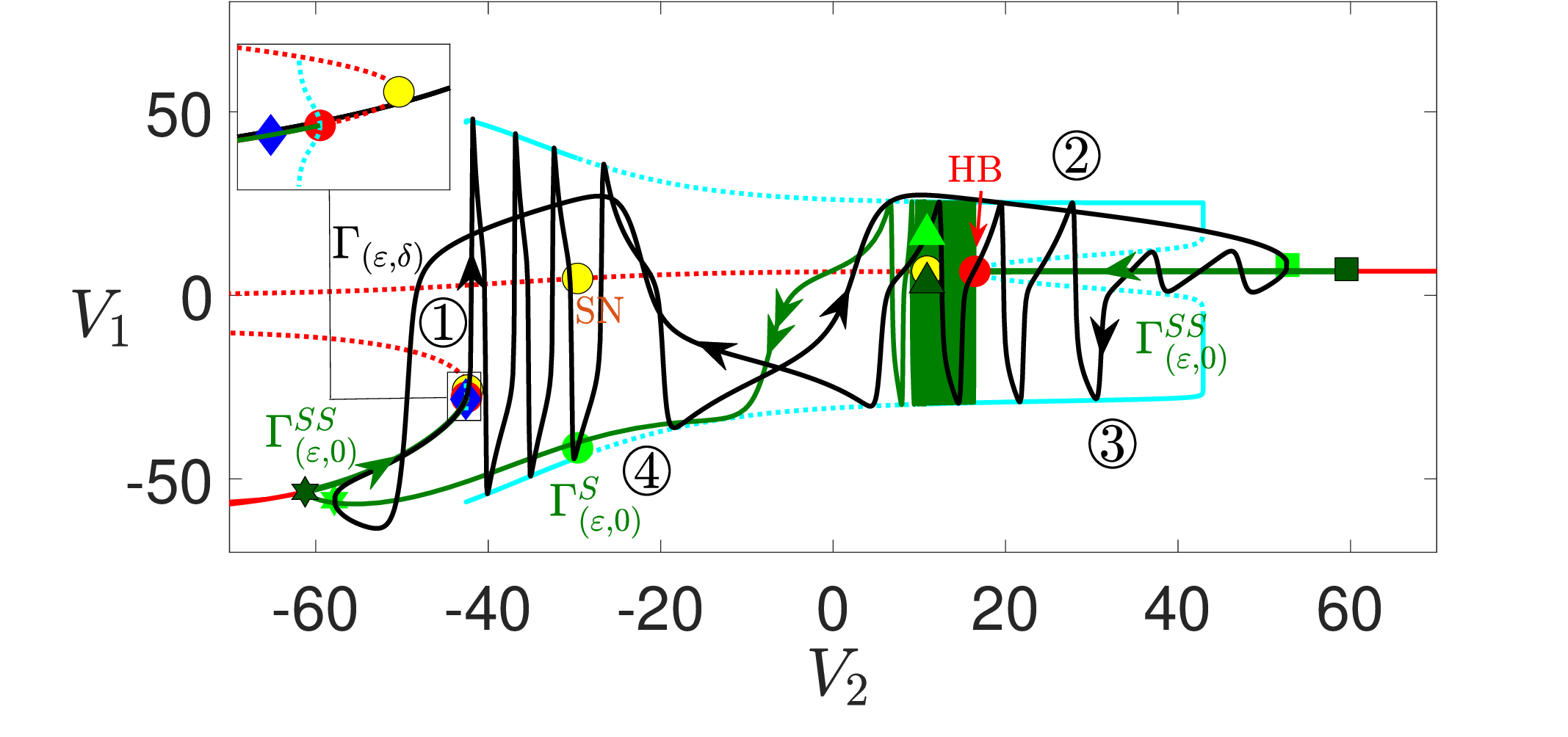} &
            \subfigimg[width=\linewidth]{\bf{\small{(B)}}}{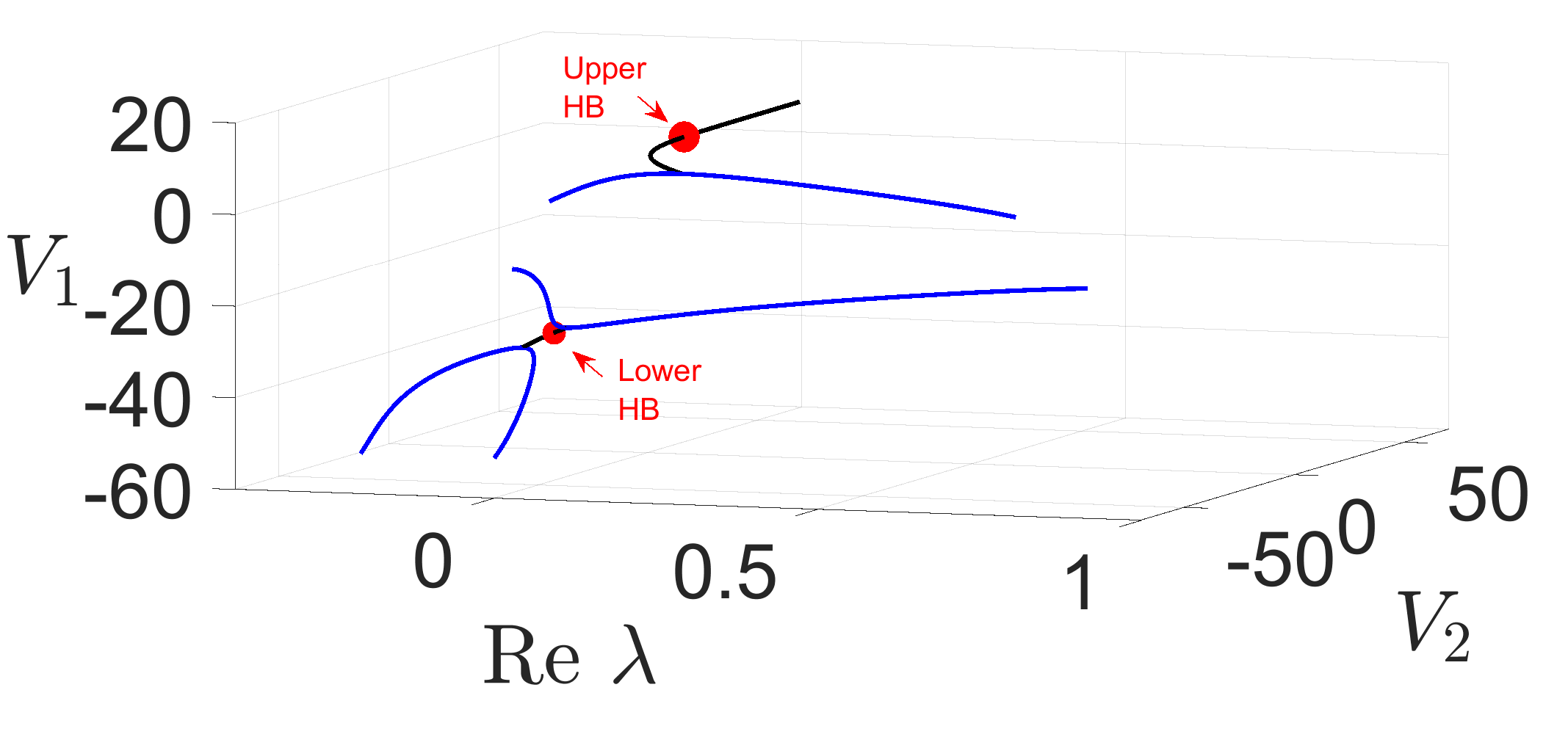}
        \end{tabular}
    \end{center}
\caption{The solutions from Figure \ref{fig:Ms-gsyn4p3}, and the bifurcation diagram of the slow layer problem \eqref{eq:slowlayer} projected onto $(V_2, V_1)$-plane. (A) Besides the trajectory (black)\RED{, the singular orbit (green)} and $\mss$ (red curve), also shown is the periodic orbit (PO) branches (solid cyan: stable; dashed cyan: unstable) born at the upper HB bifurcation. The upper two yellow circles denoting the saddle-node bifurcations of $M_{SS}$ exhibit the same $V_2$ values as the folds of the $V_2$-nullcline (green circle and green triangle). The lower yellow circle (see the inset) represents the actual fold of $\mss$ (denoted as $L_{ss}^2$) and is not a fold of the $V_2$-nullcline. Other symbols have the same meanings as in Figure \ref{fig:Ms-gsyn4p3}. (B) Real part of the eigenvalues of the upper triangular block of $J_{\rm SL}$ \eqref{eq:JSL}, the Jacobian matrix of the slow \RED{layer} problem \eqref{eq:slowlayer} along the superslow manifold $\mss$ from panel (A). The eigenvalues along $\mss$ are real when there are two branches of $\mathrm{Re}(\lambda)$ (blue curves) and complex when there is a single branch of $\mathrm{Re}(\lambda)$ (black curve). }
\label{fig:Ms-V1V2-gsyn4p3}
\end{figure*}

\begin{figure*}[!t]
	\begin{center}
		\begin{tabular}{@{}p{0.48\linewidth}@{\quad}p{0.48\linewidth}@{}}
			\subfigimg[width=\linewidth]{\bf{\small{(A)}}}{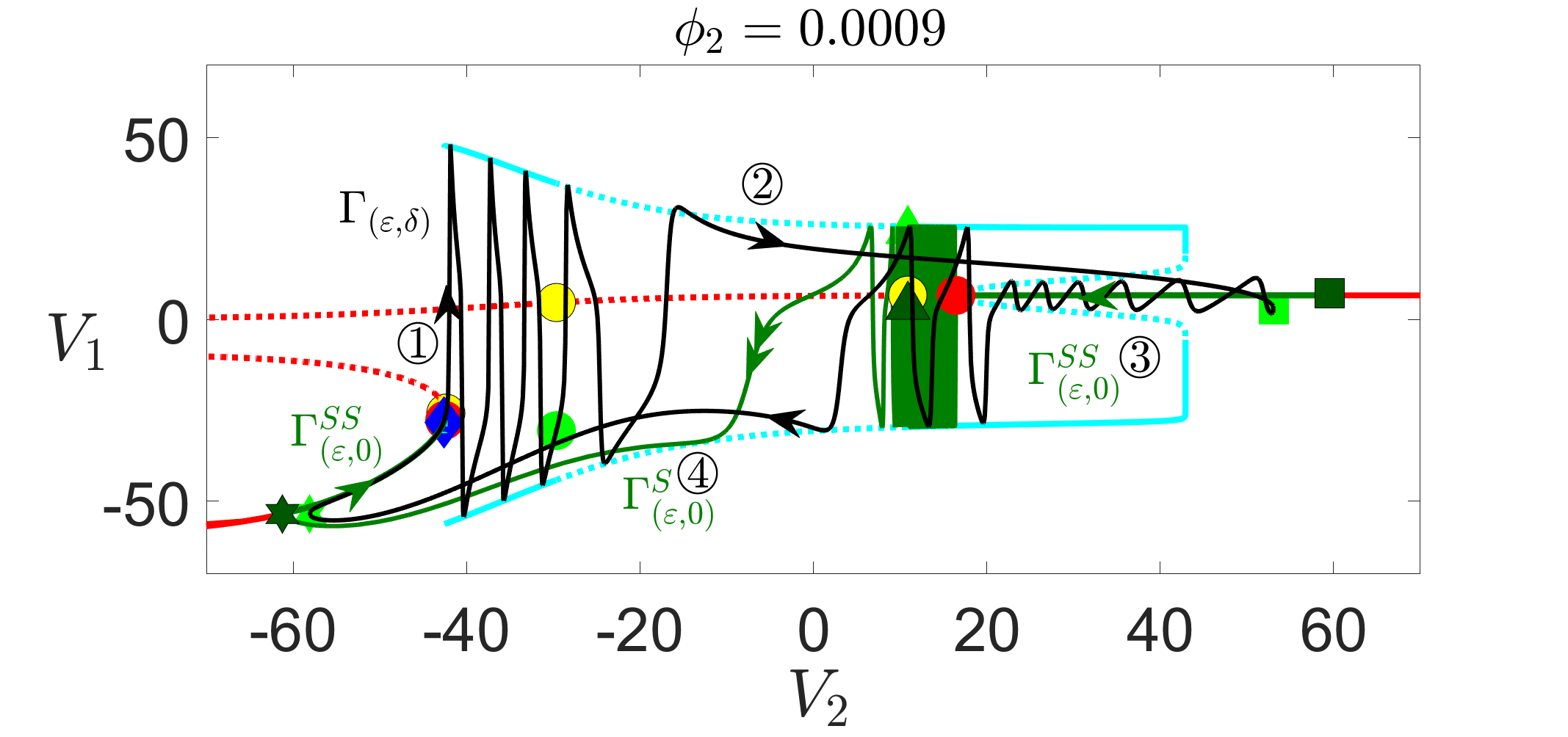} &
			\subfigimg[width=\linewidth]{\bf{\small{(B)}}}{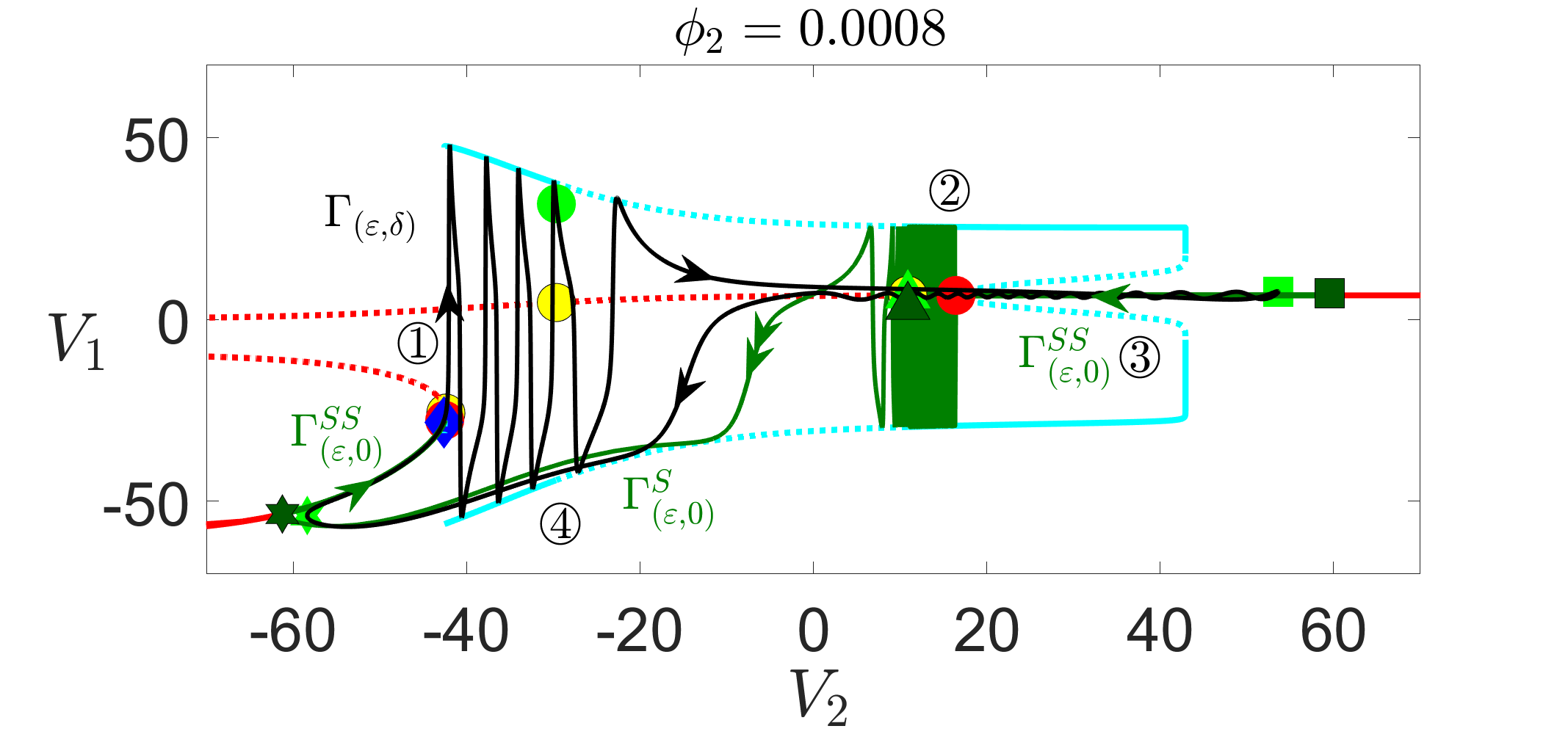}
		\end{tabular}
	\end{center}
	\caption{Projections of solutions of \eqref{eq:main} with $g_{\rm syn}=4.3$, $C_1=8$ and different $\phi_2$ values onto the $(V_2,V_1)$-space.	(A) $\phi_2=0.0009$ (B)  $\phi_2=0.0008$.  All the symbols and curves have the same meaning as in Figure \ref{fig:Ms-V1V2-gsyn4p3}.}
	\label{fig:4p3-phi2-p8-p9}
\end{figure*}

To examine how the DHB mechanism engages in organizing the MMOs, we turn to the subsystems obtained by treating $\delta$ as the main singular perturbation parameter (see \S\ref{sec:delayed-HB}). Figure \ref{fig:Ms-V1V2-gsyn4p3}A shows the projection of the solution, $\mss$ and the bifurcation diagram of the slow layer problem \eqref{eq:slowlayer} onto the $(V_2,V_1)$-plane. As discussed above, the upper $\Gamma^{SS}_{(\varepsilon,0)}$ (green curve) perturbs to SAOs followed with large spikes (black curve). Starting from the green square, the SAOs gradually decrease in magnitude as the trajectory moves towards the upper DHB. After two such oscillations, the trajectory crosses the unstable periodic orbit branch, whose amplitude also decreases as it approaches the upper DHB. Upon crossing, the trajectory undergoes a sudden jump to the outside large-amplitude periodic orbit branch, giving rise to large spikes.

\RED{Recall at default parameters, $\delta \approx 0.053$. Figure \ref{fig:4p3-phi2-p8-p9} shows that decreasing the perturbation $\delta$ improves the singular limit predictions and increases the amount of time the full system trajectory $\Gamma_{(\varepsilon,\delta)}$ spends near the attracting branch of upper $\mss$ (denoted as $\mss^a$). In particular, panel (B) indicates the trajectory with a slightly smaller perturbation $\delta\approx 0.042$ is able to pass over the upper DHB (red circle) to the repelling side of $\mss$ (denoted as $\mss^r$) and there is a delay in which the trajectory traces $\mss^r$ before it jumps away. However, we recognize that $\Gamma_{(\varepsilon,\delta)}$ under the default perturbation $\delta \approx 0.053$ does not exhibit a similar delay (see Figure \ref{fig:Ms-V1V2-gsyn4p3}A), which is contrary to what one would expect in DHB-induced MMOs. While this might initially suggest that $\delta \approx 0.053$ is distant from the singular limit, we show below that this is not the case. Moreover, this specific value of $\delta$ aligns with the perturbation used in \citep{Nan2015}, where GSPT analysis has been successfully employed to elucidate the dynamics across various $\gsyn$ values. 

It is worth emphasizing two interesting points: Firstly, the impact of decreasing $\delta$ on SAOs is non-monotonic. In certain cases, smaller perturbations can lead to fewer SAOs with even larger amplitudes. This effect is related to how the slow flow during phase \textcircled{2} approaches $\mss^a$, a detailed analysis of which is provided in subsection \ref{sec:vary_eps_delta_mmo_4p3}. Therefore, the absence of a delay in Figure \ref{fig:Ms-V1V2-gsyn4p3} does not imply a significant deviation from the singular limit. Rather, it is mainly
due to the manner in which the trajectory approaches $\mss^a$ during phase \textcircled{2}, resulting in small oscillations that cross the unstable inner periodic orbit branch before passing through the DHB. As discussed above (also see Figure \ref{fig:gsyn_4p3-phi2}), a slight reduction or increase in the perturbation size can induce a delay phenomenon. 
Secondly, the plateauing behavior of the trajectory after passing the Hopf bifurcation is somewhat different from what one would expect to see in a typical DHB fashion, which typically involves oscillations with diminishing and then increasing amplitude. This is because the variable $V_2$ switches from superslow to slow timescale at the green triangle shortly after passing the HB, and there is insufficient time for the trajectory to oscillate. Hence, the associated pattern after HB is plateauing, and the amount of time the trajectory spends near $\mss^r$ is significantly shorter than that near $\mss^a$.} 

In the following three subsections, we demonstrate how the absence of the interaction between canard and DHB mechanisms, specifically due to the lack of an upper CDH, can result in the sensitivity of the MMOs to timescale variations. To achieve this, we first explore how changes in the singular perturbation parameters $\varepsilon$ and $\delta$ can induce transitions between MMOs and non-MMOs by analyzing their impacts on the two MMO mechanisms - canard and DHB. Additionally, we provide an explanation for why the lower CDH singularity does not guarantee the occurrence of SAOs.

\subsection{Effects of varying $\varepsilon$ and $\delta$ on DHB and FSN points}\label{sec:vary-eps-delta-4p3}

When $\gsyn=4.3$, the CDH singularity only exists on the lower fold surface of $\ms$. We demonstrate that this leads to different effects of $\varepsilon$ on the upper and lower DHB points, respectively (see Figure~\ref{fig:two-par-4p3}). We also examine the effect of $\delta$ on the lower FSN points (see Figure \ref{fig:lower_fsn_4p3}). 

\subsubsection{Effect of $\varepsilon$ on DHB}
\label{sec:effect_on_DHB}

\begin{figure*}[!htp]
    \begin{center}
        \begin{tabular}{@{}p{0.48\linewidth}@{\quad}p{0.48\linewidth}@{}}
            \subfigimg[width=\linewidth]{\bf{\small{(A)}}}{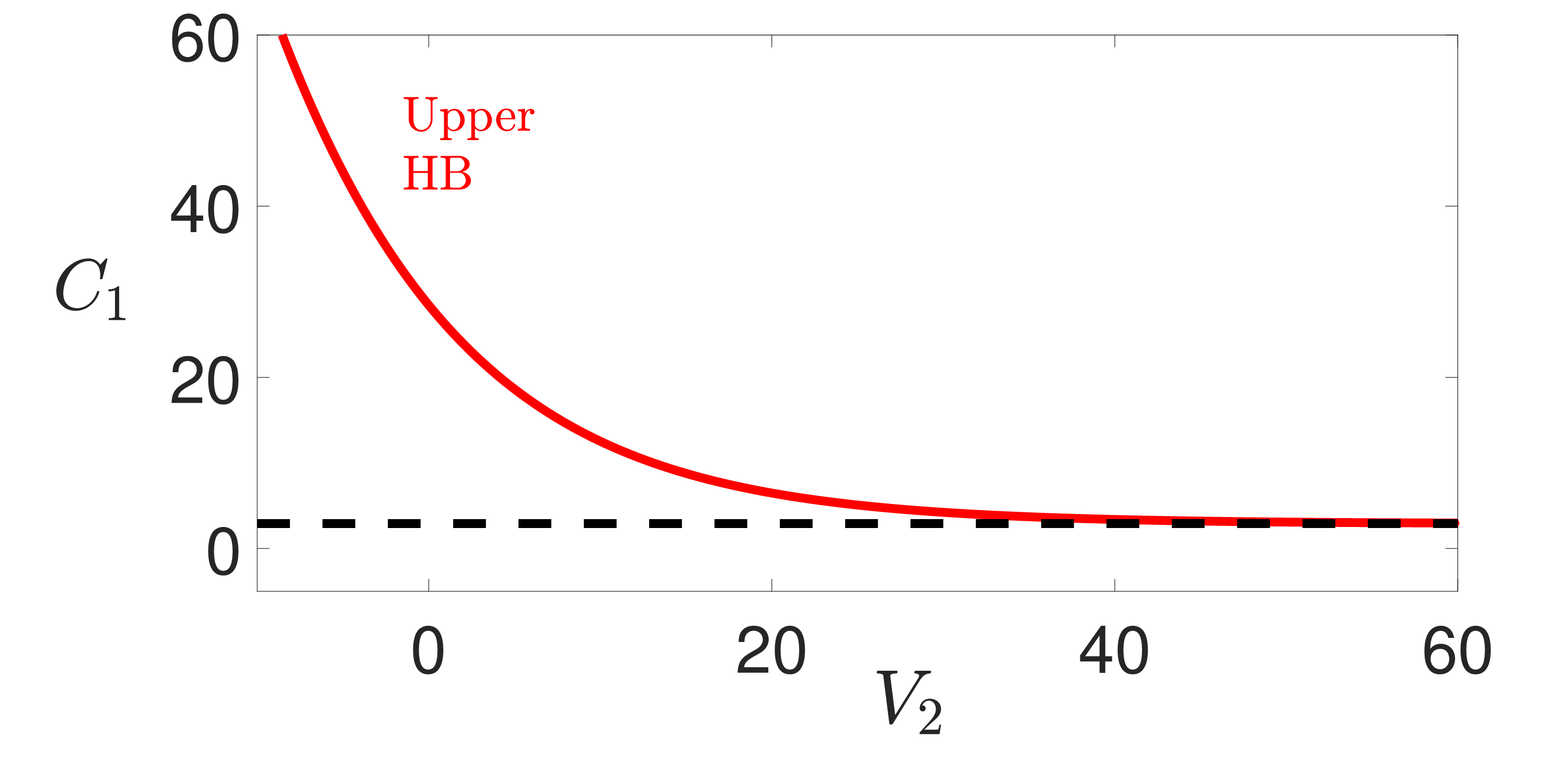}
            &\subfigimg[width=\linewidth]{\bf{\small{(B)}}}{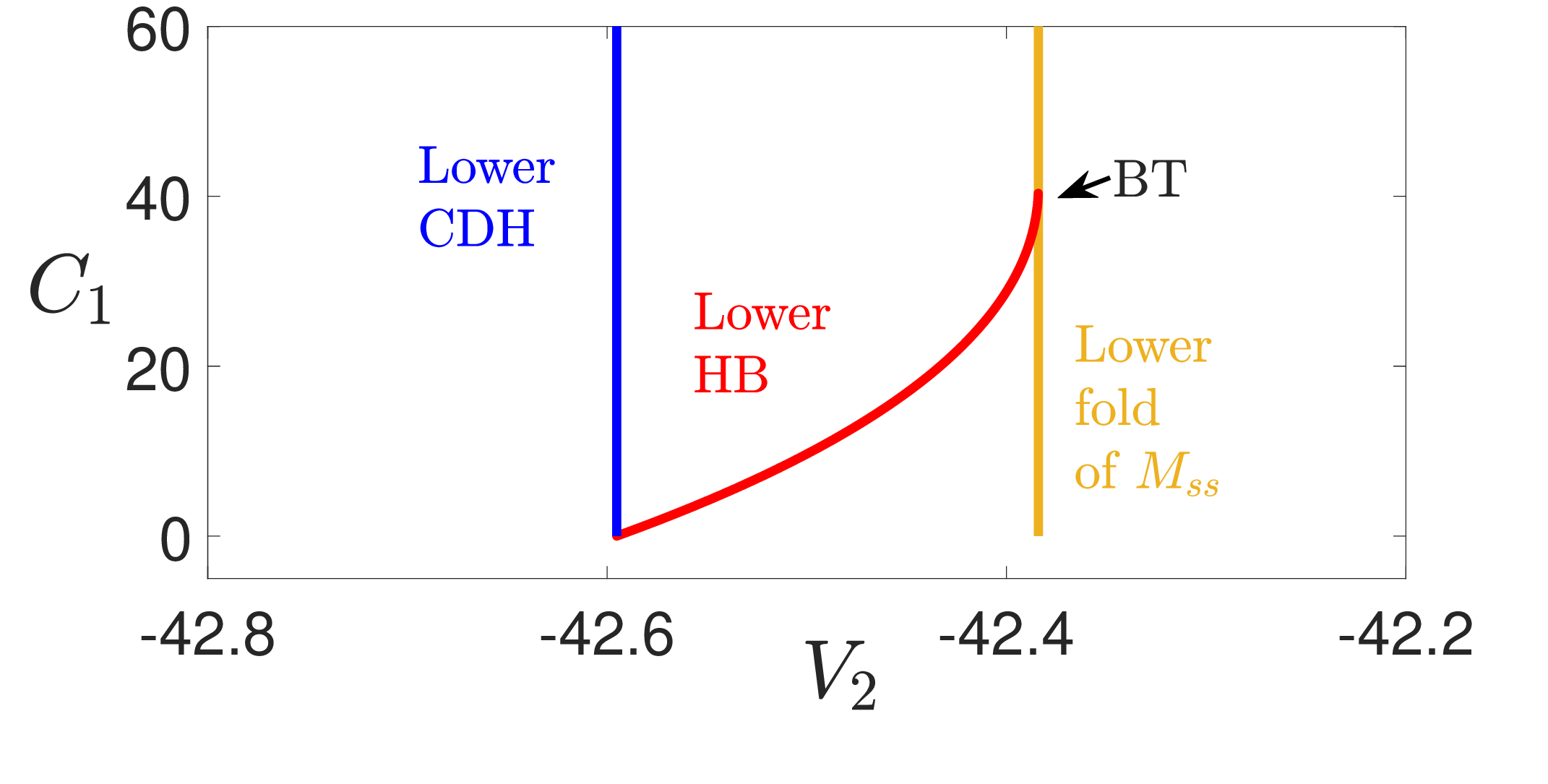}
        \end{tabular}
    \end{center}
\caption{2-parameter bifurcation curves for \eqref{eq:slowlayer} projected onto ($V_2,C_1$)-space when $g_{\rm syn}=4.3$. (A) The bifurcation curve of the upper HB (Figure \ref{fig:Ms-gsyn4p3}, upper red circle), with a horizontal asymptote at $C_1=2.9097$. (B) The bifurcation curves of the lower HB (red curve) and the lower fold $L_{ss}^2$ of $\mss$ (yellow curve), which meet at a Bogdanov-Takens (BT) bifurcation. As $C_1\to 0$ (i.e., $\varepsilon\to 0$), the lower HB converges to the lower CDH (blue). }
\label{fig:two-par-4p3}
\end{figure*}

The effects of $\varepsilon$ on the DHBs $\mss^H$ are summarized by the 2-parameter bifurcation diagrams of \eqref{eq:slowlayer} projected onto $(V_2,C_1)$-space (see Figure \ref{fig:two-par-4p3}). As $C_1$ decreases (or equivalently, as $\varepsilon$ decreases), the upper Hopf moves to larger $V_2$ values and eventually vanishes for $\varepsilon$ small enough (see Figure \ref{fig:two-par-4p3}A). \RED{This explains why the upper singular orbit $\Gamma^{SS}_{(\varepsilon,0)}$ in Figure \ref{fig:singular-orbit}A for $\varepsilon\neq 0$ switches to a continuum of spikes in Figure \ref{fig:singular-orbit}B as $\varepsilon\to 0$.} On the other hand, since there is a CDH on the lower $\ls$, the lower Hopf will converge to that CDH as $\varepsilon\to 0$ (see Figure \ref{fig:two-par-4p3}B and recall Remark \ref{rem:dh}). When $\varepsilon$ increases, the lower Hopf and the lower fold of $\mss$ meet and coalesce at a Bogdanov-Takens (BT) bifurcation. After the BT bifurcation, the Hopf bifurcation disappears. Unlike the upper DHB, the lower DHB is close to the actual fold of $\mss$ (also see Figure \ref{fig:two-par-4p3}A, lower yellow circle).   

\subsubsection{Effect of $\delta$ on FSN points}

\begin{figure*}[!t]
    \begin{center}
        \begin{tabular}{@{}p{0.48\linewidth}@{\quad}p{0.48\linewidth}@{}}
            \subfigimg[width=\linewidth]{\bf{\small{(A)}}}{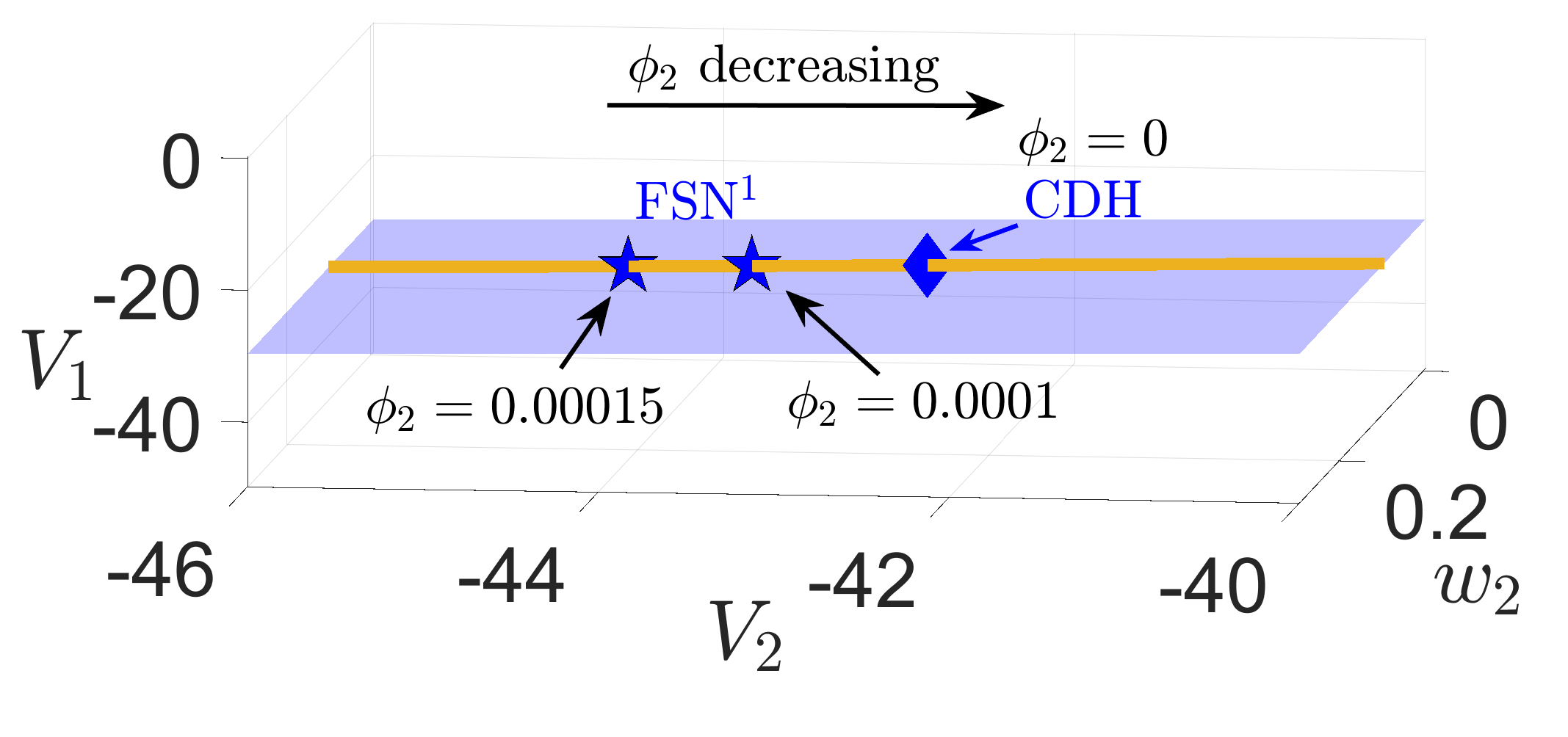}
            &\subfigimg[width=\linewidth]{\bf{\small{(B)}}}{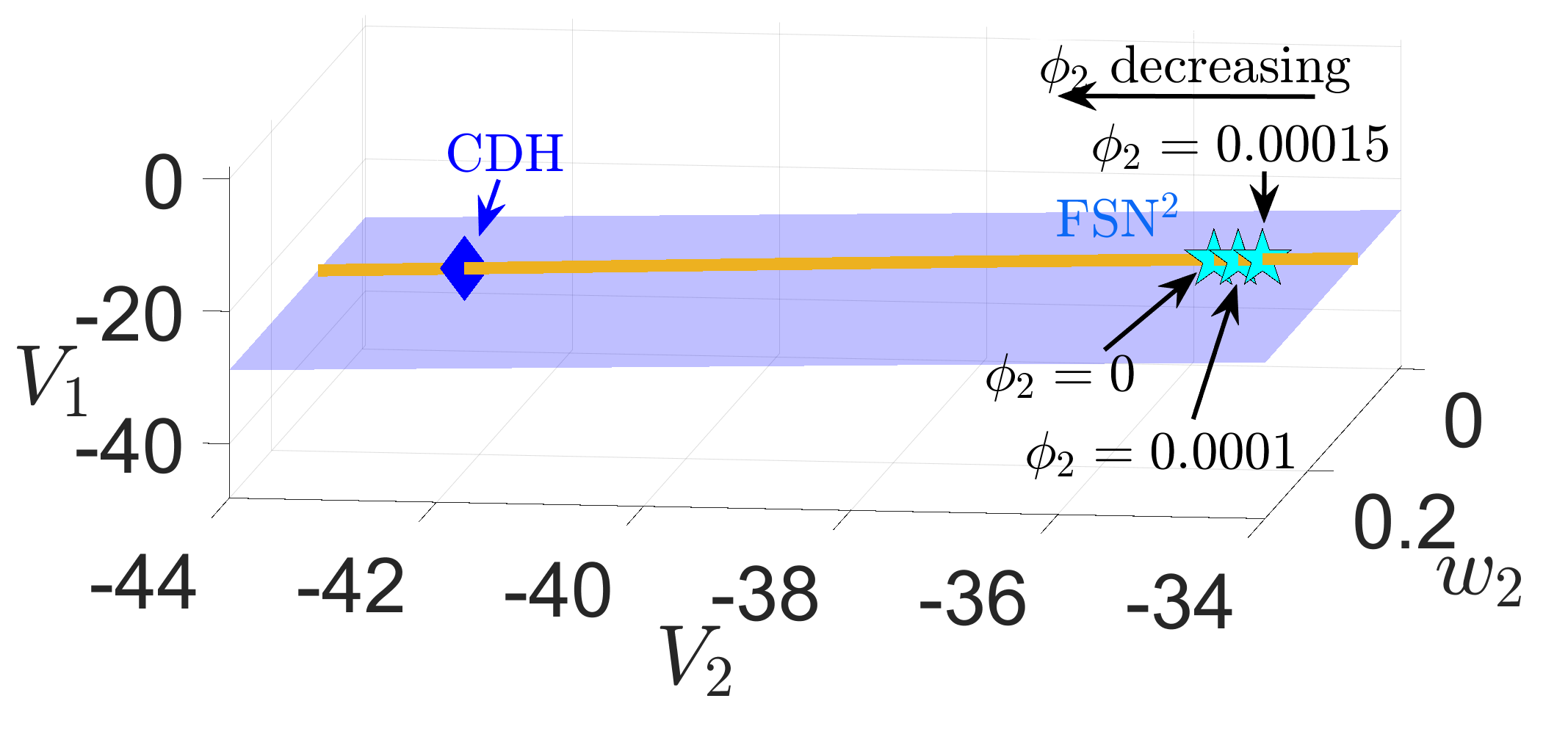}
        \end{tabular}
    \end{center}
\caption{Relation of $\fsn^1$ \eqref{FSN-1way}, $\fsn^2$ \eqref{FSN-2way} and the CDH point on the lower fold for $g_{\rm syn}=4.3$ and various values of $\phi_2$. (A) $\fsn^1$ (blue star) converges to the $\rm CDH$ (blue diamond) as $\phi_2\to 0$ (equivalently, $\delta\to 0$). (B) $\fsn^2$ (cyan star) are located far away from the $\rm CDH$ and converge to the leftmost cyan star at the intersection of \RED{the folded singularity curve $\mathcal{M}$} and $Q(V_1,V_2,w_2)=0$ (see \eqref{eq:FSN2-double-limit}). The yellow curve denotes the lower \RED{folded singularity curve} without showing stability and types. 
}
\label{fig:lower_fsn_4p3}
\end{figure*}

There is no CDH or $\fsn^1$ on the upper fold, hence, we only examine the effect of $\delta$ on FSN singularities on the lower $\ls$. Figure \ref{fig:lower_fsn_4p3}A shows as $\delta\to 0$ (or equivalently, $\phi_2\to 0$), the $\fsn^1$ singularity converges to the lower CDH as demonstrated in our analysis (see \eqref{eq:FSN1-double-limit}). On the other hand,  the $\fsn^2$ singularity is significantly distant from the CDH (see Figure \ref{fig:lower_fsn_4p3}B and recall the condition \eqref{eq:FSN2-double-limit}).

\subsection{Why there are no SAOs near the lower CDH}\label{sec:noSTOs-CDH}

Before examining the effects of varying $\varepsilon$ and $\delta$ on \RED{MMO} dynamics based on their impact on DHB and FSN points, we discuss briefly in this subsection why there are no SAOs near the lower CDH. 

As discussed above (also see Remarks \ref{rm:fsnI} and \ref{rem:dh}), the lower CDH is $O(\delta)$ close to an $\mathrm{FSN}^1$ and $O(\varepsilon)$ close to a DHB. One would naturally expect to observe SAOs arising from Canard and/or DHB mechanisms near this lower CDH. Nonetheless, there is no SAOs near the lower fold regardless of $\varepsilon$ and $\delta$ values. Below we explain why none of the two mechanisms produces SAOs.

As discussed earlier, the reason why there is no canard-induced SAOs when $\varepsilon$ and $\delta$ are at their default values is because the trajectory crosses the fold surface at a regular jump point near which the folded singularities including the CDH are folded focus. This remains to be the case as $\varepsilon$ varies or as $\delta$ increases. While as $\delta$ decreases the trajectory will follow $\mss$ more closely and hence cross the fold surface somewhere near or at a folded node, that folded node is very close to an FSN where canard theory breaks down. As a result, the existence of canard solutions for smaller $\delta$ is not guaranteed. 

On the other hand, with default $\varepsilon$ and $\delta$ values, the trajectory jumps away from the lower fold before reaching the Hopf bifurcation and hence there is no DHB-induced small oscillations. Increasing $\delta$ will make the DHB less relevant, whereas increasing $\varepsilon$ will move the DHB further away from the lower $\ls$ and eventually vanish upon crossing the actual fold of $\mss$ (see Figure \ref{fig:two-par-4p3}B). Thus, we do not expect to detect MMOs with increased $\varepsilon$ or $\delta$. As $\varepsilon$ or $\delta$ \RED{decreases}, the trajectory should pass closer to the lower DHB point. This is because decreasing $\varepsilon$ moves the Hopf bifurcation closer to the CDH singularity (see Figure \ref{fig:two-par-4p3}B) and reducing $\delta$ pushes the trajectory to travel along $\mss$ more closely. However, the reason that no SAOs are induced by the passage through the lower HB is that this HB is relatively close to the actual fold of $\mss$, i.e., close to a double zero eigenvalue at a BT bifurcation of the $(V_1,w_1,V_2)$ subsystem (see Figure \ref{fig:Ms-V1V2-gsyn4p3}A and Figure \ref{fig:two-par-4p3}B). As a result, the branch of unstable small-amplitude periodic orbits born at the lower HB is almost invisible (see the inset of Figure \ref{fig:Ms-V1V2-gsyn4p3}A) and there is only a small region of $\mss$ along which the Jacobian matrix $J_{\rm SL}$ \eqref{eq:JSL} of the slow layer problem \eqref{eq:slowlayer} has complex eigenvalues. Figure \ref{fig:Ms-V1V2-gsyn4p3}B shows the real part of the first two eigenvalues $\lambda$ of \eqref{eq:JSL} evaluated along $\mss$, excluding the third eigenvalue given by $f_{2V_2}$. The eigenvalues are real when there are two branches of curves for $\rm Re\,\lambda$ and become complex when the curves coalesce and become a single branch. In panel (B), the eigenvalues on the stable lower branch of $\mss$ are initially real and negative. That is, the trajectory approaches the attracting $\mss$ along stable nodes of the slow layer problem. As the superslow flow brings the trajectory towards the Hopf bifurcation, the eigenvalues become complex. However, this region of complex eigenvalues is short and $\lambda$ becomes real again shortly after. As a result, the trajectory has insufficient time to oscillate and we do not observe any small-amplitude oscillations before the trajectory jumps up to the outer periodic orbit branch.  

\subsection{Effects of varying $\varepsilon$ and $\delta$ on MMOs}\label{sec:vary_eps_delta_mmo_4p3}

\begin{figure*}[!t]
    \begin{center}
        \begin{tabular}{@{}p{0.48\linewidth}@{\quad}p{0.48\linewidth}@{}}
            \subfigimg[width=\linewidth]{\bf{\small{(A)}}}{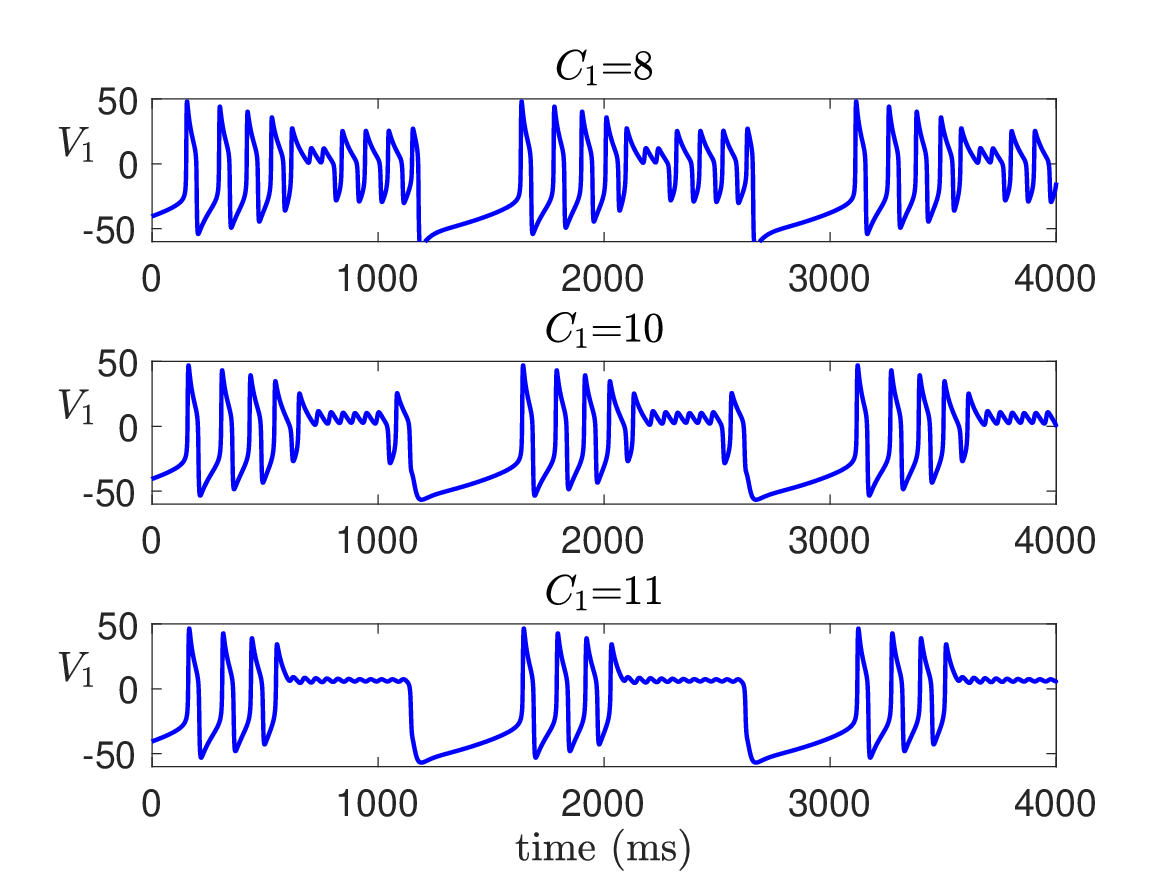}&
            \subfigimg[width=\linewidth]{\bf{\small{(B)}}}{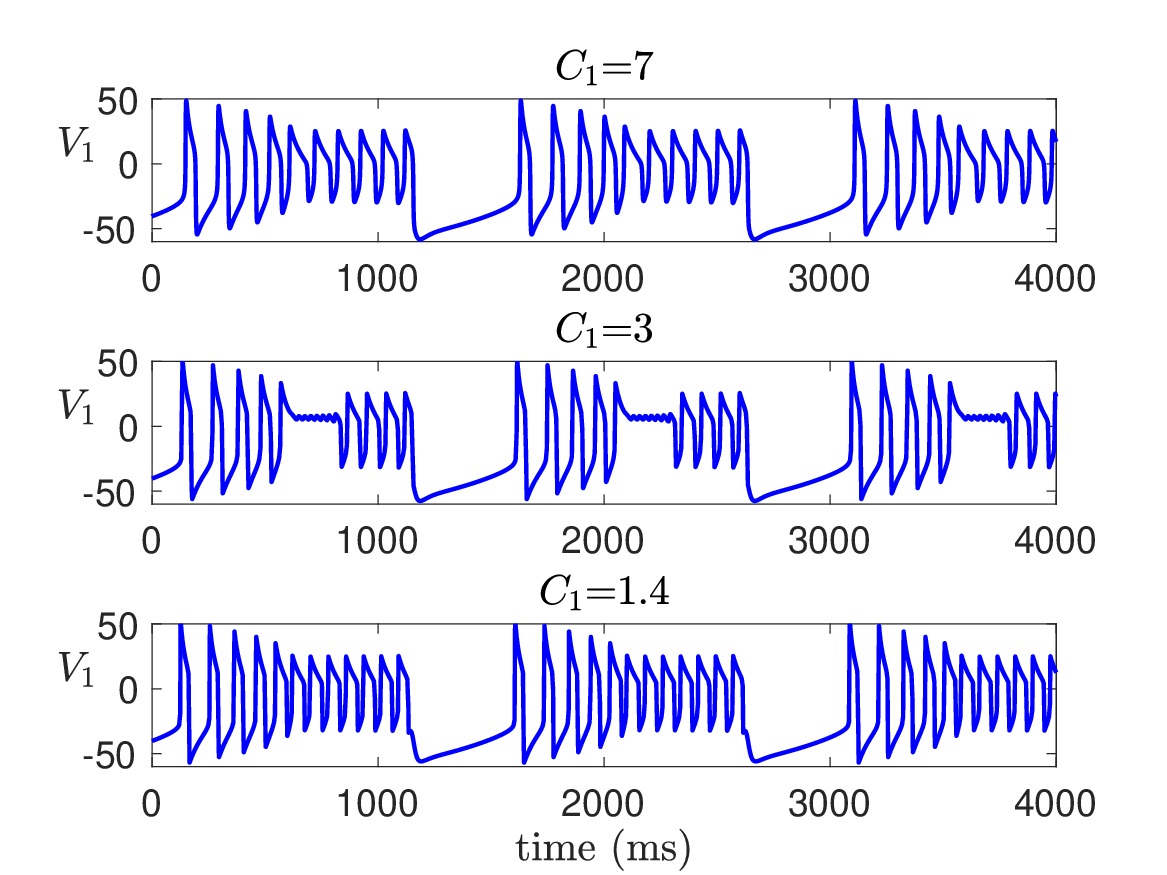}
        \end{tabular}
    \end{center}
\caption{Effect of varying $C_1$ (or $\varepsilon$) on time traces of solutions of \eqref{eq:main} for $g_{\rm syn}=4.3, \phi_2=0.001$ and other parameters given in Table \ref{tab:par}. (A) Increasing $C_1$ from its default value $8$ (equivalently, increasing $\varepsilon$ from $0.1$) preserves MMOs. The number of small oscillations increases with $C_1$. (B) Decreasing $C_1$ from the default value leads to a transition from MMOs ($C_1=8$) to non-MMOs ($C_1=7$) to MMOs again ($C_1=3$) to non-MMOs ($C_1=1.4$).}
\label{fig:gsyn_4p3_varyC1}
\end{figure*}

Recall when $\gsyn=4.3$, the only mechanism available for MMOs is the DHB mechanism. While there exist folded node singularities, they do not play any significant role in generating MMOs. In this subsection, we explore the effects of $\varepsilon$ and $\delta$ on the dynamics of the full system \eqref{eq:slow} by mainly examining their effects on the upper DHB around which SAOs are observed. As before, we will vary $\varepsilon$ by changing $C_1$ and vary $\delta$ by changing $\phi_2$. 
Recall that increasing (resp., decreasing) $C_1$ or $\varepsilon$ slows down (resp., speeds up) the fast variable $V_1$, whereas increasing (resp., decreasing) $\phi_2$ or $\delta$ speeds up (resp., slows down) the superslow variable $w_2$. Other (slow) variables are not affected. 

\begin{figure*}[!t]
    \begin{center}
        \begin{tabular}{@{}p{0.47\linewidth}@{\quad}p{0.47\linewidth}@{}}
            \subfigimg[width=\linewidth]{\bf{\small{(A)}}}{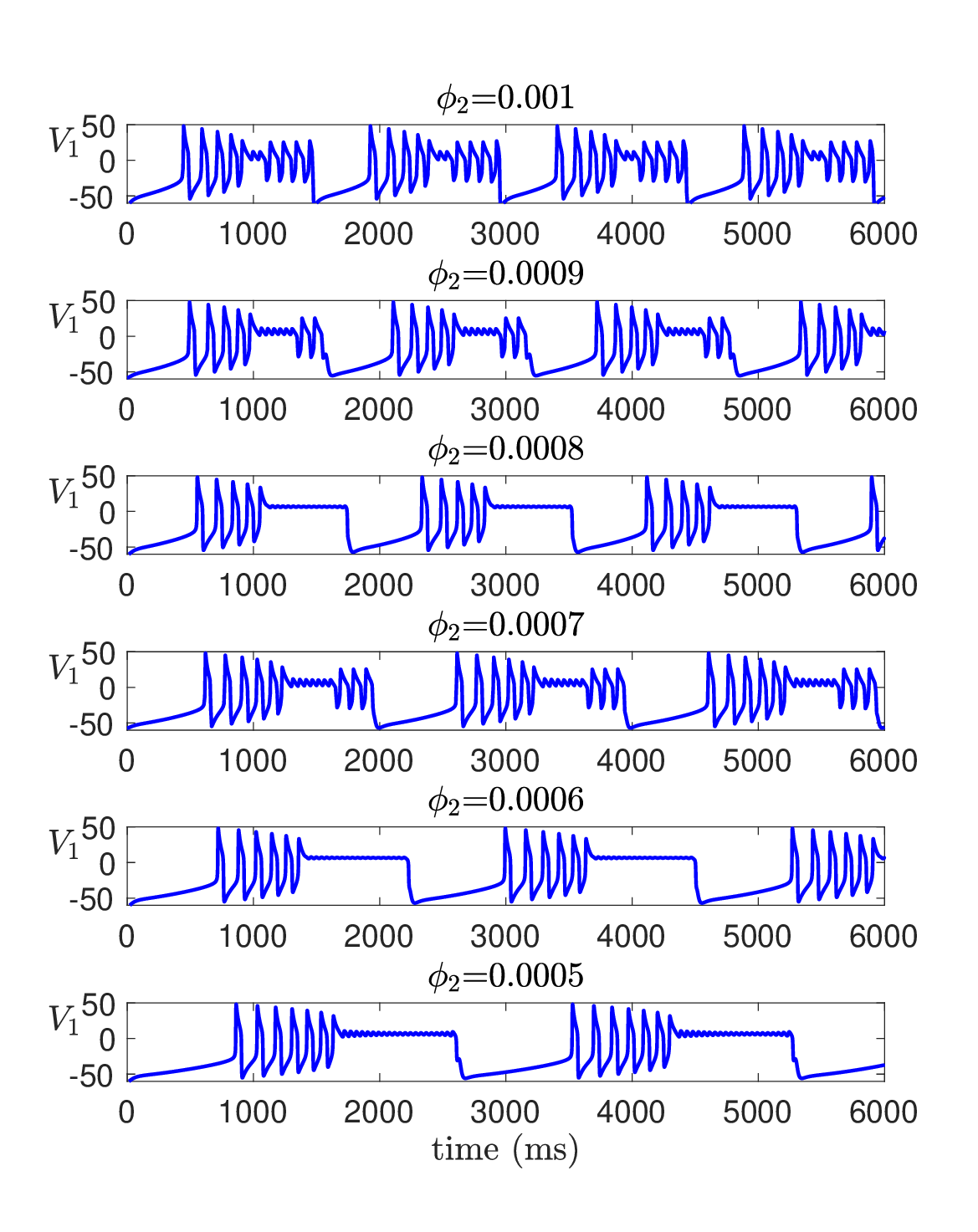}&
            \subfigimg[width=\linewidth]{\bf{\small{(B)}}}{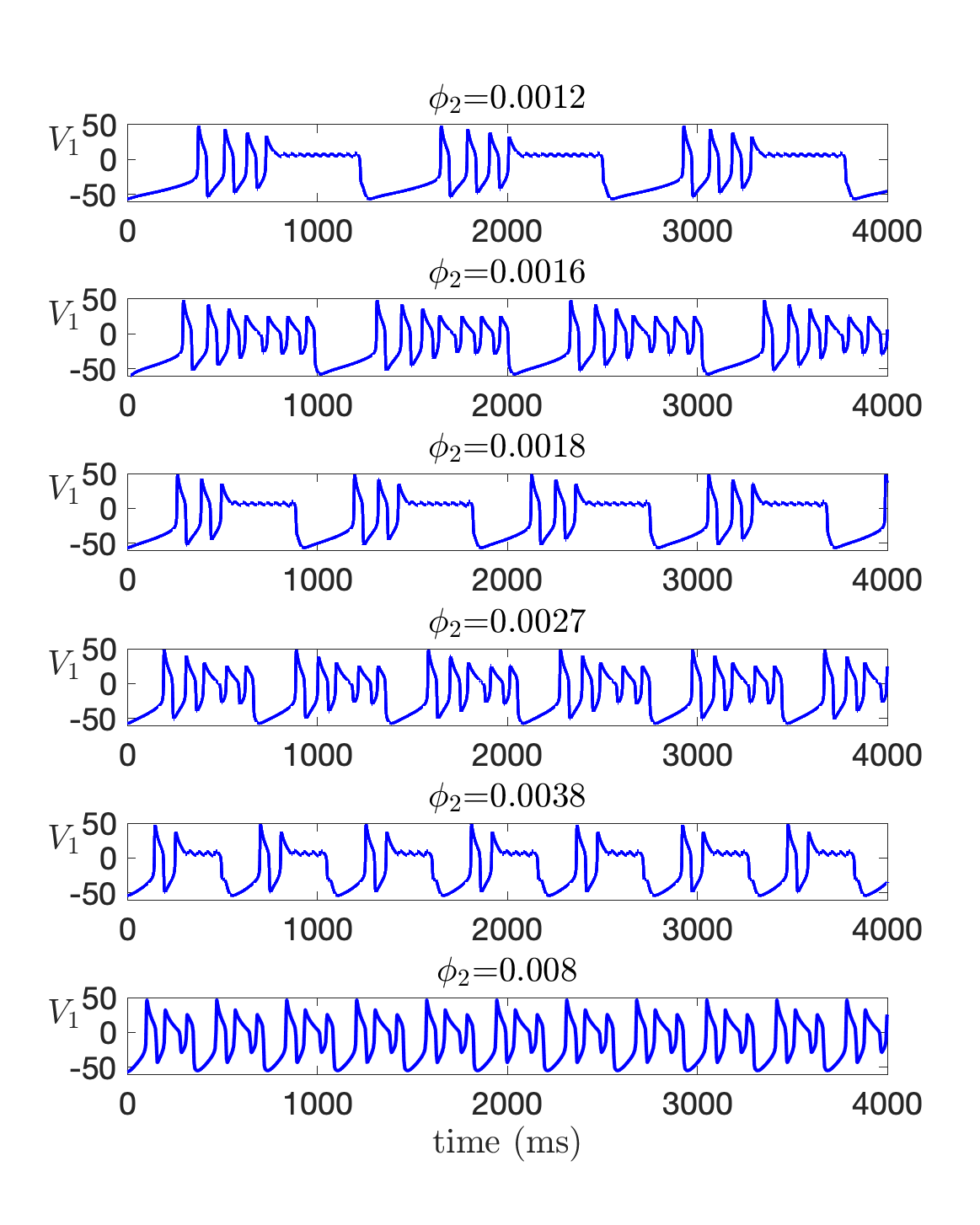}
        \end{tabular}
    \end{center}
\caption{Effect of varying $\phi_2$ (or $\delta$) on time traces of solutions of \eqref{eq:main} for $g_{\rm syn}=4.3, C_1=8$ and other parameters given in Table \ref{tab:par}. (A) Decreasing $\phi_2$ preserves MMOs, but the amplitude and the number of the small oscillations alternates between increasing and decreasing with $\phi_2$. 
(B) Increasing $\phi_2$ leads to a total of five transitions between MMOs and non-MMOs. 
\RED{MMOs exist for $\phi_2< 0.0016$ (i.e., $\delta \leq O(\varepsilon)$) and are completely lost for $\phi_2 > 0.007$ (i.e., $\delta \geq O(\varepsilon^{\frac{1}{3}}))$.}
}
\label{fig:gsyn_4p3-phi2}
\end{figure*}

Figure \ref{fig:c1-phi2}A summarizes the effects of $(C_1, \phi_2)$ on MMOs when $\gsyn=4.3$. Our findings suggest that MMOs with only the DHB mechanism are robust to changes that slow down either the fast variable or the superslow variable, but they are vulnerable to perturbations that speed up either timescale \RED{to a degree where $\delta$ is approximately greater than $O(\varepsilon)$}. 
Specifically, we observe the following:
\begin{itemize}
    \item For fixed $\delta=0.053$ at $\phi_2=0.001$, slowing down the fast variable $V_1$ by increasing $C_1$ from its default value $C_1=8$ (Figure \ref{fig:c1-phi2}A, vertical line above the red star) preserves the MMOs (also see Figure  \ref{fig:gsyn_4p3_varyC1}A, from top to bottom row). Moreover, we observe more characteristics of DHBs in the SAOs as $\varepsilon$ increases.
    \item For fixed $\delta$ at $\phi_2=0.001$, speeding up the fast variable $V_1$ via decreasing $C_1$ leads to a total of 3 transitions between MMOs and non-MMOs (see Figure \ref{fig:gsyn_4p3_varyC1}B). These transitions correspond to the 3 crossings between the vertical black line and the yellow/blue boundary in Figure \ref{fig:c1-phi2}.
    \item For fixed $\varepsilon=0.1$ at $C_1=8$, slowing down the superslow variable $w_2$ by decreasing $\phi_2$ from the default value $\phi_2=0.001$ (Figure \ref{fig:c1-phi2}A, red star) preserves MMOs. However, the number of small oscillations in the MMOs does not exhibit a simple monotonic increase or decrease, but rather shows an alternating pattern of increase and decrease as $\phi_2$ decreases. Additionally, the amplitudes of the small oscillations also display a similar non-monotonic behavior (see Figure \ref{fig:gsyn_4p3-phi2}A).
    \item For fixed 
    $\varepsilon$ at $C_1=8$, speeding up the superslow variable $w_2$ by increasing $\phi_2$ leads to a total of five transitions between MMOs and non-MMOs (see Figure \ref{fig:gsyn_4p3-phi2}B). These transitions correspond to the five crossings between the horizontal black line and the yellow/blue boundary in Figure \ref{fig:c1-phi2}A. 
\end{itemize}

Next, we discuss the above four scenarios separately. 

\begin{figure*}[!htp]
    \begin{center}
        \begin{tabular}{@{}p{0.33\linewidth}@{\quad}p{0.33\linewidth}p{0.33\linewidth}@{}}
            \subfigimg[width=\linewidth]{}{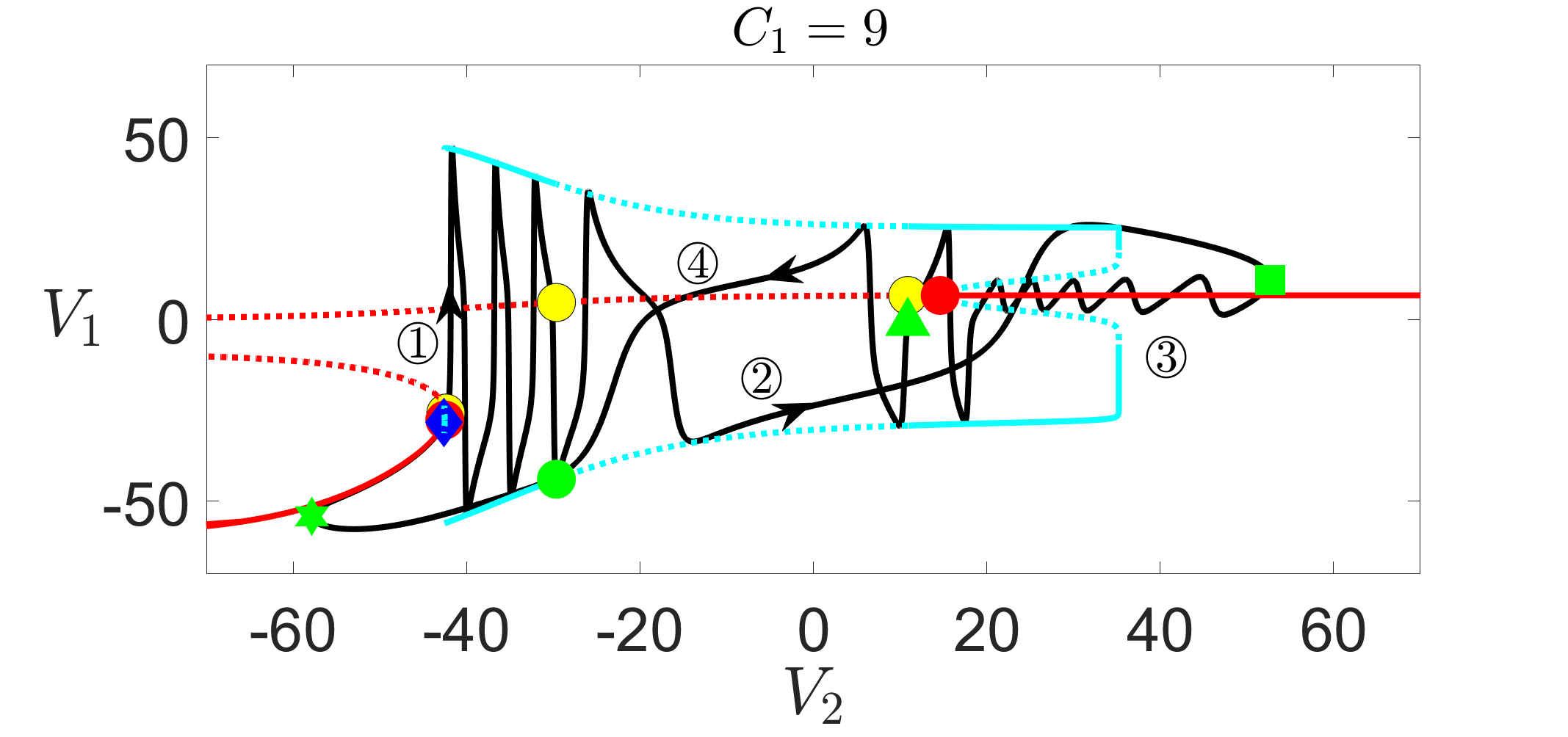}&
            \subfigimg[width=\linewidth]{}{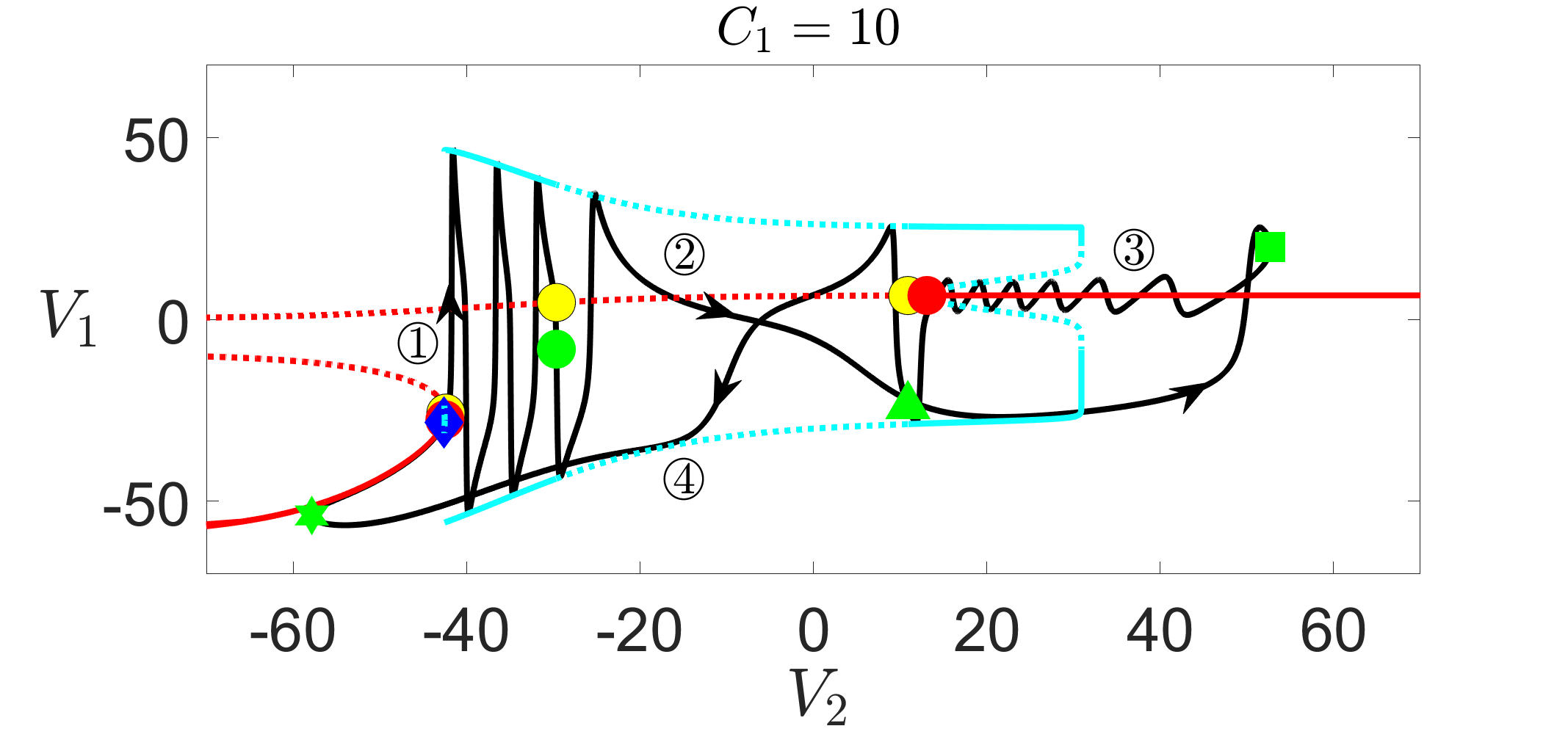}&
            \subfigimg[width=\linewidth]{}{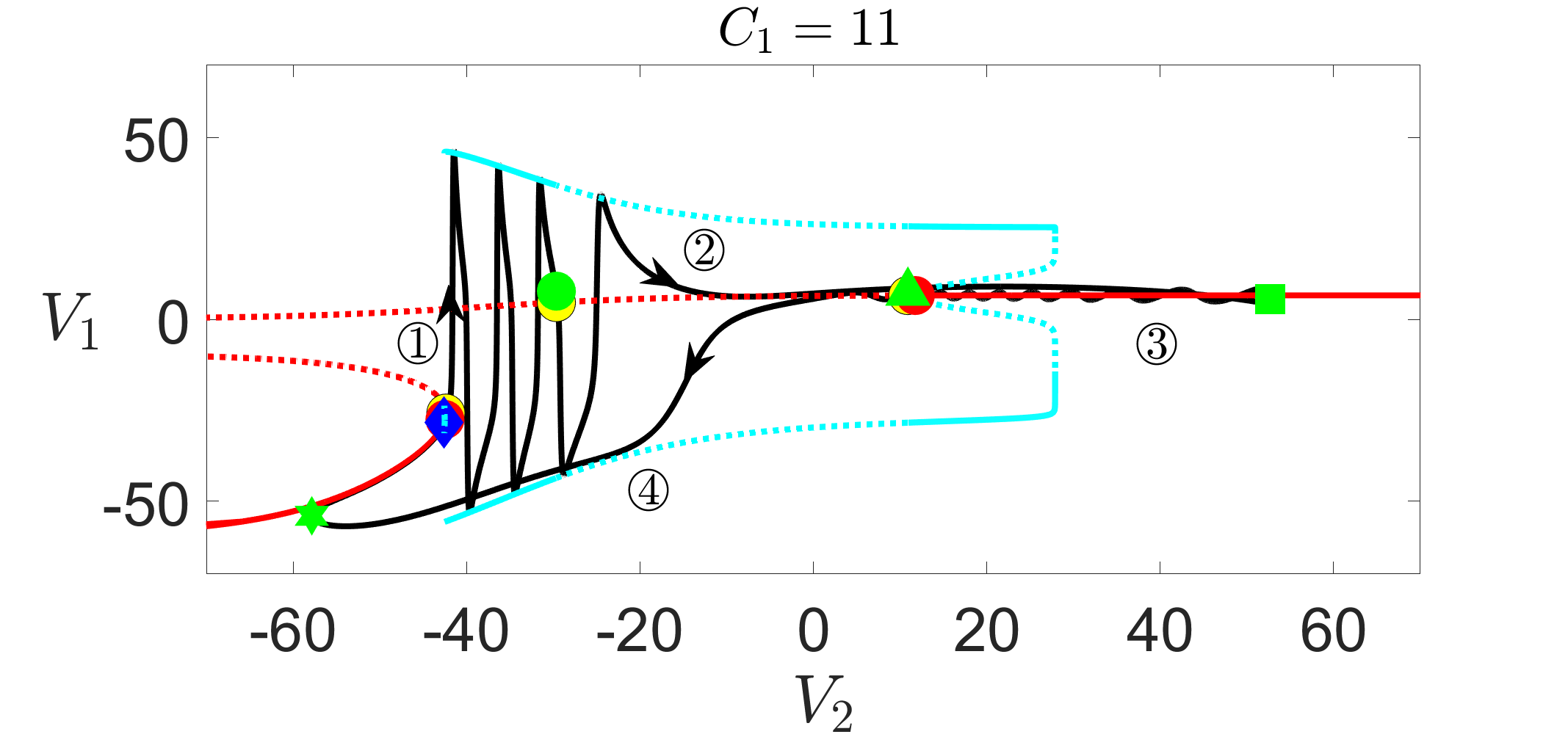}
        \end{tabular}
    \end{center}
\caption{Solutions of \eqref{eq:main} for $\gsyn=4.3$ and various values of $C_1$ and bifurcation diagrams of the slow layer problem \eqref{eq:slowlayer}, projected to $(V_2, V_1)$-space. From left to right, $C_1=9$ ($\epsilon=0.1125$), $C_1=10$ ($\epsilon=0.125$), $C_1=11$ ($\epsilon=0.1375$).  
Increasing $C_1$ moves the upper HB bifurcation (red circle) to lower $V_2$ values and hence increases the number of small oscillations surrounding $M_{ss}$. Color codings and symbols are the same as in Figure \ref{fig:Ms-V1V2-gsyn4p3}A.}
\label{fig:gsyn4p3-bif-largerC1}
\end{figure*}

\subsubsection{MMOs are robust to increasing $C_1$}

Increasing $C_1$ slows down the fast variable $V_1$ and hence moves the three timescale (1F, 2S, 1SS) problem closer to (3S, 1SS) separation. As a result, the critical manifold $\ms$ and the folded singularities become less relevant with increased $C_1$. Nonetheless, this does not affect the existence of MMOs, as they arise from the upper Hopf bifurcation. As $C_1$ is increased, the upper HB moves to smaller $V_2$ values (see Figure \ref{fig:two-par-4p3}A). This change allows the trajectory to travel a longer distance along the stable part of $\mss$ during phase \textcircled{3} and generate more SAOs, as shown in Figure \ref{fig:gsyn4p3-bif-largerC1}. \RED{To simplify the presentation, we omit singular orbits and focus only on $\mss$ and bifurcation diagrams that are essential for organizing SAOs in the full model, as demonstrated earlier. }

Recall that the small oscillations occurring during phase \textcircled{3} switch to large spikes upon crossing the inner unstable periodic orbit branch that is born at the upper HB (see Figure \ref{fig:Ms-V1V2-gsyn4p3}A).
Increasing $C_1$ causes the upper HB in $(V_2, V_1)$-space (Figure \ref{fig:gsyn4p3-bif-largerC1}, red circle) to move to the left and become further away from the maximum of $V_2$ (green square) which remains unchanged. As a result, the trajectory with larger $C_1$ begins small oscillations with decaying amplitude at a greater distance from the Hopf bifurcation point. This greater distance results in the trajectory crossing the inner unstable periodic orbit at a smaller value of $V_2$ (compare Figures \ref{fig:Ms-V1V2-gsyn4p3} and \ref{fig:gsyn4p3-bif-largerC1}). As $C_1$ is increased to $11$, the trajectory passes over the HB to $\mss^r$ and experiences a delay before jumping away (see Figure \ref{fig:gsyn4p3-bif-largerC1}, right panel). \RED{As explained earlier, the absence of small oscillations with growing amplitudes after the upper HB is due to the slow jump down of $V_2$ at the green triangle.}

\subsubsection{Decreasing $C_1$ leads to three MMOs/non-MMOs transitions}

\begin{figure*}[!htp]
    \begin{center}
        \begin{tabular}{@{}p{0.4\linewidth}@{\quad}p{0.4\linewidth}@{}}
            \subfigimg[width=\linewidth]{\bf{\small{(A)}}}{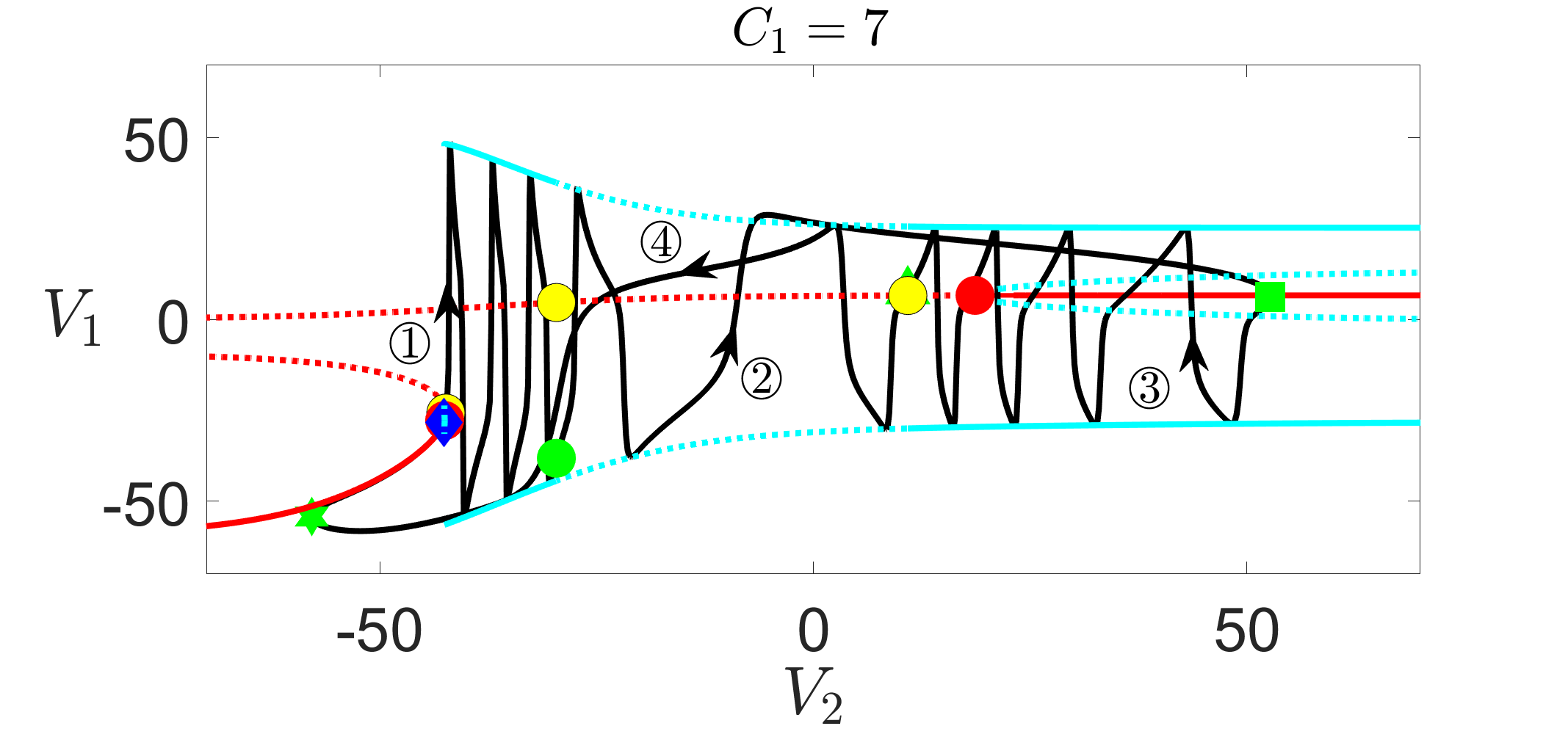}&
            \subfigimg[width=\linewidth]{\bf{\small{(B)}}}{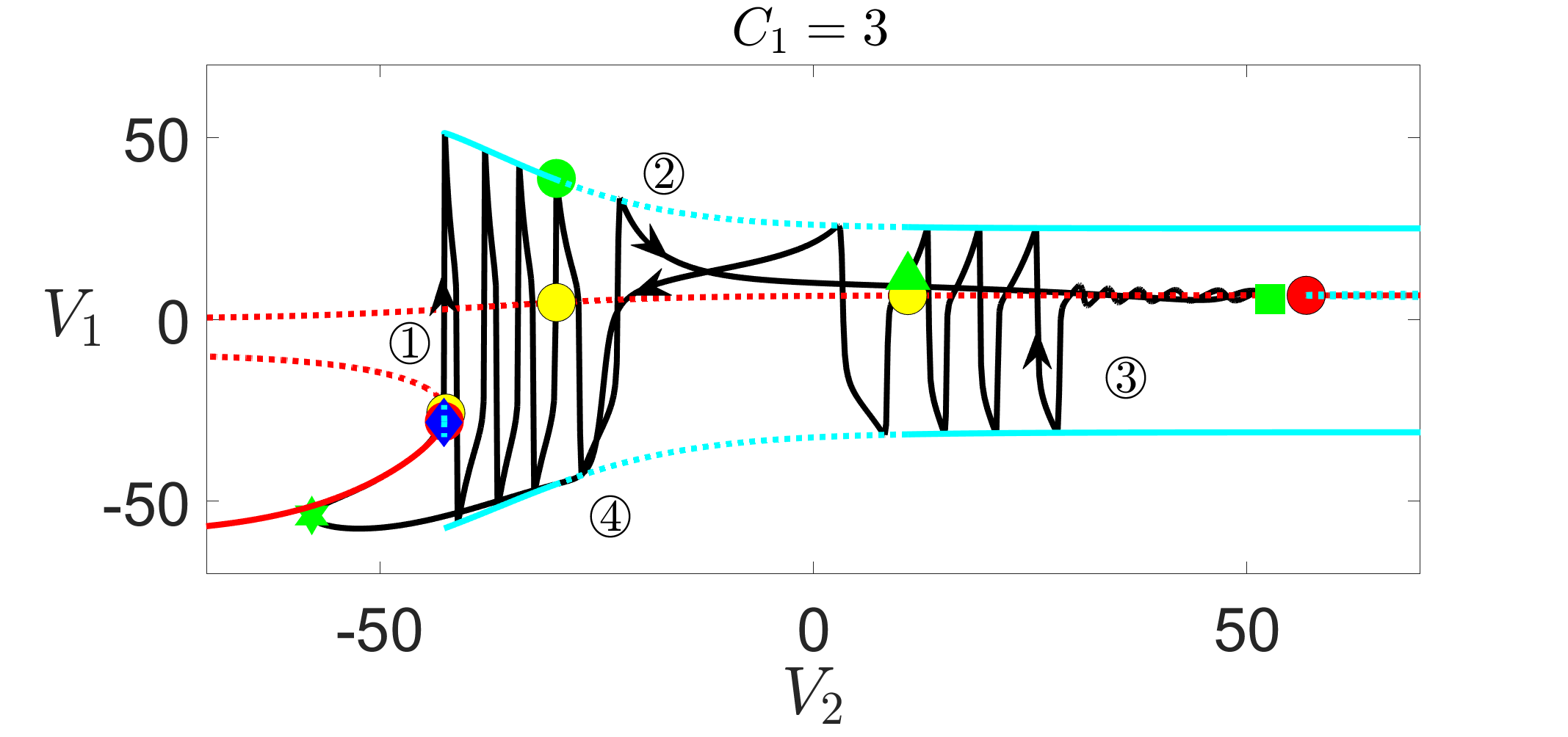}\\
            \subfigimg[width=\linewidth]{\bf{\small{(C)}}}{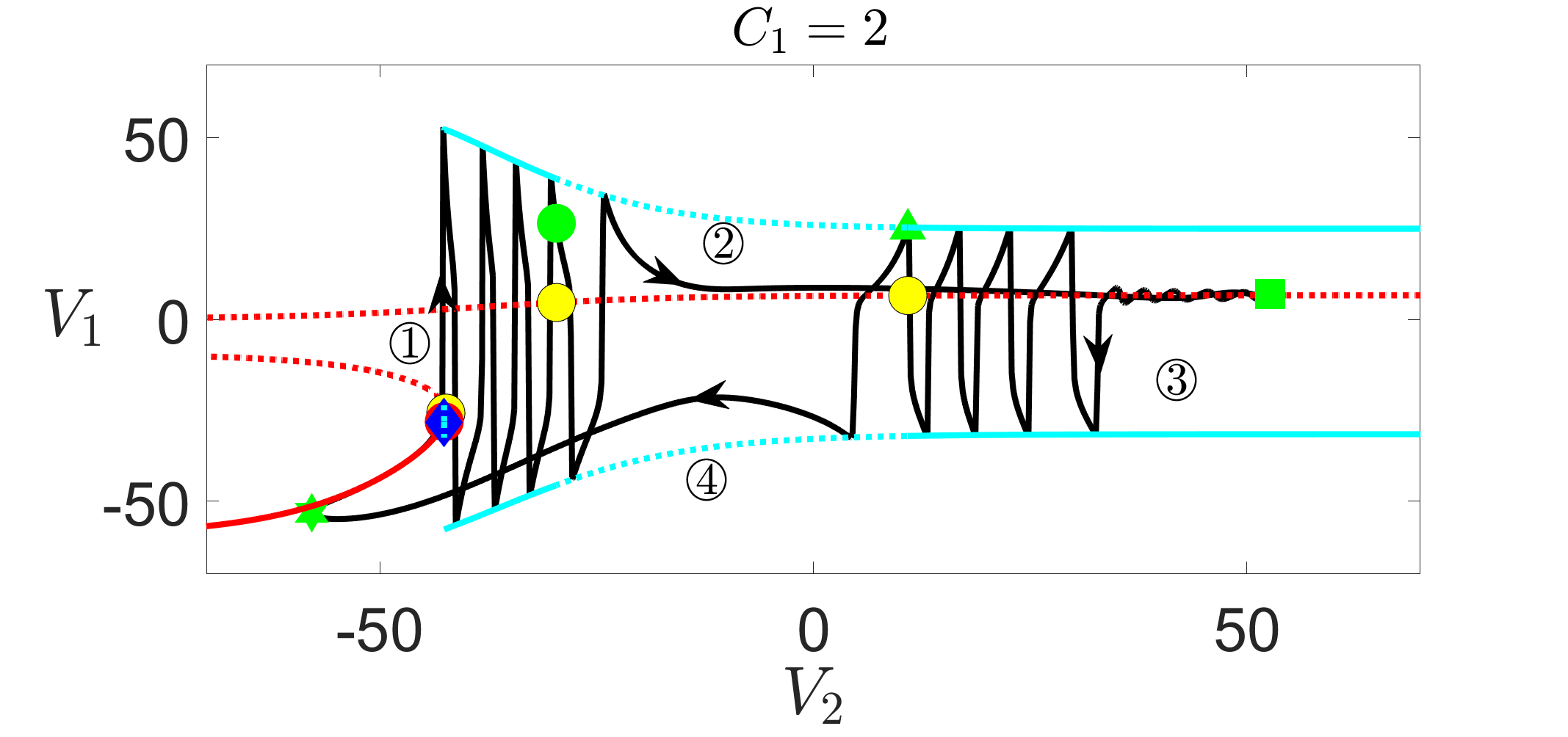}&
            \subfigimg[width=\linewidth]{\bf{\small{(D)}}}{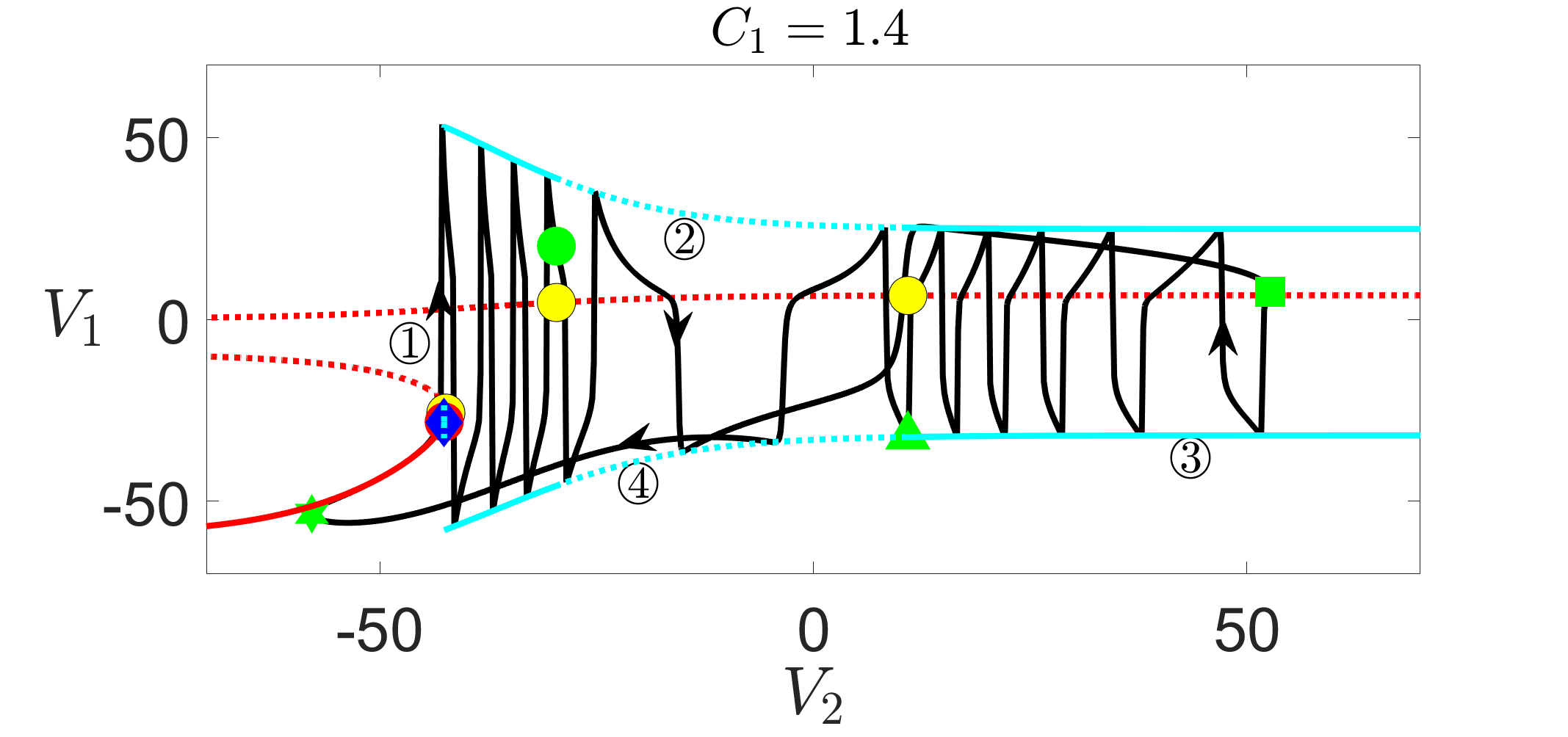}			
        \end{tabular}
    \end{center}
\caption{Projections of the full system solutions  and bifurcation diagrams of \eqref{eq:slowlayer} for $\gsyn=4.3$ and (A) $C_1=7$ ($\epsilon =0.0875 $),
(B) $C_1=3$ ($\epsilon =0.0375$), (C) $C_1=2$ ($\epsilon =0.025$), (D) $C_1=1.4$ ($\epsilon =0.0175$). The upper DHB (red circle) moves to larger $V_2$ values with decreased $C_1$ and eventually vanishes for $C_1<2.9$. Color codings and symbols have the same meaning as in Figure \ref{fig:Ms-V1V2-gsyn4p3}A.}
\label{fig:de-c1}
\end{figure*}

Contrary to the preservation of MMOs with increasing $C_1$, MMOs appear to be sensitive to the decrease of $C_1$. Specifically, speeding up the fast variable $V_1$ by decreasing $C_1$ leads to transitions from MMOs (e.g., at $C_1=8$) to non-MMOs (e.g., $C_1=7$), back to MMOs at $C_1\in (1.46, 4.2)$, and then to non-MMOs for $C_1<1.46$.

The initial decrease of $C_1$ (e.g., from $C_1=8$ to $C_1=7$, see Figure \ref{fig:de-c1}A) results in the loss of SAOs and thus a transition to non-MMOs due to the opposite effect of the mechanism discussed in the case of increasing $C_1$. Interestingly, a further decrease of $C_1$ to a range of $C_1\in(1.46, 4.2)$, which causes the upper HB to cross the green square and eventually vanish, results in the recovery of MMOs characterized by small oscillations with increasing amplitude (see Figure \ref{fig:de-c1}B and C). This is because, for $C_1\in (1.46, 4.2)$ and $V_2$ near the green square, \eqref{eq:slowlayer} exhibits either unstable periodic orbits with negligible amplitudes or saddle foci \RED{equilibrium} characterized by a negative real eigenvalue and complex eigenvalues with positive real parts. As a result, the SAOs during phase \textcircled{3} grow in amplitude as the trajectory spirals away from $\mss$. When $C_1$ is reduced to be below $1.46$, the voltage $V_1$ exhibits rapid spikes that occur immediately after the green square, \RED{consistent with the singular orbits shown in Figure \ref{fig:singular-orbit}B and C}.
As a result, there is no more MMOs (see Figure \ref{fig:de-c1}D, $C_1=1.4$).

\begin{figure*}[!t]
    \begin{center}
        \includegraphics[width=0.8\linewidth]{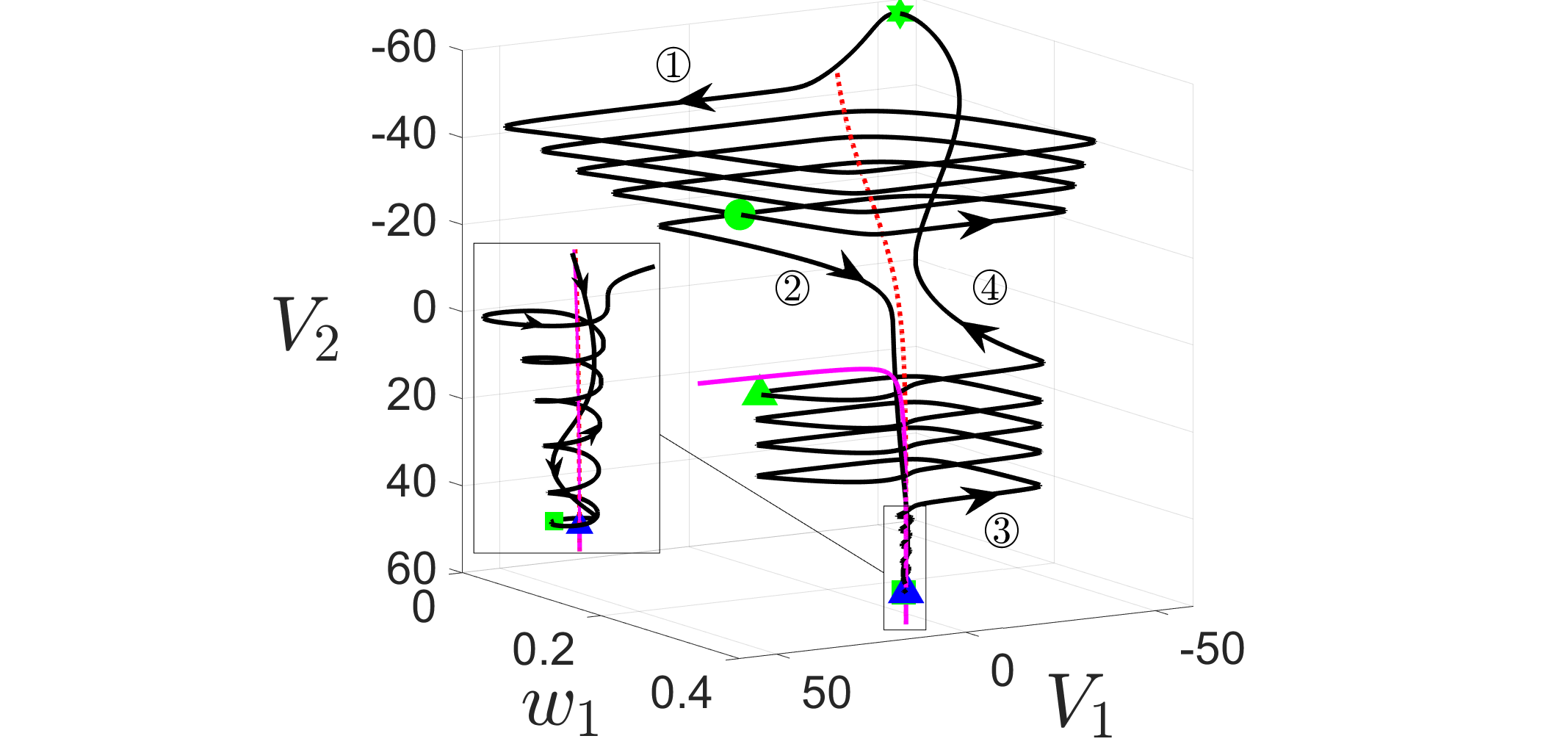}
    \end{center}
\caption{The solution (black curve) from Figure \ref{fig:de-c1}C projected onto $(V_1, w_1,V_2)$ space, together with $\mss$ (red curve). The blue triangle denotes a saddle-focus equilibrium on $\mss$ with maximum $V_2$, which has one negative real eigenvalue and a pair of complex eigenvalues with positive real parts ($0.0027\pm 0.31i$). A stable manifold branch associated with the saddle-focus at the blue triangle is denoted by the magenta curve. }
\label{fig:4p3-saddlefocus}
\end{figure*}

To better understand why the SAOs for $C_1\in (1.46, 4.2)$ grow in amplitude (see Figure \ref{fig:de-c1})B and C), we project the trajectory when $C_1=2$ onto $(V_1, w_1, V_2)$ space (see Figure \ref{fig:4p3-saddlefocus}). The blue triangle denotes a saddle focus of the $(V_1,w_1,V_2)$ subsystem near the green square. During phase \textcircled{2} from the green circle to the green square, the trajectory travels towards the saddle-focus (blue triangle) along its stable manifold (magenta curve) on the slow timescale. After phase \textcircled{3} begins at the green square, SAOs grow in amplitude as the trajectory moves upwards and spirals away from the equilibrium curve along its unstable manifold (not shown). Similar dynamical behaviors have also been observed near a subcritical Hopf-homoclinic bifurcation \cite{GW2000} and a singular Hopf bifurcation in two-timescale settings \citep{Desroches2012}.

\subsubsection{Decreasing $\phi_2$ preserves MMOs but causes non-monotonic effects on SAOs}

Decreasing \RED{the perturbation parameter $\phi_2$ (i.e., reducing $\delta$)} moves the system closer to the 3-slow/1-superslow splitting and manifests the DHB mechanism. This preserves the DHB-induced MMOs as expected. However, a non-monotonic effect on the small-amplitude oscillations is also observed, as shown in Figure \ref{fig:gsyn_4p3-phi2}A, where the amplitude and the number of SAOs exhibit an alternation between increase and decrease. 

Smaller perturbations should cause the solution trajectory to follow $\mss$ more closely and at a slower rate. Intuitively, one may expect that this leads to an increase in the number of SAOs and a decrease in their amplitudes. Indeed, we observe such changes as $\phi_2$ decreases from $0.001$ to $0.0008$, as shown in Figure \ref{fig:gsyn_4p3-phi2}A (top three rows) and also in Figure \ref{fig:4p3-phi2-p8-p9} as we discussed before. However, to our surprise, we find that for $\phi_2=0.0007$, the MMOs exhibit less SAOs with larger amplitudes than the SAOs at $\phi_2=0.0008$ (see Figure \ref{fig:gsyn_4p3-phi2}A, the third and the fourth row). The number of SAOs increases and their amplitudes decrease again as $\phi_2$ decreases from $0.0007$ to $0.0006$. This alternating pattern of changes in the number and amplitude of SAOs repeats as $\phi_2$ continues to decrease (see Figure \ref{fig:gsyn_4p3-phi2}A). 

Our analysis reveals that as $\phi_2$ decreases, there will be more SAOs with smaller amplitudes if no additional big (full) spike is generated. However, if an additional full spike is gained during the process of decreasing $\phi_2$, the changes to the SAOs will be reversed; that is, there will be fewer SAOs and they will exhibit larger amplitudes. This is because the additional spike before SAOs can push the trajectory away from $\mss$ at the beginning of phase \textcircled{3}, leading to fewer SAOs with larger amplitudes. \RED{Hence, the amplitude and number of SAOs in the full system are not only determined by the size of the perturbation but also by how the flow approaches $\mss^a$ during phase \textcircled{2}. As previously discussed, infinitely many singular orbit segments can be constructed during this phase. This leads to different ways for the full trajectory to reach $\mss^a$ under varying perturbations.} In a sense, this alternating pattern of changes in SAOs occurs due to a spike-adding like mechanism. We refer the readers to Appendix \ref{app:de-phi2-4p3} for a more detailed discussion on why decreasing $\phi_2$ causes non-monotonic effects on SAOs.

\subsubsection{Increasing $\phi_2$ leads to five MMOs/non-MMOs transitions}

When $\phi_2$ increases, the superslow variable $w_2$ speeds up, moving the system closer to the 1-fast/3-slow splitting and making the DHB mechanism less relevant. Since there is no canard mechanism, MMOs should be eliminated for $\phi_2$ large enough. Indeed, we observe a total of five transitions between MMOs and non-MMOs as $\phi_2$ increases, and eventually, MMOs are lost for $\phi_2>0.007$ \RED{(i.e., $\delta\geq O(\varepsilon^{\frac{1}{3}}))$}.
The mechanism driving these MMOs/non-MMOs transitions over the increase of $\phi_2$ is similar to the mechanism underlying the non-monotonic effects on SAOs when $\phi_2$ is decreased. 
Specifically, if no additional big (full) spike before phase \textcircled{3} is lost with the increase of $\phi_2$, there will be fewer SAOs with larger amplitudes or no MMOs as one would naturally anticipate (e.g., when $\phi_2$ increases from $0.0012$ to $0.0016$, see Figure \ref{fig:gsyn_4p3-phi2}B, top two rows). In contrast, if one full spike is lost during the process of increasing $\phi_2$, changes to the SAOs will be reversed such that there will be more SAOs with smaller amplitude (e.g., when $\phi_2$ increases from $0.001$ to $0.0012$, see the top row in Figure \ref{fig:gsyn_4p3-phi2}A and B). Eventually, MMOs will be completely lost when $\phi_2\geq 0.008$ 
for which the HB is no longer relevant. 

For a more detailed discussion, we refer the readers to Appendix \ref{app:in-phi2-4p3}. 

\section{Analysis of MMOs when $g_{\rm syn}=4.4\,\rm mS/cm^2$}\label{4p4}

This section explores MMOs that occur when an upper CDH singularity is present. In this scenario, we show in subsection \ref{sec:soln-4p4} that both canard and DHB mechanisms coexist and interact to produce MMOs that exhibit significant robustness to timescale variations as shown in Figure \ref{fig:c1-phi2}B. We explain the robust occurrence of MMOs in subsection \ref{sec:effect-timescale-4p4} and show that the two MMO mechanisms can be modulated by adjusting $\epsilon$ and $\delta$. Specifically, increasing $\epsilon$ manifests the DHB-like characteristics, while an increase in $\delta$ leads to dominance of the canard mechanism. 

\begin{figure*}[!t] 
    \begin{center}
        \includegraphics[width=0.8\linewidth]{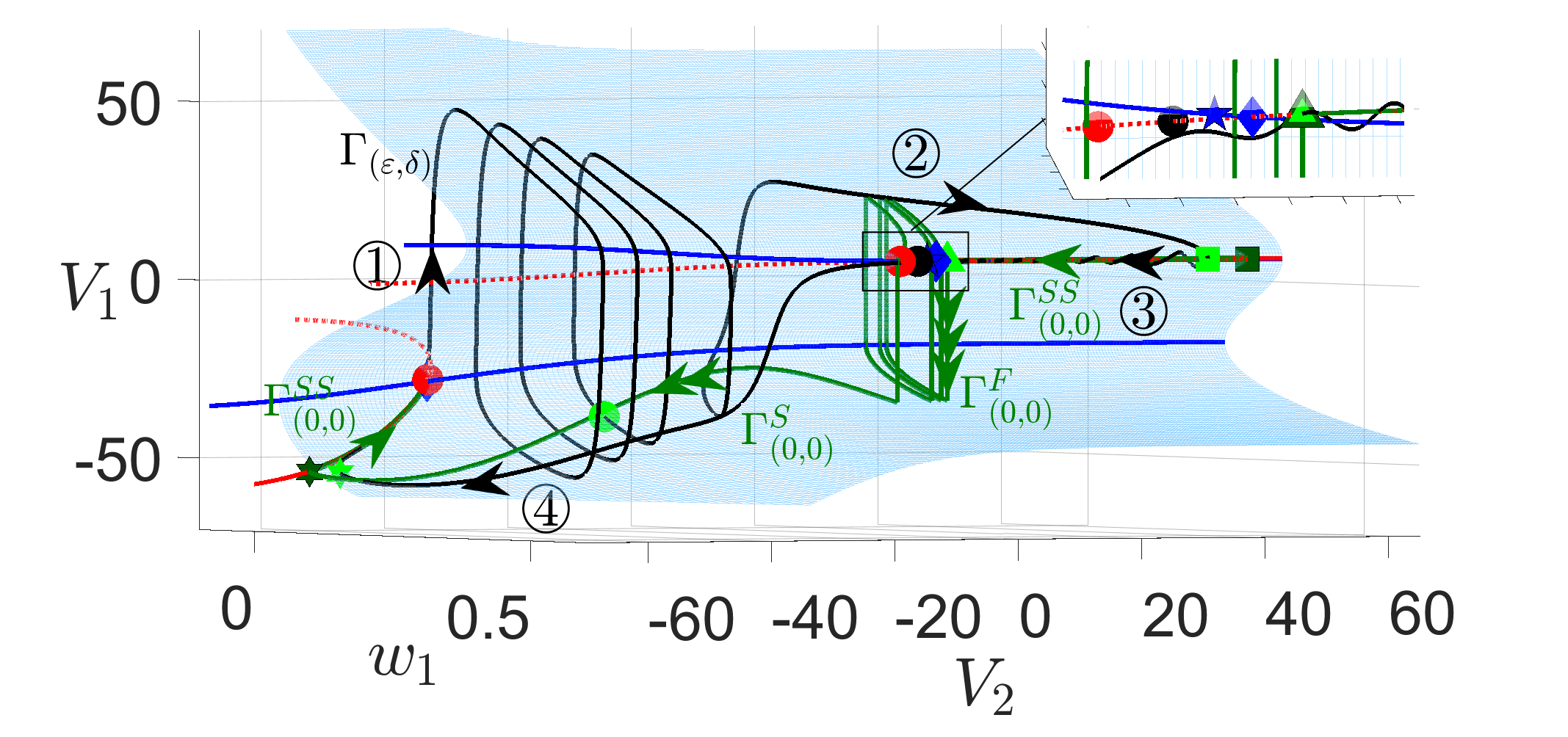}
    \end{center}
\caption{Projection of a solution trajectory (black curve) of system \eqref{eq:main} for $g_{\rm syn}=4.4$ from Figure \ref{fig:simu}B to $(V_1,w_1,V_2)$-space. Also shown are \RED{parts of the singular orbit $\Gamma^{F}_{(0,0)} \cup \Gamma^{S}_{(0,0)} \cup \Gamma^{SS}_{(0,0)}$ from Figure \ref{fig:singular-orbit}D (green curve),} the critical manifold $M_S$ (blue surface), folds $\ls$ of $\ms$ (blue curves) and the superslow manifold $M_{SS}$ (red curves). The solid (resp. dashed) red curves represent the attracting (resp. repelling) branches of $\mss$. 
Green symbols mark the transitions between slow and superslow pieces of the $(V_2,w_2)$ oscillations as in Figure \ref{fig:transition}. The blue diamonds denote CDH singularities - the intersection of $\ls$ and $\mss$. The upper CDH is a folded node and the lower CDH (almost overlapping with the lower DHB denoted by a red circle) is a folded focus. The upper CDH point is $O(\delta)$ close to the folded saddle-node singularity $\fsn^1$ (blue star) and $O(\varepsilon)$ close to the upper DHB (red circle). The black circle denotes the isolated ordinary singularity of the full system, whose type is a saddle-focus.  } 
\label{fig:Ms-gsyn4p4}
\end{figure*}

\subsection{Relation of the trajectory to $\ms$ and $\mss$}\label{sec:soln-4p4}

\begin{figure*}[!htp]
\begin{center}
\begin{tabular}{@{}p{0.48\linewidth}@{\quad}p{0.48\linewidth}@{}}
\subfigimg[width=\linewidth]{\bf{\small{}}}{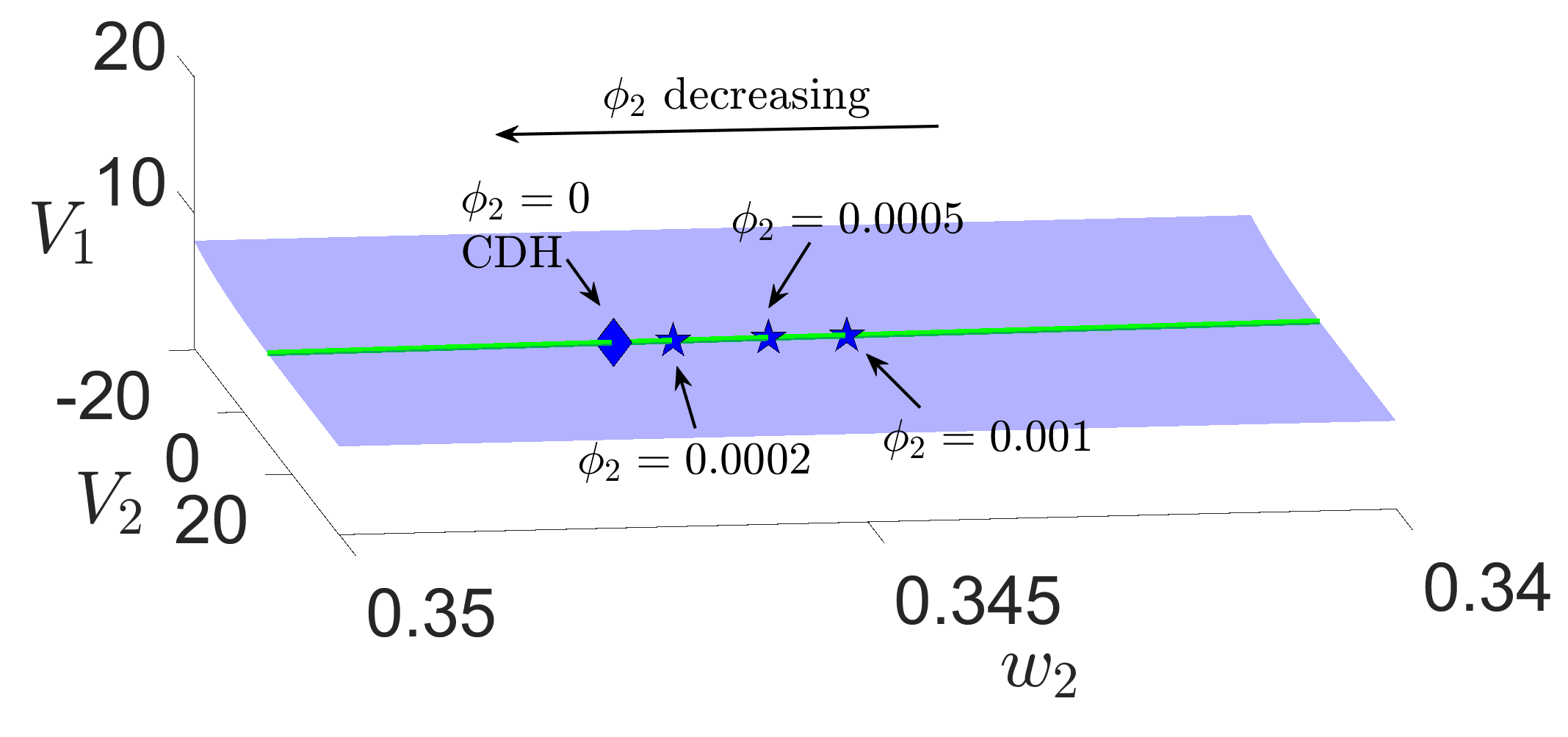}&
\subfigimg[width=\linewidth]{\bf{\small{}}}{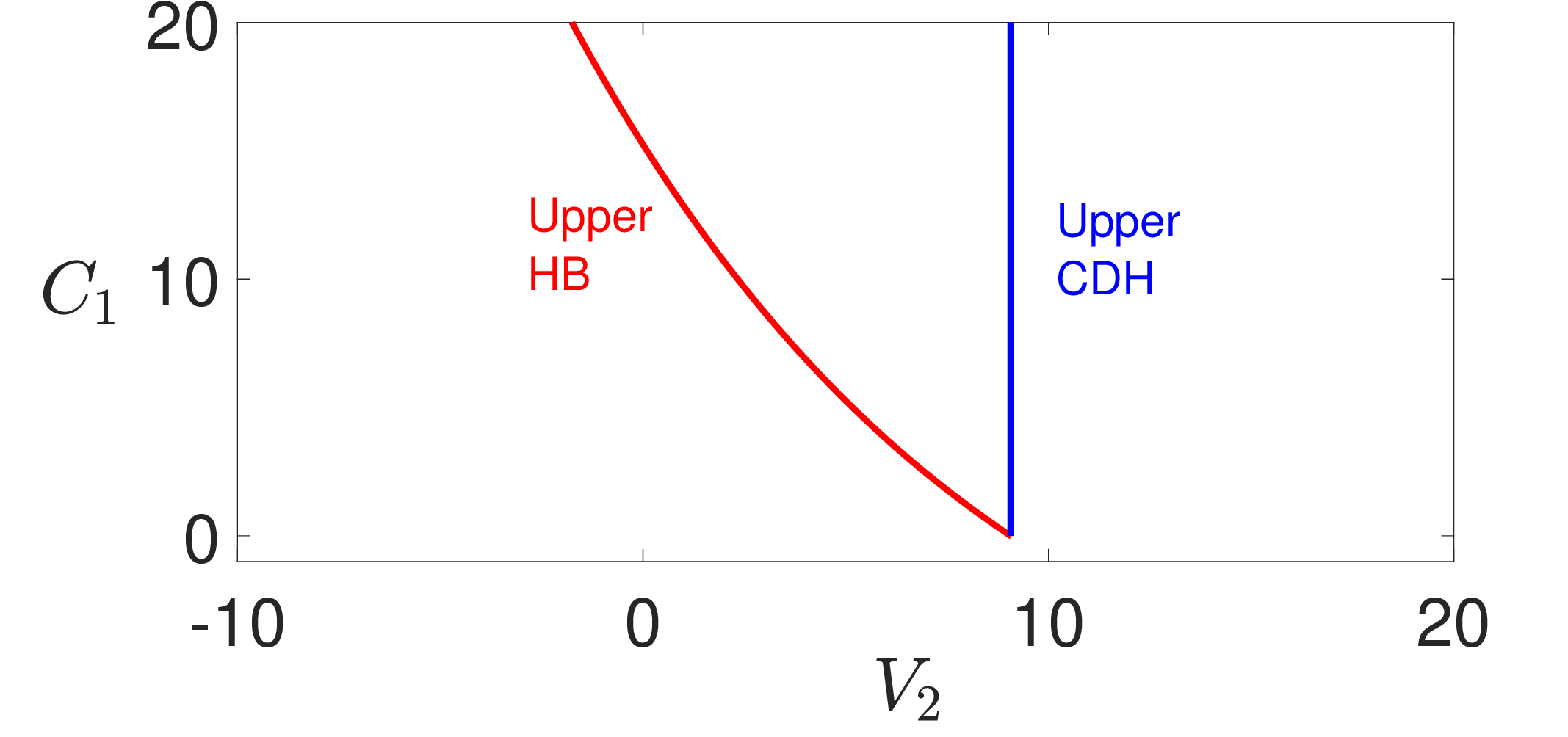}
\end{tabular}
\end{center}
\caption{The upper $\rm CDH$ represents the interaction of canard and DHB mechanisms when $g_{\rm syn}=4.4$. (Left) The upper $\fsn^1$ singularity (blue star) converges to the upper $\rm CDH$ (blue diamond) as $\phi_2\to 0$ (or equivalently $\delta\to 0$). (Right) As $C_1$ (or $\epsilon$) decreases, the upper HB (red curve) moves to larger $V_2$ values and approaches the upper $\rm CDH$ (blue curve) in the singular limit $\varepsilon\to 0$.}
\label{fig:effects-eps-hb-4p4}
\end{figure*}

The solution of \eqref{eq:main} for $g_{\rm syn}=4.4$ in \RED{Figure \ref{fig:simu}} is projected onto the $(V_1, w_1,V_2)-$space, together with \RED{the singular orbit from Figure \ref{fig:singular-orbit}D (green curve),} critical manifold $\ms$ (blue surface), fold $\ls$ (blue curve) and the superslow manifold $\mss$ (red curves) (see Figure \ref{fig:Ms-gsyn4p4}). As the coupling strength $\gsyn$ increases from $4.3$ to $4.4$, two new features regarding the upper $\mss$ emerge. First, the stability of the upper $\mss$ changes at a fold point $L^1_{ss}$ (the green triangle) rather than a DHB. Specifically, as $V_2$ increases, the equilibrium along the upper $\mss$ switches from unstable saddle-focus (characterized by one positive real eigenvalue and a pair of complex eigenvalues whose real parts are negative) to stable focus (with one negative real eigenvalue and a pair of complex eigenvalues whose real parts are negative). Second and more importantly, the upper fold  $\ls$ now intersects the upper branch of $\mss$ at a CDH singularity (the blue diamond), which is a folded node at the default parameter values given in Table \ref{tab:par}. 
As proved in subsection \ref{subsec:interaction}, this CDH point is located $O(\delta)$ close to a folded saddle-node singularity FSN$^1$ (upper blue star) and $O(\varepsilon)$ close to a HB (upper red circle). This is further confirmed in Figure \ref{fig:effects-eps-hb-4p4}, which shows that the upper $\fsn^1$ and HB point converge to the same CDH point on the upper $\ls$ in the double singular limit $(\varepsilon, \delta) \to (0,0)$.

Throughout the remainder of this section, we concentrate on elucidating the emergence and robustness of small amplitude oscillations (SAOs) in the vicinity of the upper CDH (see Figure \ref{fig:Ms-gsyn4p4}). \RED{ During phase \textcircled{3}, the singular orbit (green) traces $\mss^a$ and jumps down at the fold point.  The full trajectory $\Gamma_{(\varepsilon,\delta)}$ does not immediately jump at the fold. Instead, there is a delay before $\Gamma_{(\varepsilon,\delta)}$ jumps to the lower attracting manifold $\ms^L$. Small amplitude oscillations are observed as the trajectory passes through the neighborhood of the canard and DHB points. For smaller $\varepsilon$ or $\delta$ perturbations (e.g., top row in Figure \ref{fig:4p4_increasing_c1} and Figure \ref{fig:4p4_increasing_phi2}), the delay is more substantial, but the small oscillations are very small and difficult to observe due to the stronger attraction to $\mss^a$ during phase \textcircled{3}.} We omit explaining the solution dynamics of \eqref{eq:main} for $g_{\rm syn}=4.4$ during other phases as they are similar to those observed when $g_{\rm syn}=4.3$. In particular, the absence of SAOs near the lower CDH point is due to the same mechanism as discussed in subsection \ref{sec:noSTOs-CDH}.  

At default parameters, $\delta=O(\varepsilon)$. The emergence of SAOs for $g_{\rm syn}=4.4$ is governed by both the canard dynamics due to folded node singularities and the slow passage effects associated with the DHB on the upper fold $\ls$, as seen in examples considered in \citep{Vo2013,Letson2017}. To understand why folded node singularities play an important role in the occurrence of SAOs, we view the trajectory and folded singularities projected onto $(V_1, V_2,w_2)$-space, as illustrated in Figure \ref{fig:proj-v1v2w2-gsyn4p4}, as well as draw the funnel volume associated with the folded node singularity curve on the upper fold (see Figure \ref{fig:funnel-gsyn4p4}). On the other hand, to grasp the importance of the DHB towards the generation of SAOs, we examine the trajectory on the $(V_1,V_2)$ projection, which includes the periodic orbit branches born at the upper Hopf bifurcation (see Figure \ref{fig:Mss-gsyn4p4}). \RED{For simplicity, we omit depicting singular orbits in the following figures and focus only on $\ms$, $\mss$ and their relations to the full trajectory.}

\subsection{MMOs are organized by both canard and DHB mechanisms}

\subsubsection{Canard Dynamics} 
\begin{figure}[!htp] 
    \centering
   \includegraphics[width=\linewidth]{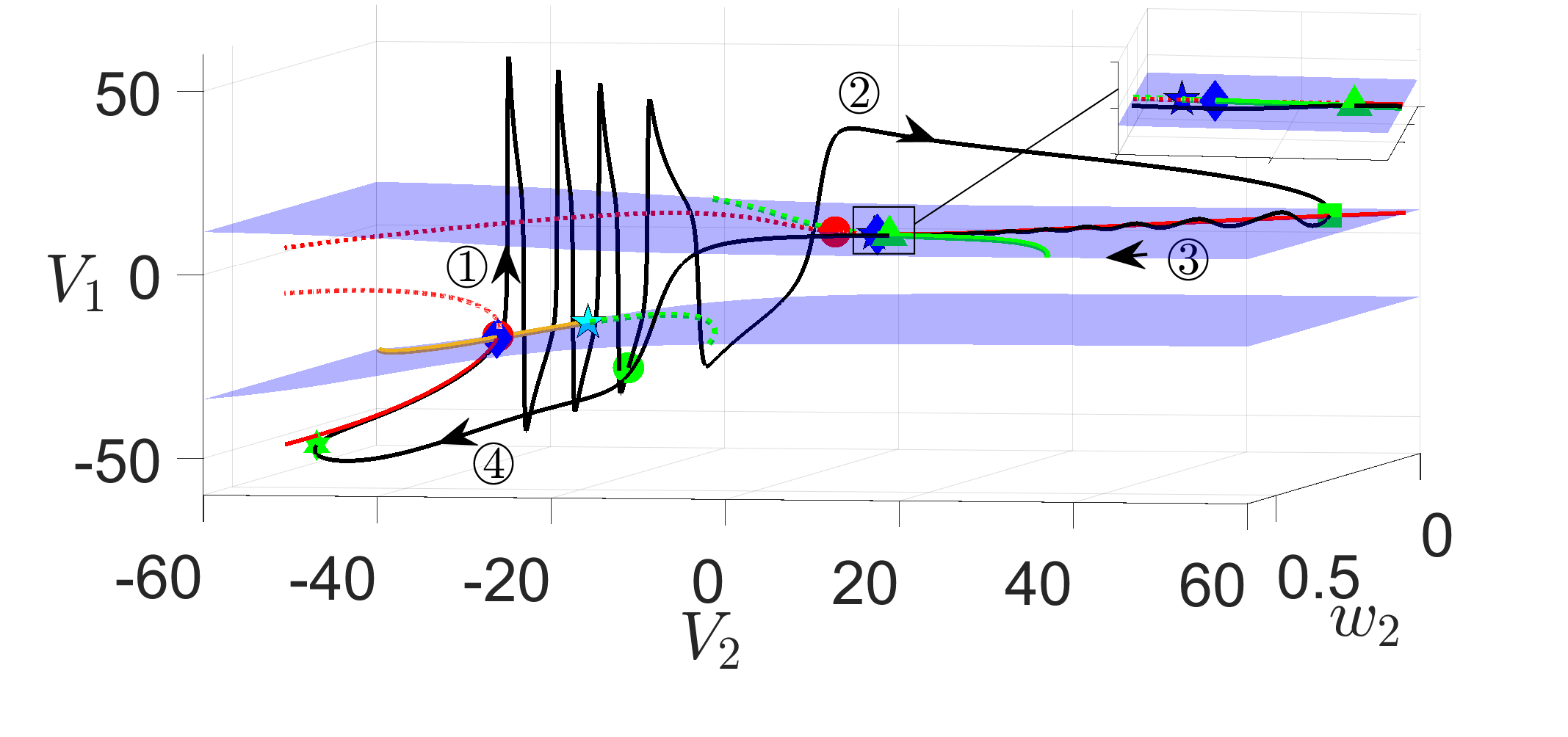}
    \caption{Projection of the solution from Figure \ref{fig:Ms-gsyn4p4}
    onto ($w_2,V_2,V_1$)-space. Also shown are the projections of the fold surface $\ls$ (blue surfaces), superslow manifold $\mss$ and folded singularities \RED{$\mathcal{M}$} (curves of green and yellow). Other color coding and symbols have the same meaning as in Figures \ref{fig:weak-eigen-direc} and \ref{fig:Ms-gsyn4p4}.}
    \label{fig:proj-v1v2w2-gsyn4p4}
\end{figure}

Figure \ref{fig:proj-v1v2w2-gsyn4p4} shows the solution trajectory of \eqref{eq:main} projected onto $(V_1, V_2,w_2)$-space when $g_{\rm syn}=4.4$. The upper folded singularity curve \RED{$\mathcal{M}$} comprises folded singularities of two types: folded nodes (solid green) and folded saddles (dashed green). The blue star denotes the saddle-node bifurcation $\fsn^1$, which occurs $O(\delta)$ away from the upper CDH at the blue diamond. Different from the previous case ($\gsyn=4.3$) where \RED{folded singularities} did not contribute to the SAO dynamics, the current solution trajectory crosses the upper \RED{$\mathcal{M}$} at a folded node, which can play a critical role in organizing the MMOs.

\begin{figure}[!htp]
    \begin{center}
         \includegraphics[width=0.9\linewidth]{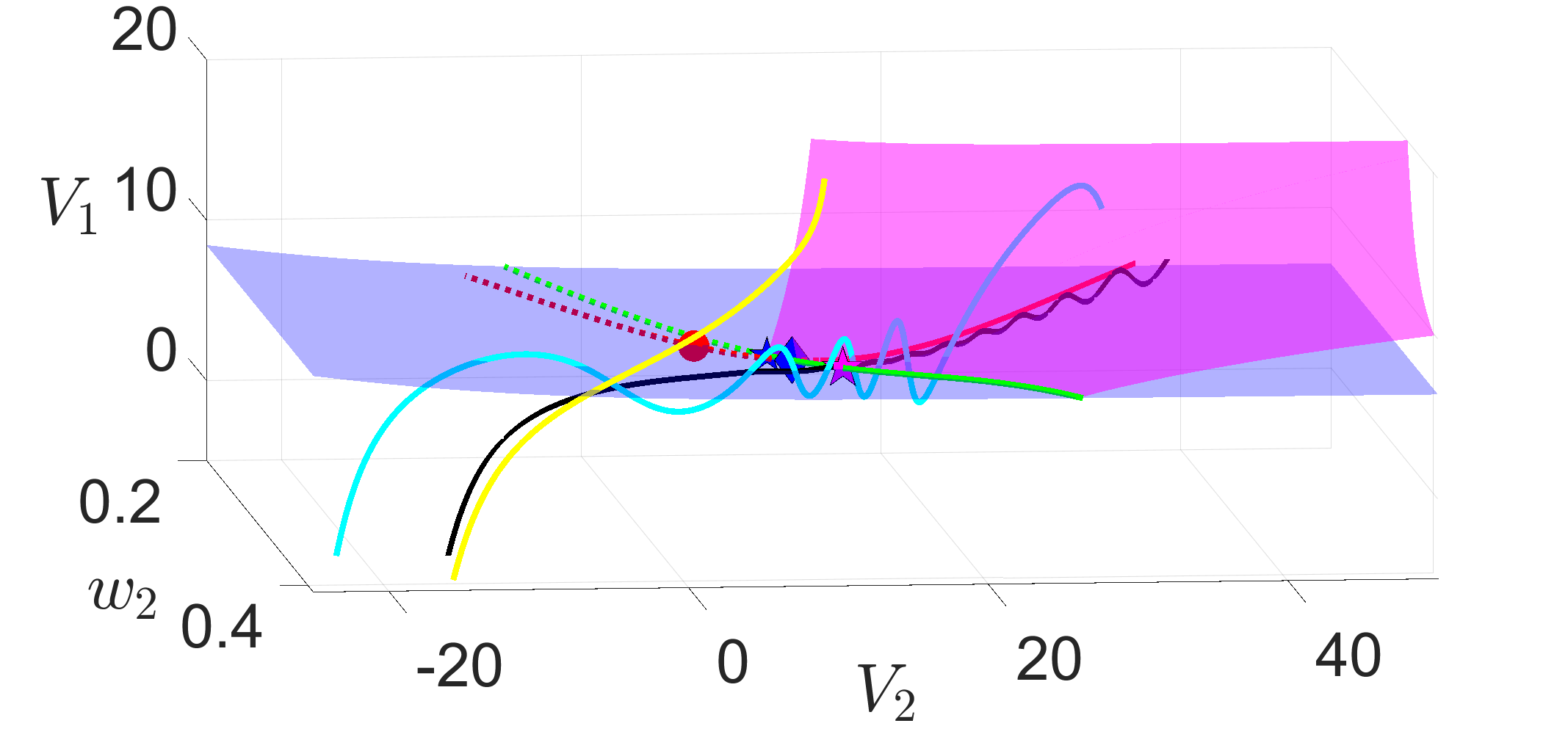}\\
            \includegraphics[width=0.9\linewidth]{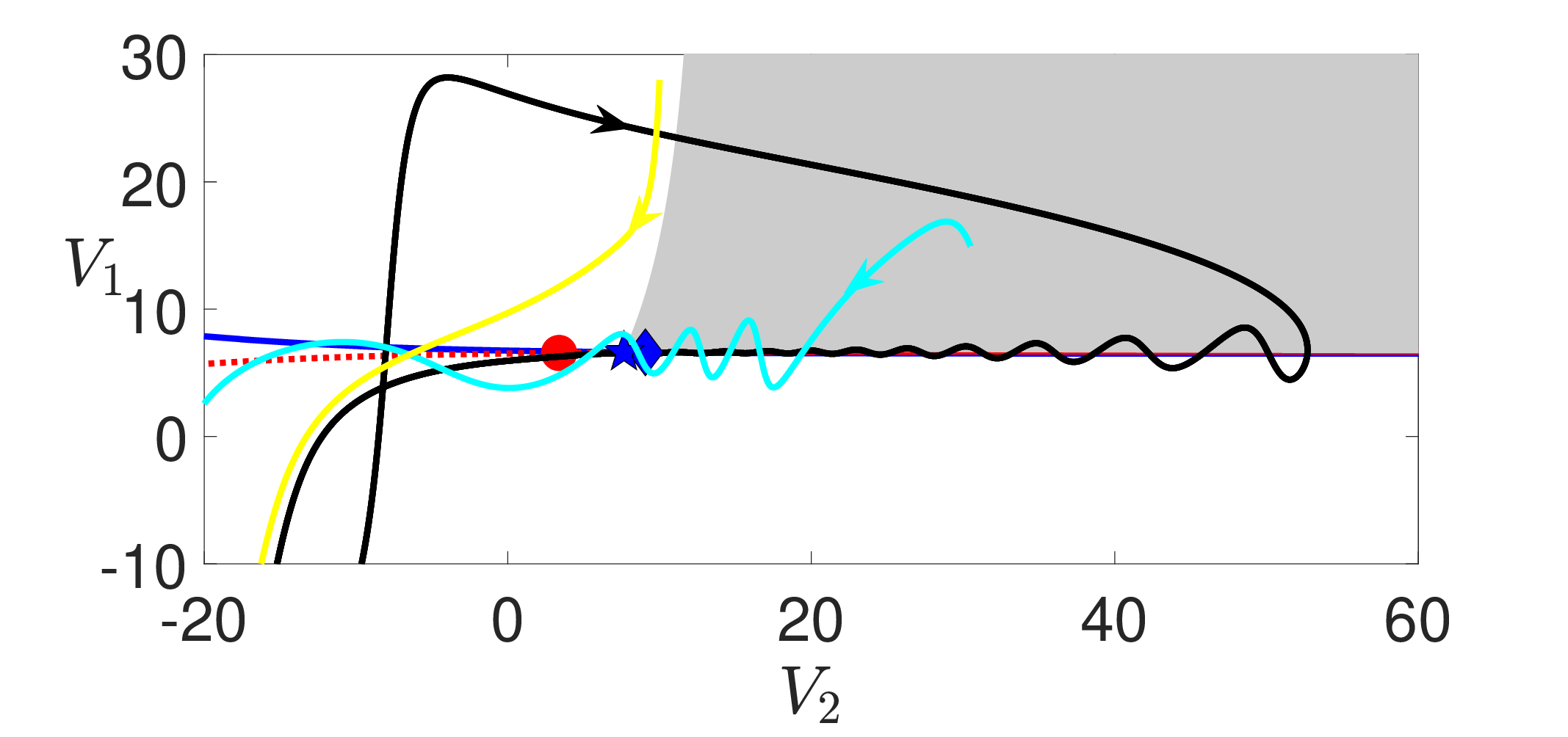}
    \end{center}
\caption{Projections of the solution from Figure \ref{fig:Ms-gsyn4p4} (black) and the singular 3D funnel volume corresponding to the curve of folded nodes onto $(w_2,V_2,V_1)$- space (top) and $(V_2,V_1)$- plane (bottom). The 3D funnel volume is bounded between the singular canard surface (magenta surface) and \RED{the fold surface $\ls$}. The black, cyan and yellow curves denote solutions of \eqref{eq:main} with different initial conditions. Other color coding and symbols have the same meaning as in Figure \ref{fig:proj-v1v2w2-gsyn4p4}.}
\label{fig:funnel-gsyn4p4}
\end{figure}

To confirm that the MMOs exhibit the characteristics of canard dynamics due to the folded node, we plot the funnel volume corresponding to the folded node singularity curve on the upper fold $\ls$ (see Figure \ref{fig:funnel-gsyn4p4}). As discussed in \S\ref{subsec:canard}, the funnel for a folded node point represents a two-dimensional trapping region on $\ms$. The 2D funnels for all folded nodes on the \RED{upper fold} together form a three-dimensional funnel volume. In the $(V_1, V_2, w_2)$ projection (see the \RED{top} panel of Figure \ref{fig:funnel-gsyn4p4}), the funnel volume is bounded by the singular strong canard surface (shown in magenta) and the upper fold surface (in blue). The \RED{bottom} panel shows the projection onto the $(V_1,V_2)$-plane, where the funnel is indicated by the shaded region. Trajectories initiating inside the funnel (e.g., the cyan and black curves) are filtered through the CDH region and there exist SAOs, whereas trajectories starting outside the funnel (e.g., the yellow curve) cross the fold $\ls$ at a regular jump point and there are no SAOs. These observations suggest that the MMOs for $g_{\rm syn}=4.4$ exhibit canard-like features and are organized by the canard mechanism.

Next, we elucidate that the DHB mechanism also contributes to the occurrence of SAOs in the MMO solution.

\subsubsection{Delayed Hopf bifurcation mechanism} 

\begin{figure}[!htp] 
    \begin{center}
        \includegraphics[width=\linewidth]{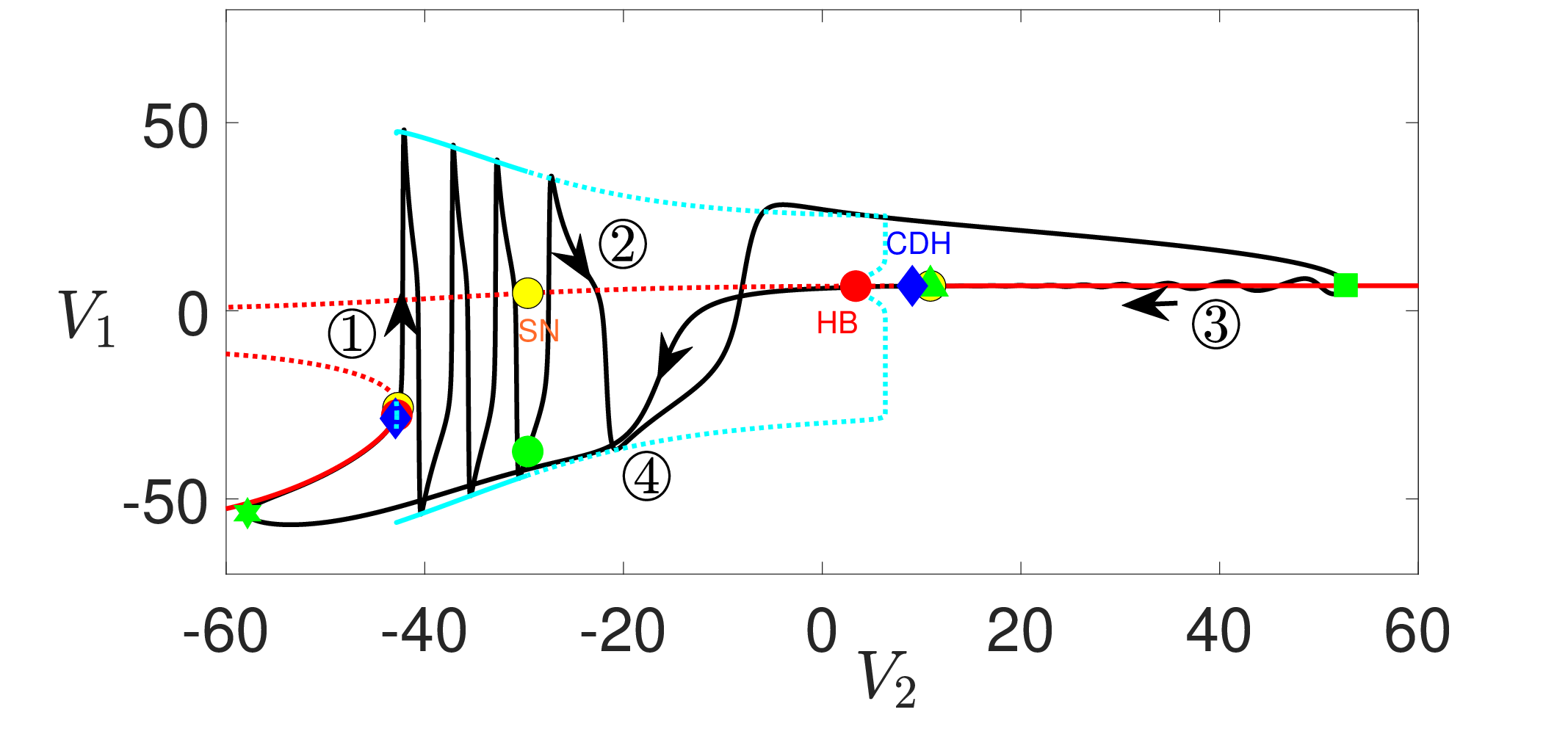}
    \end{center}
\caption{Projection of the solution trajectory and superslow manifold $M_{SS}$ (red curve) from Figure \ref{fig:Ms-gsyn4p4} onto ($V_2,V_1$)-plane. 
The solid (resp. dashed) cyan curves denote stable (resp. unstable) periodic orbit branches. Yellow circles represent saddle-node bifurcations of $\mss$, in which the upper two have the same $V_2$ values as the folds of the $V_2$-nullcline at green circle and triangle. Other colors and symbols are the same as in Figure \ref{fig:Ms-gsyn4p4}.}
\label{fig:Mss-gsyn4p4}
\end{figure}

Figure \ref{fig:Mss-gsyn4p4} shows the projection of the solution trajectory and the bifurcation diagram of the slow layer problem \eqref{eq:slowlayer} onto the $(V_2,V_1)$-plane.  Starting at the green square, the trajectory exhibits SAOs as it follows the upper branch of $\mss$ towards the left. As the trajectory passes through the attracting region of $\mss$, the oscillation amplitude decays in a typical DHB fashion. Moreover, the orbit experiences a delay along the repelling $\mss$ for an amount of time as it passes through the HB. It is worth noting that after the trajectory enters the unstable part of $\mss$, there are no symmetric oscillations with respect to the DHB, i.e., there are no SAOs with increasing amplitudes. This is due to the fact that $V_2$ switches from a superslow to a slow timescale at the green triangle, as discussed earlier when $g_{\rm syn}=4.3$. 

Thus, for $\delta = \mathcal{O}(\epsilon)$, the MMOs for $g_{\rm syn}=4.4$ exhibit characteristics of both the canard and DHB mechanisms. To further confirm this, we performed two perturbations on the system. Firstly, we increased $\varepsilon$ by raising the value of $C_1$ to make the canard mechanism less relevant. This is because increasing $C_1$ slows down the evolution speed of $V_1$, which in turn drives the three-timescale system \eqref{eq:main} closer to (3S, 1SS) splitting.
We observed that MMOs persisted for $C_1$ as large as 80 ($\varepsilon=O(1)$), at which the folded singularities were no longer relevant (see Figure \ref{fig:4p4_increasing_c1}). Secondly, we increased $\delta$ by raising $\phi_2$ to drive the system closer to (1F, 3S) splitting, and observed that SAOs persisted for $\phi_2$ as large as 0.01 ($\delta=O(1)$), at which the DHB mechanism was no longer relevant. 
The persistence of SAOs even when one of the two mechanisms vanishes further highlights the coexistence and interplay of the canard and DHB mechanisms for supporting SAOs. 

\subsection{Effects of varying $\epsilon$ and $\delta$ on MMOs}\label{sec:effect-timescale-4p4}

Unlike the sensitivity of MMOs at $g_{\rm syn}=4.3$ to timescale variations, the interaction of canard and DHB mechanisms due to the existence of the upper CDH when $g_{\rm syn}=4.4$ makes MMOs much more robust, as discussed above and illustrated in Figure \ref{fig:c1-phi2}B. 
Specifically, MMOs persist over biologically relevant ranges of $C_1\in (0.1, 80)$ and $\phi_2\in(1e-4, 0.01)$ (also see Figures \ref{fig:4p4_increasing_c1} and \ref{fig:4p4_increasing_phi2}). 

The robustness of MMOs or SAOs to decreasing $\varepsilon$ or $\delta$ is expected, as it moves either the DHB point or the folded saddle-node singularity closer to the CDH, causing them to move into the midst of the small oscillations (see Figure \ref{fig:4p4_increasing_c1}A and Figure \ref{fig:4p4_increasing_phi2}A). In this subsection, we only explain the effects of increasing $\varepsilon$ or $\delta$ on the features of small oscillations \RED{within} MMOs, as decreasing them yields analogous effects but reversed. \RED{We find that MMOs with $\delta=\mathcal{O}(\varepsilon)$ exhibit both canard- and DHB-like features, whereas by tuning $\delta\geq \mathcal{O}(\sqrt{\varepsilon})$, DHB-like features diminish and the canard mechanism dominates.
}

Below we summarize the effects of increasing the singular perturbations:
\begin{itemize}
\item For fixed $\phi_2=0.001$, increasing $\varepsilon$ (i.e., increasing $C_1$) enhances the DHB-like features of the SAOs (see Figure \ref{fig:4p4_increasing_c1}). 
\item For fixed $C_1=8$, increasing $\delta$ (i.e., increasing $\phi_2$) makes the canard mechanism dominant and causes DHB-like features to gradually vanish (see Figure \ref{fig:4p4_increasing_phi2}). 
\end{itemize}



\subsubsection{Increasing $C_1$ makes DHB dominate}

\begin{figure*}[!htp]
    \begin{center}
        \begin{tabular}
        {@{}p{0.48\linewidth}@{\quad}p{0.48\linewidth}@{}}{\subfigimg[width=\linewidth]{\bf{\small{(A1)}}}{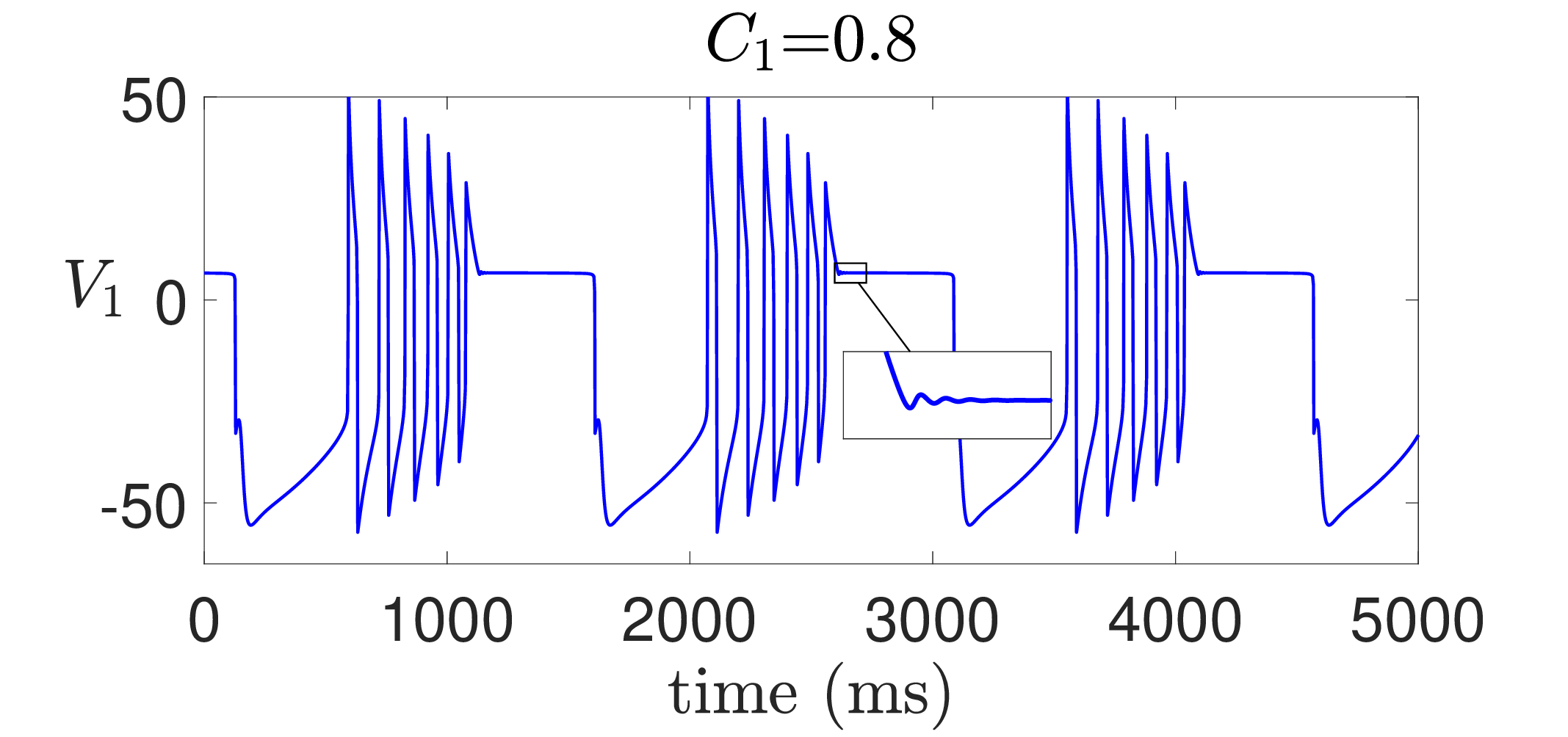}}&{\subfigimg[width=\linewidth]{\bf{\small{(A2)}}}{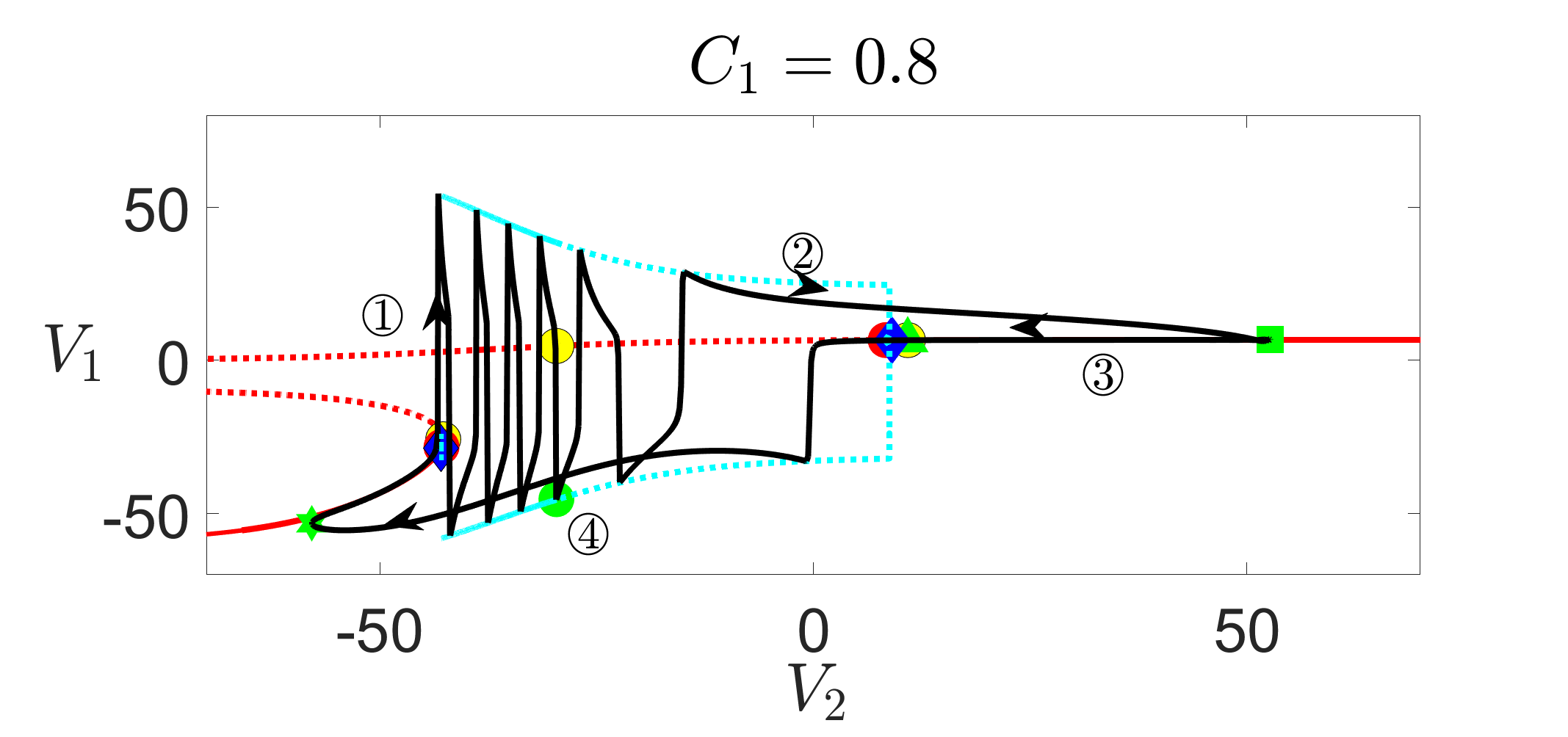}}\\
            \subfigimg[width=\linewidth]{\bf{\small{(B1)}}}{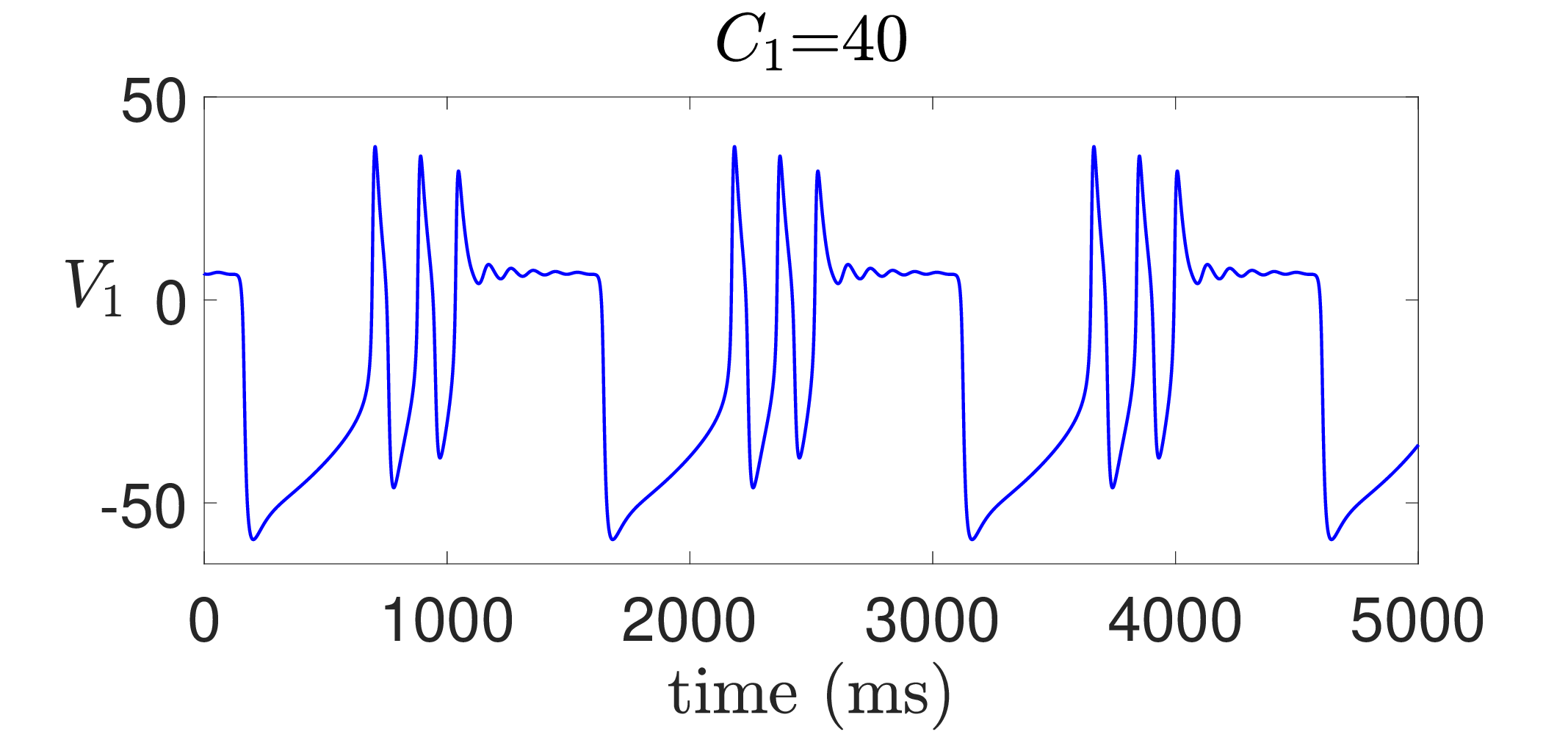}&
            \subfigimg[width=\linewidth]{\bf{\small{(B2)}}}{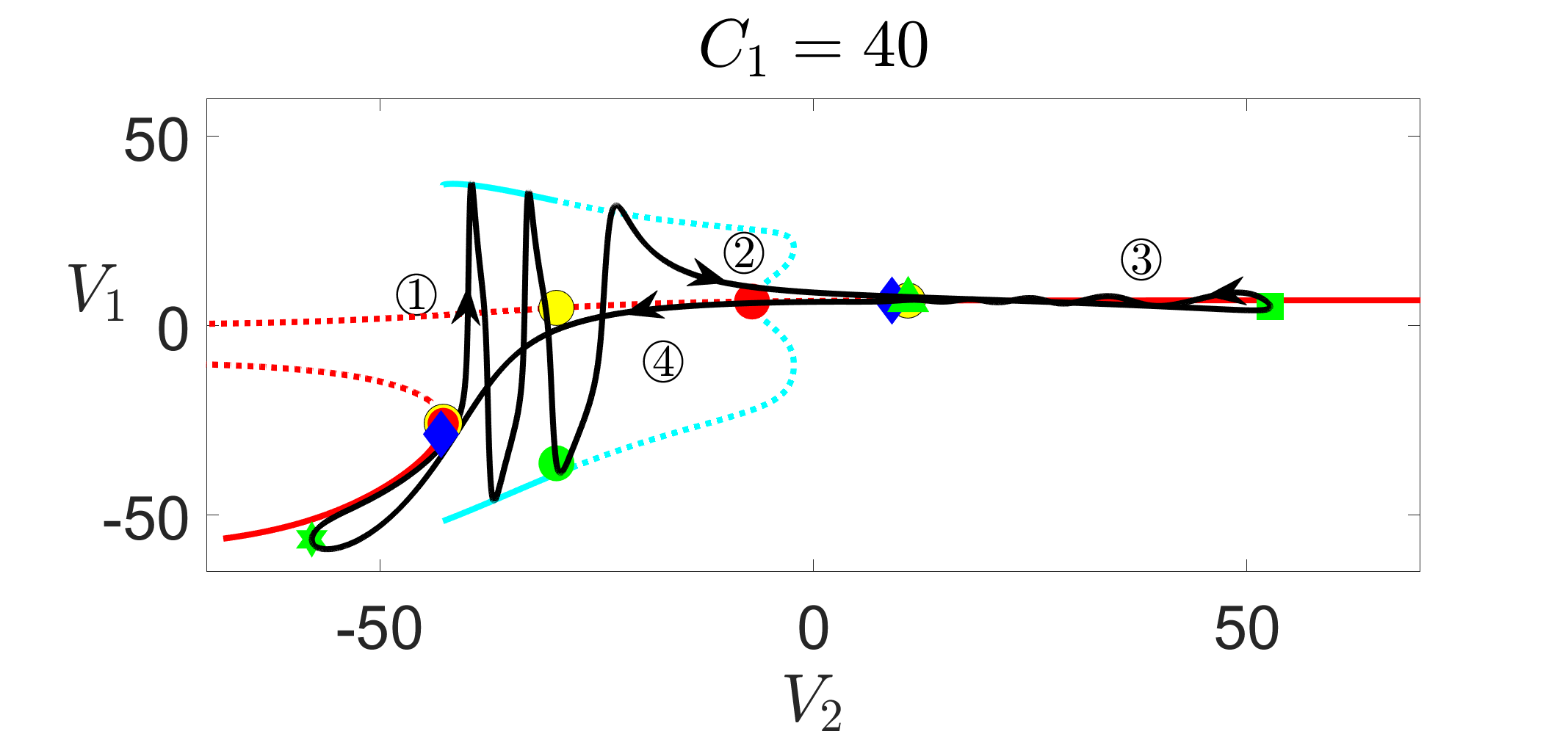}\\
            \subfigimg[width=\linewidth]{\bf{\small{(C1)}}}{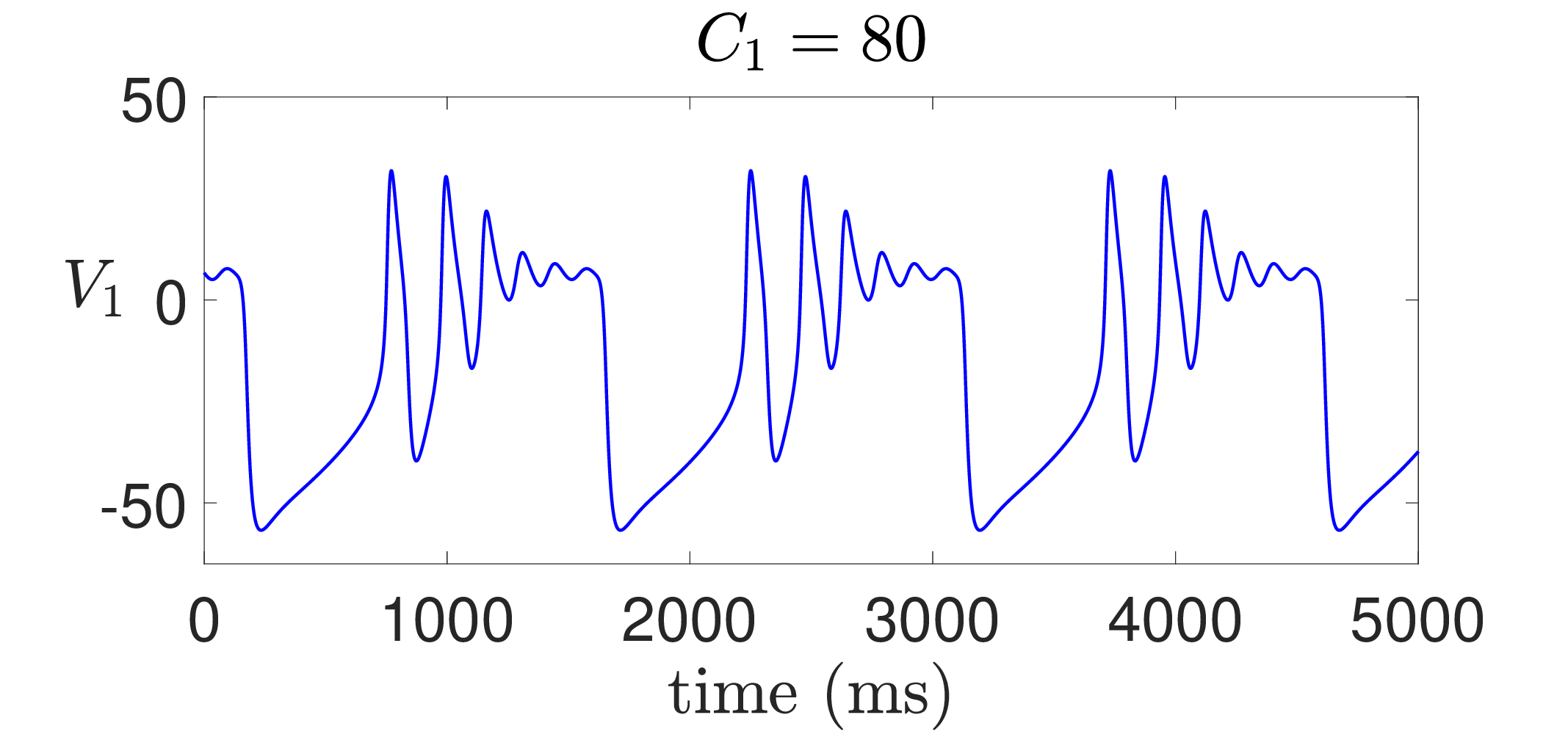}&
            \subfigimg[width=\linewidth]{\bf{\small{(C2)}}}{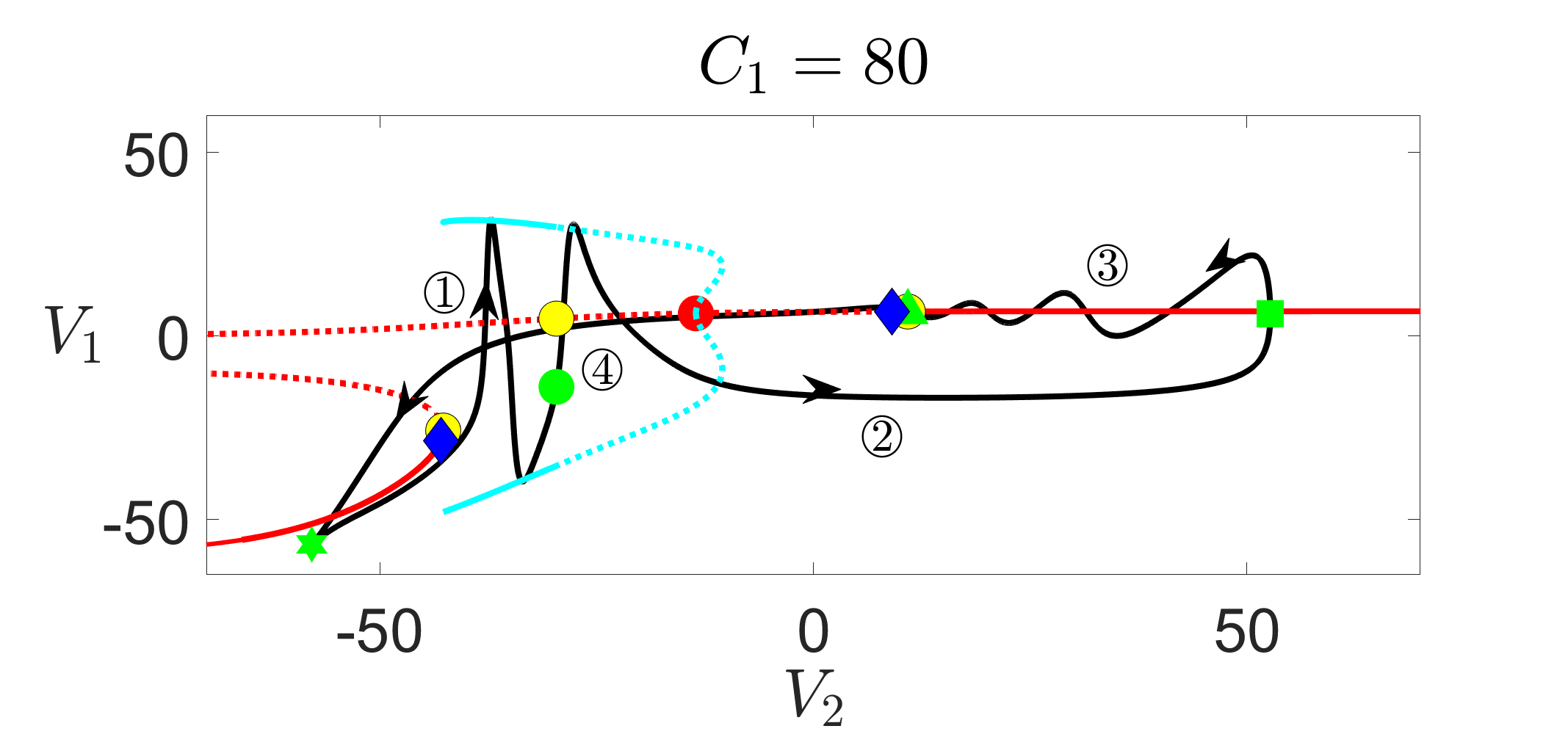}
        \end{tabular}
    \end{center}
\caption{Time traces of the solutions of \eqref{eq:main} (left panels) and their projections onto the corresponding bifurcation diagrams in $(V_2,V_1)$-space (right panels) for $g_{\rm syn}=4.4$ and (A) $C_1=0.8$ $(\epsilon=0.01)$, (B) $C_1=40$ $(\epsilon=0.5)$, and (C) $C_1=80$ ($\epsilon=1$). Increasing $C_1$ moves the upper DHB (red circle) to the lower $V_2$ and manifests DHB-like features. Color coding and symbols of the bifurcation diagrams on the right panel are the same as in Figure \ref{fig:Ms-gsyn4p4}. }
\label{fig:4p4_increasing_c1} 
\end{figure*}

Increasing $C_1$ drives the three-timescale system \eqref{eq:main} closer to (3S, 1SS) splitting. As a result, the critical manifold $M_s$ and folded singularities including the CDH point become less meaningful and eventually irrelevant for sufficiently large $C_1$ (e.g., $C_1=80$). Despite this, MMOs persist due to the existence of the DHB mechanism. Moreover, the DHB mechanism becomes more dominant in controlling the features of the small oscillations as $C_1$ increases, while the influence of the canard mechanism becomes less significant.

Figure \ref{fig:4p4_increasing_c1} shows the effect of increasing $C_1$ on voltage traces of the full system and the bifurcation diagrams of the fast subsystem \eqref{eq:slowlayer}. As $C_1$ increases from panel (A2) to panel (C2), the upper DHB (red circle) moves away from the upper CDH (blue diamond) to smaller $V_2$ values, whereas the CDH points and slow/superslow timescale transitions denoted by yellow and green symbols all remain unaffected by $C_1$. As a result, the trajectory with larger $C_1$ begins small oscillations at a larger distance from the Hopf bifurcation point, similar to what we observed in the case of $\gsyn=4.3$. Furthermore, we have noticed that trajectories for larger $C_1$ exhibit more pronounced DHB-like characteristics, including SAOs with decreasing amplitudes and a more extended travel distance along the unstable branch of $\mss$. As discussed before, due to a switch of the $V_2$ timescale at the green triangle, there are no oscillations with growing amplitudes as one would expect in a typical DHB fashion. 

\subsubsection{Increasing $\phi_2$ makes the canard mechanism dominate}

\begin{figure*}[!t]
    \begin{center}
        \begin{tabular}{@{}p{0.48\linewidth}@{\quad}p{0.48\linewidth}@{}}
            \subfigimg[width=\linewidth]{\bf{\small{(A1)}}}{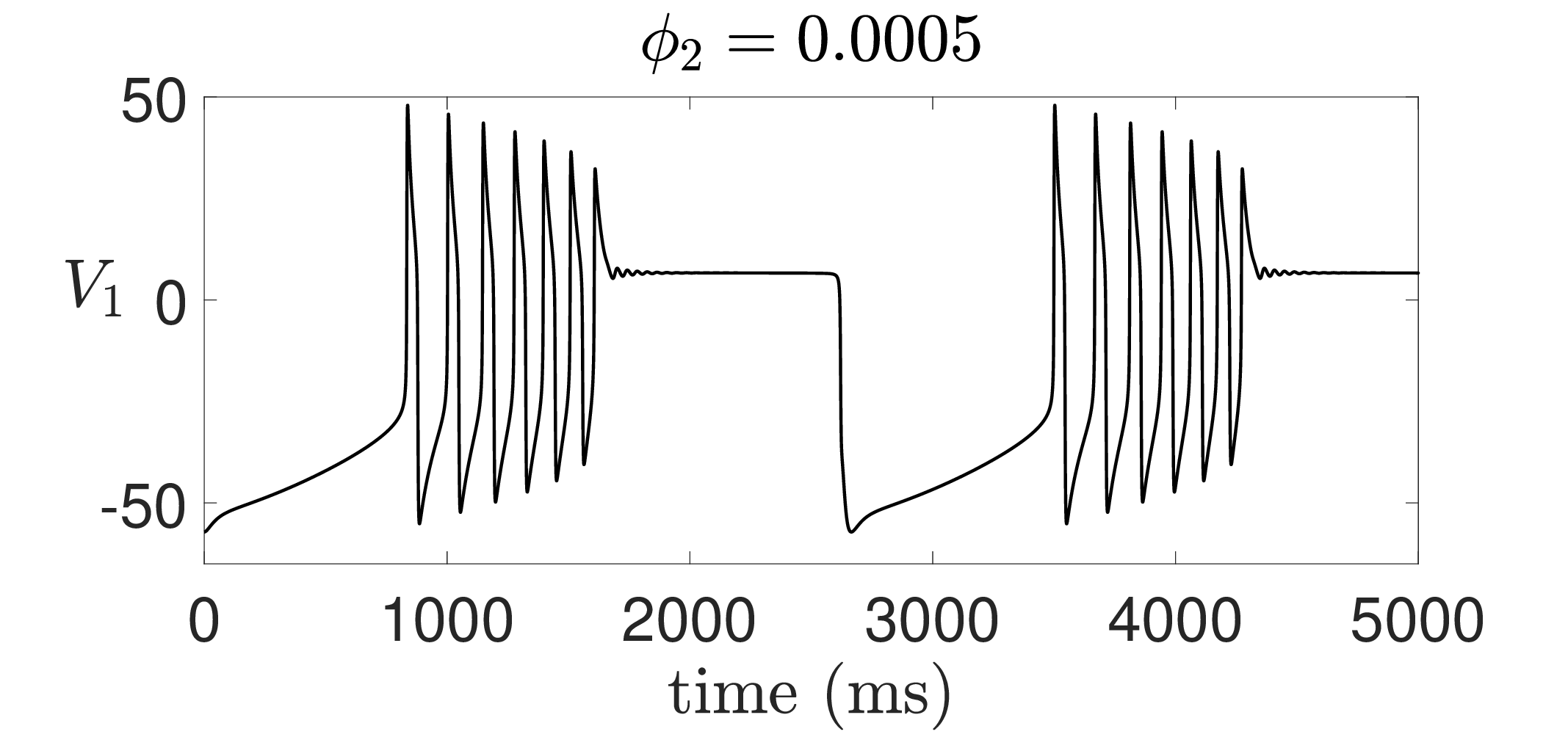}&
            \subfigimg[width=\linewidth]{\bf{\small{(A2)}}}{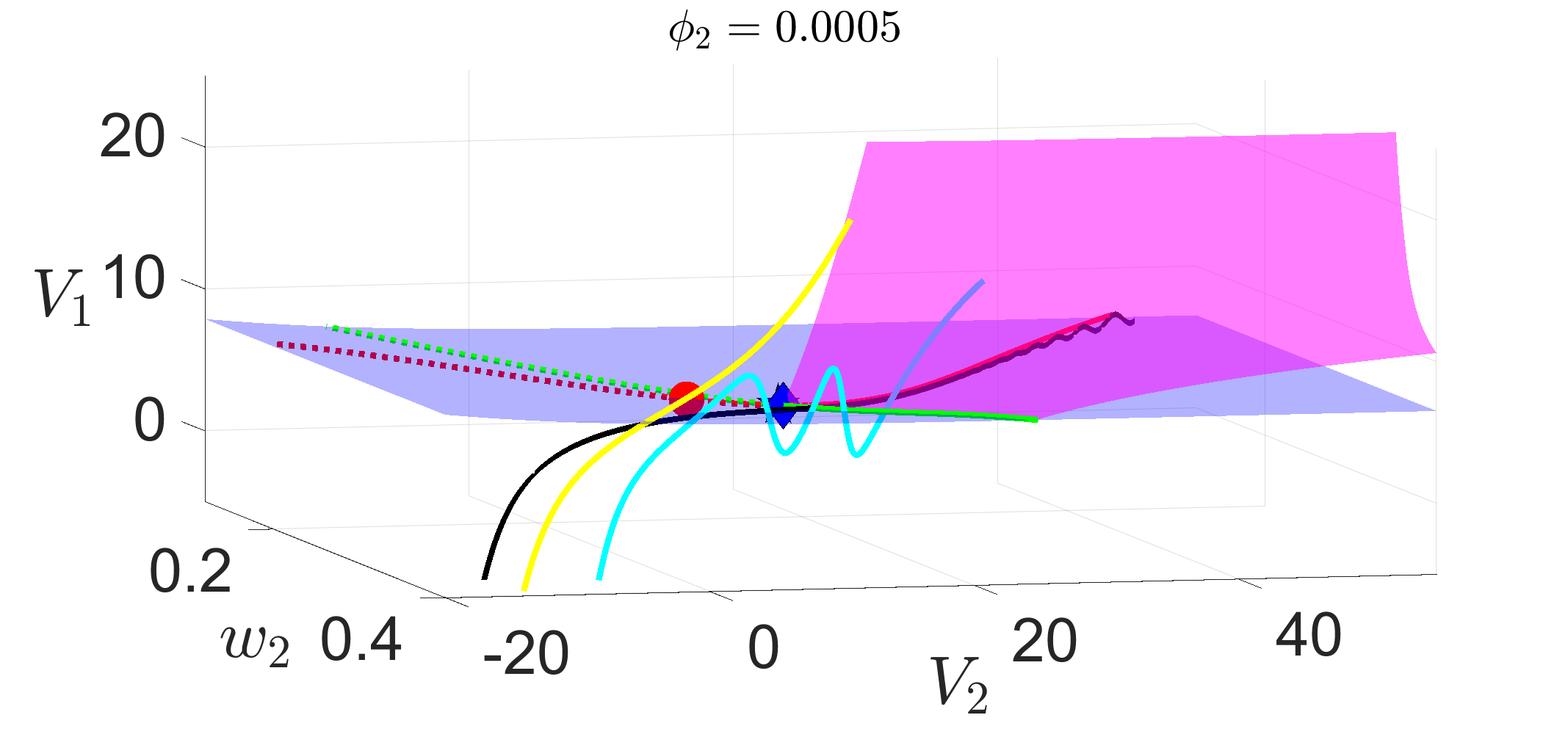}\\
            \subfigimg[width=\linewidth]{\bf{\small{(B1)}}}{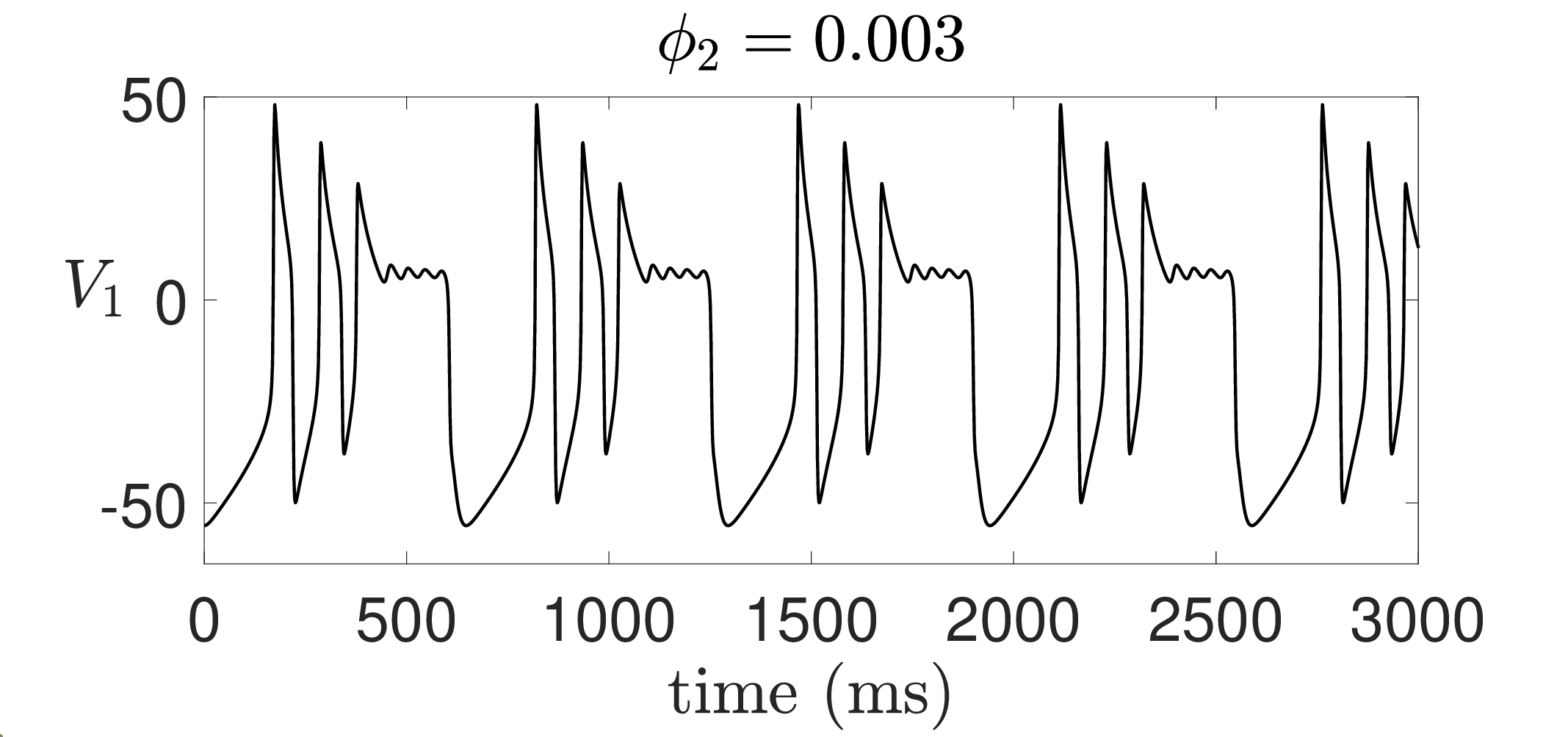}&
            \subfigimg[width=\linewidth]{\bf{\small{(B2)}}}{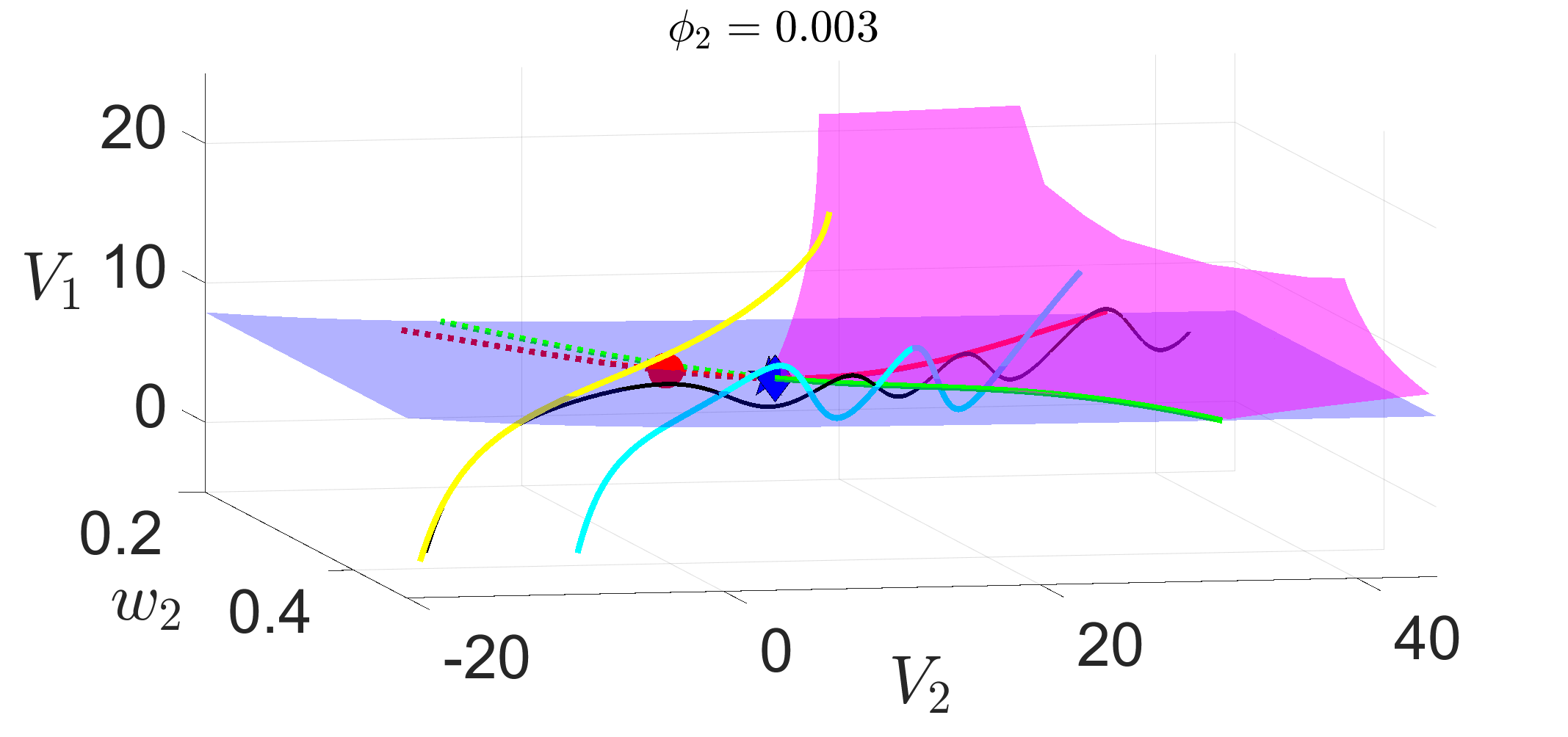}\\
            \subfigimg[width=\linewidth]{\bf{\small{(C1)}}}{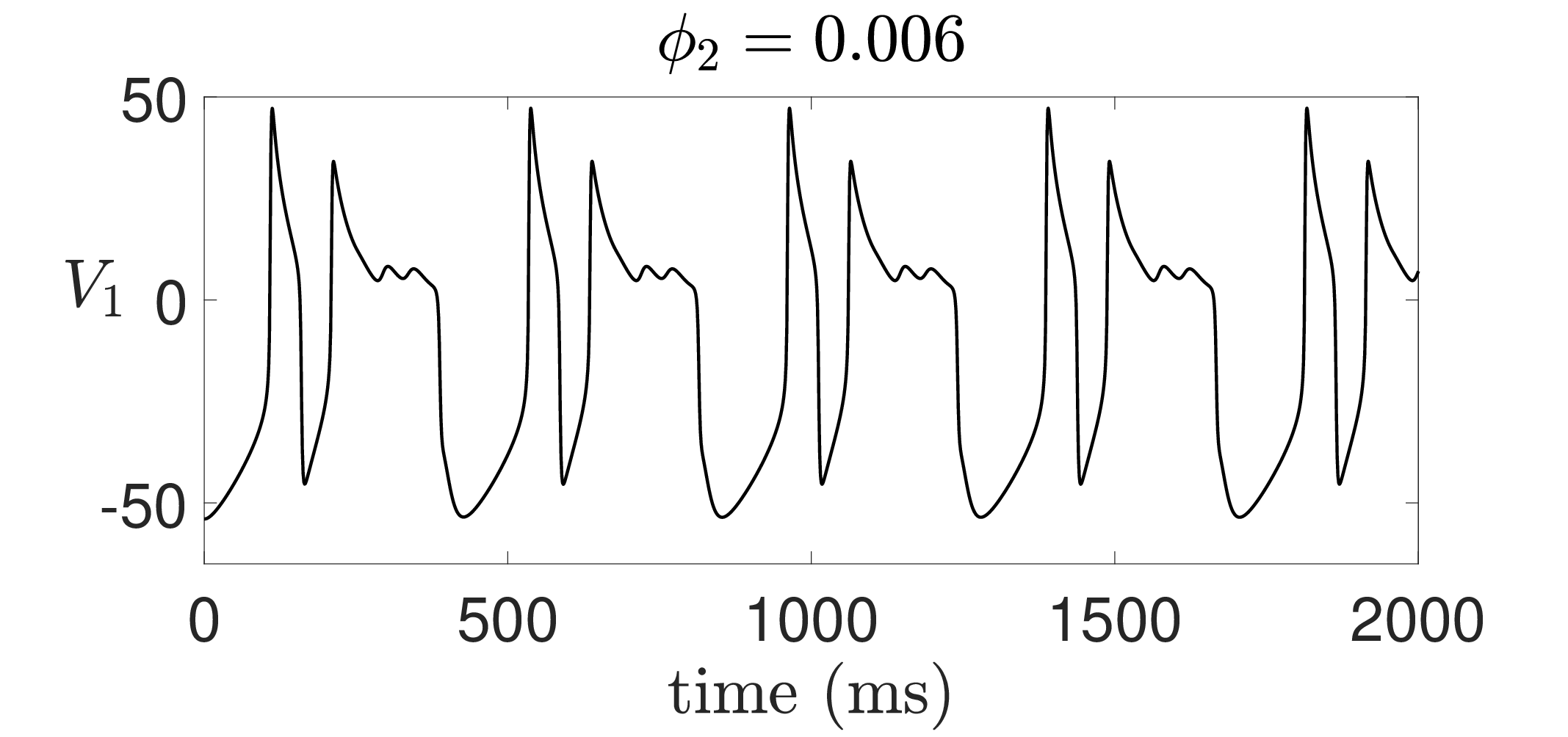}&
            \subfigimg[width=\linewidth]{\bf{\small{(C2)}}}{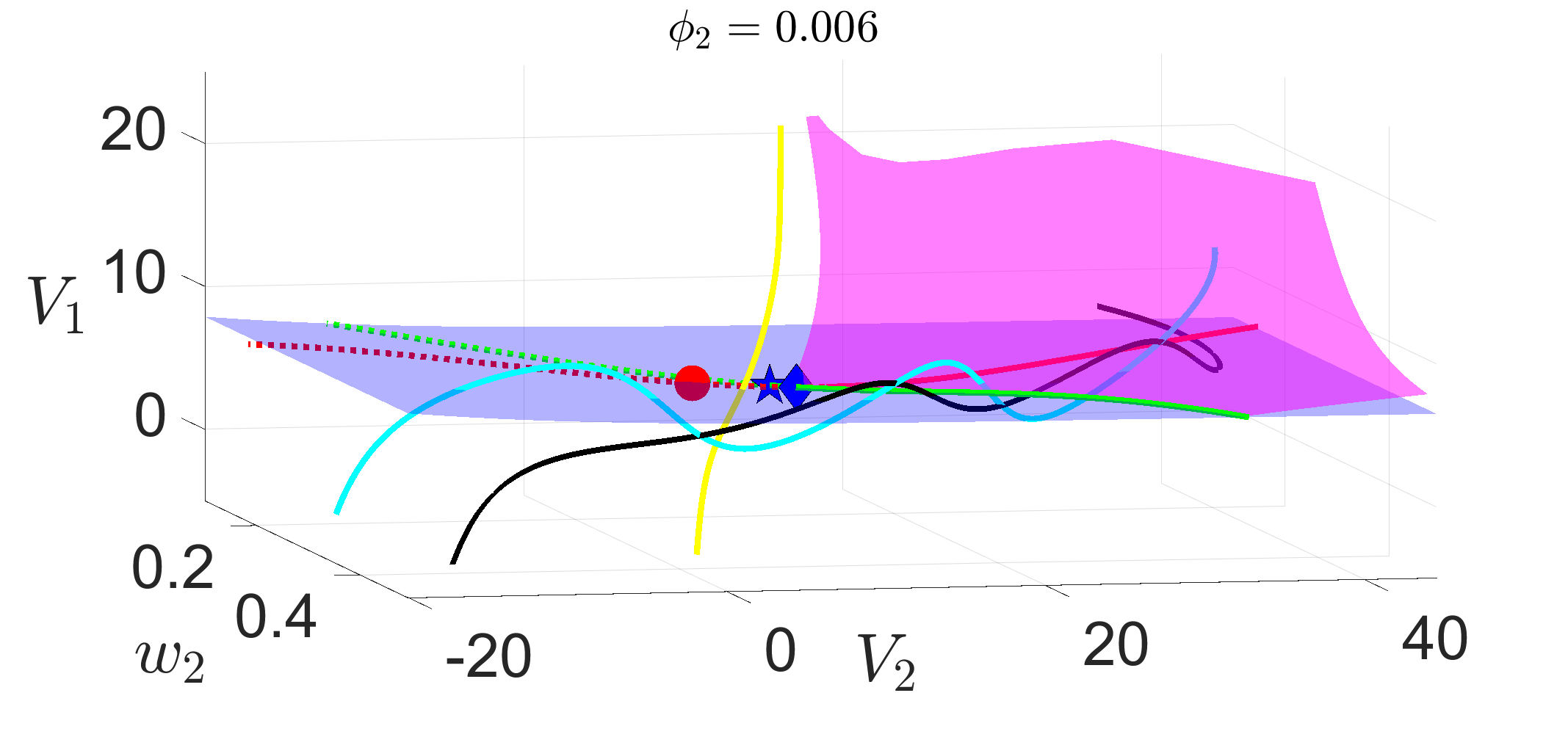}\\
            \subfigimg[width=\linewidth]{\bf{\small{(D1)}}}{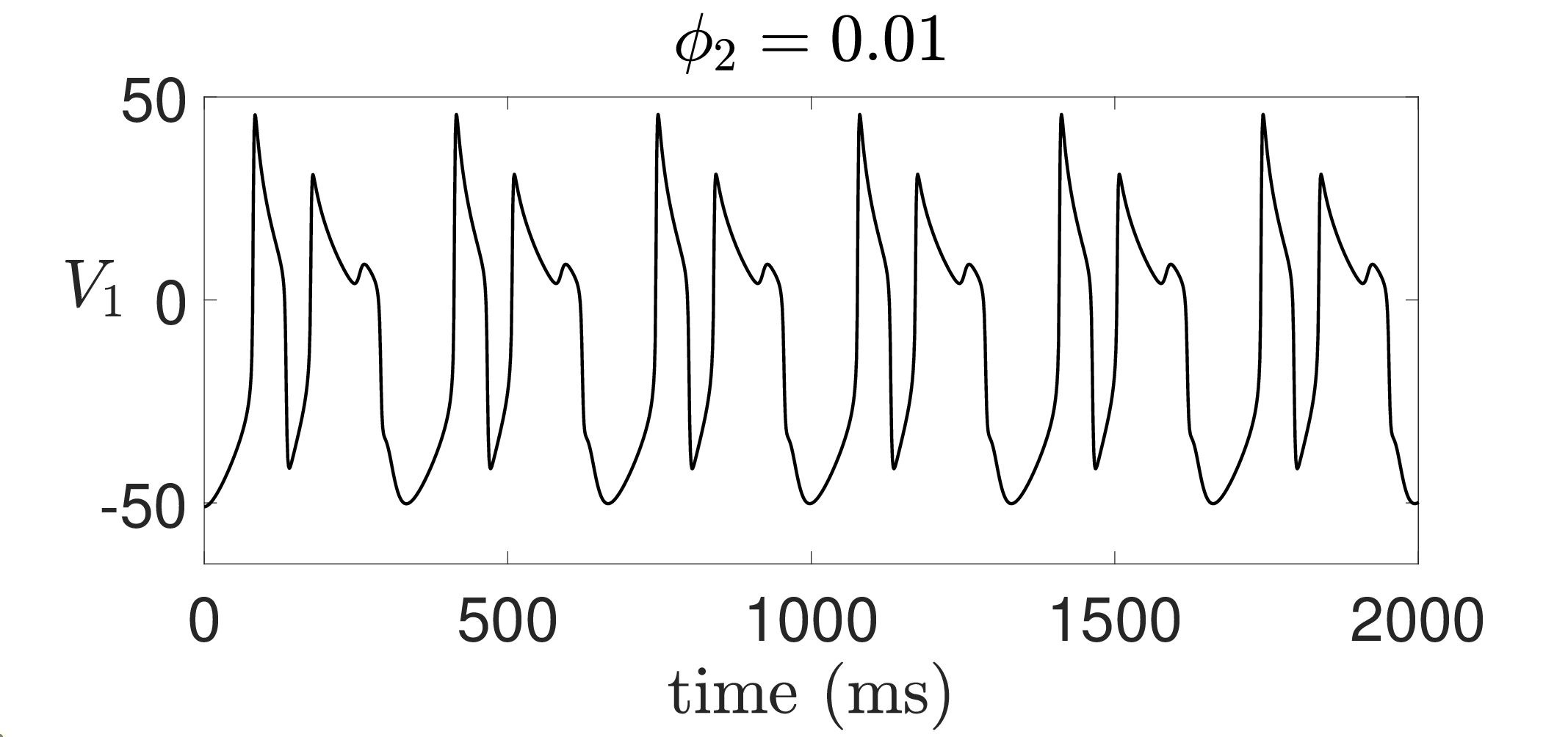}&
            \subfigimg[width=\linewidth]{\bf{\small{(D2)}}}{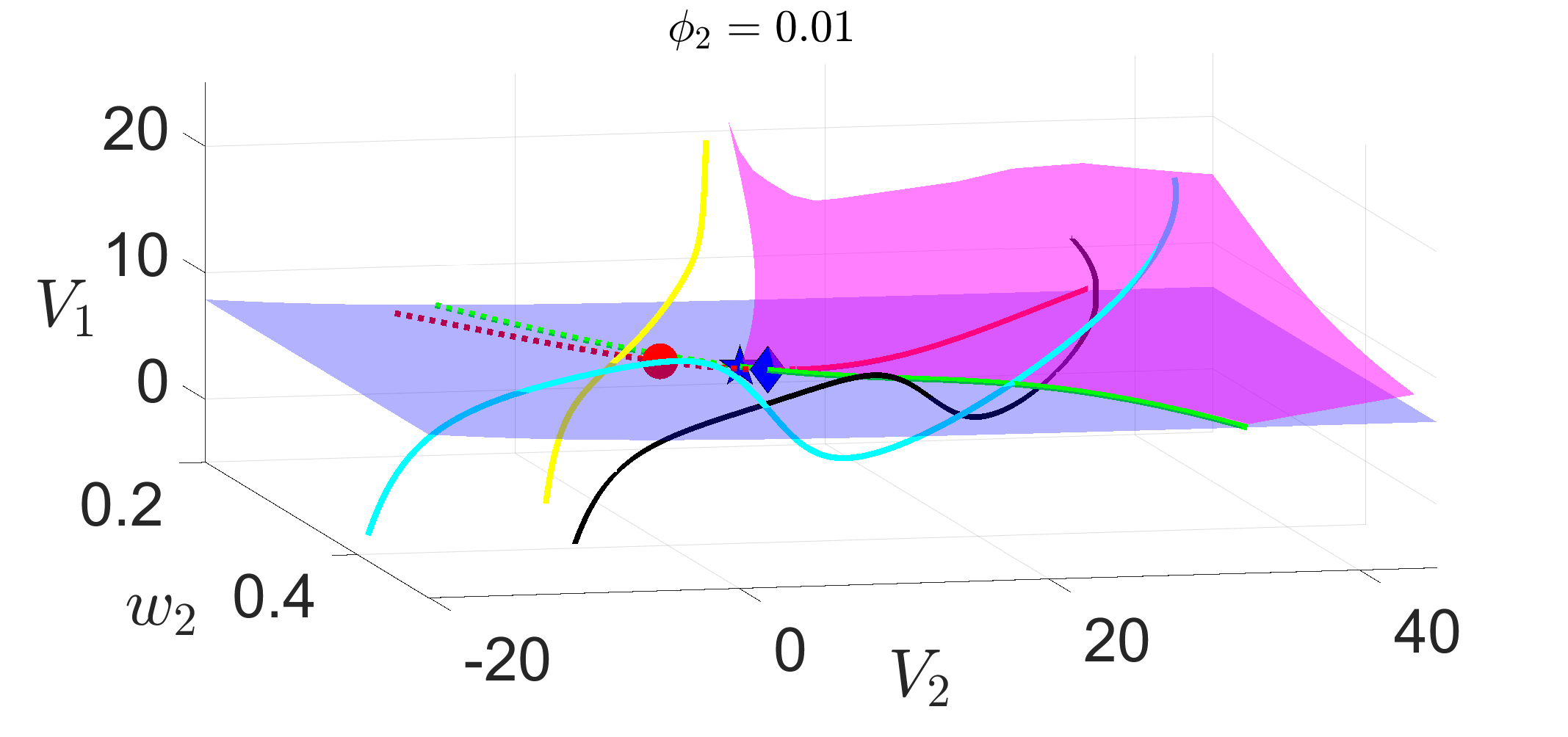}
        \end{tabular}
    \end{center}
\caption{Time traces of the solutions of \eqref{eq:main} (left panels) and their projections (black curves) onto $(w_2,V_2,V_1)$-space (right panels) for $g_{\rm syn}=4.4$ and (A) $\phi_2=0.0005$ (i.e., $\delta=0.0265 = \frac{1}{3}\epsilon$), (B) $\phi_2=0.003$ (i.e., $\delta=0.159 = \frac{1}{2} \sqrt\epsilon$), (C) $\phi_2=0.006$ (i.e., $\delta=0.318 = \sqrt\epsilon$) and (D) $\phi_2=0.01$ (i.e., $\delta = 0.56 = \epsilon^{\frac{1}{4}}$). Increasing $\phi_2$ preserves MMOs with more canard-like feature. Codings and symbols on the right panels are the same as in Figure \ref{fig:funnel-gsyn4p4}.}
\label{fig:4p4_increasing_phi2}
\end{figure*}

Increasing $\phi_2$ speeds up the superslow variable $w_2$ and hence  drives the three-timescale system \eqref{eq:main} closer to (1F, 3S) splitting. As a result, the superslow manifold $\mss$ and the DHB points become less relevant and eventually no longer meaningful for sufficiently large $\phi_2$ (e.g., $\phi_2=0.01$). Nonetheless, MMOs continue to persist due to the existence of the canard mechanism. Figure \ref{fig:4p4_increasing_phi2} reveals the presence of canard-like features for $\phi_2$ values across a wide range (from $0.0005$ to $0.01$). That is, trajectories within the funnel volume (cyan and black curves) exhibit SAOs near the CDH,  whereas those outside the funnel (yellow curves) display no oscillations. 

Moreover, as $\phi_2$ increases, small oscillations tend to pull away from $\mss$ and lose their DHB-like features (i.e., oscillations with decaying amplitude and a delay after passing the HB). Specifically, MMOs with $\delta=O(\varepsilon)$ (e.g., Figures~\ref{fig:proj-v1v2w2-gsyn4p4} and \ref{fig:4p4_increasing_phi2}A) exhibit both canard and DHB characteristics. When $\delta\approx \frac{1}{2}\sqrt{\varepsilon}$ (Figure \ref{fig:4p4_increasing_phi2}B), there are still some DHB-like features. Further increasing $\delta$ to  $\delta=0.006\approx \sqrt{\varepsilon}$ or $\delta=0.01\approx \varepsilon^{\frac{1}{4}}$ (Figure \ref{fig:4p4_increasing_phi2}C and D), the trajectories no longer closely follow $\mss$ and the amplitudes of small oscillations become almost constant, which reflects the absence of DHB-like features.

In summary, increasing $\phi_2$ makes the canard mechanism dominate and the SAOs exhibit fewer DHB-like features. Conversely, decreasing $\phi_2$ brings the solution and $\mss$ closer together and amplifies the DHB characteristics of the sustained MMOs (see Figure \ref{fig:4p4_increasing_phi2}). 

\section{Discussion}\label{sec:discussion}

Mixed mode oscillations (MMOs) are commonly exhibited in dynamical systems that involve multiple timescales. These complex oscillatory dynamics have been observed in various areas of applications and are well studied in two-timescale settings \citep{Brons2006,SW2001,Wechselberger2005,Baer1989, Neishtadt1987, Neishtadt1988, RW2008,Hayes2016,Desroches2012}. In contrast, progress on MMOs in three-timescale problems has been made only in the recent past (see e.g., \citep{Jalics2010,Vo2013,Letson2017,Krupa2008b, Krupa2012, Kak2022, Kak2023a, Kak2023b}). In this work, we have contributed to the investigation of MMOs in three-timescale settings by considering a four-dimensional model system of coupled Morris-Lecar neurons that exhibit three distinct timescales. We have investigated two types of MMO solutions obtained with different synaptic strengths ($\gsyn=4.3\,\rm mS/cm^2$ and $\gsyn=4.4\,\rm mS/cm^2$). Applying geometric singular perturbation theory and bifurcation analysis \citep{Fenichel1979,Nan2015}, we have revealed that the two MMOs exhibit different mechanisms, leading to remarkably different sensitivities to variations in timescales (see Figure \ref{fig:c1-phi2}). In particular, the MMO solution at $\gsyn=4.3$ only involves one mechanism (the delayed Hopf (DHB)) and appears to be vulnerable to certain timescale perturbations, whereas for $\gsyn=4.4$, two separate mechanisms (canard and DHB) coexist and interact to produce more robust MMOs.  

The existence of three distinct timescales leads to two important subjects: the critical (or slow) manifold $\ms$ and the superslow manifold $\mss$. The point where a fold of the critical manifold $\ls$ intersects the superslow manifold $\mss$ is referred to as the canard-delayed-Hopf (CDH) singularity, which naturally arises in problems that involve three different timescales. Ref. \onlinecite{Letson2017} considered a common scenario in three-timescale systems where the CDH singularity exists and proved the existence of canard solutions near the CDH singularity for sufficiently small $\varepsilon$ and $\delta$. Moreover, small-amplitude oscillations (SAOs) constantly occur in the vicinity of a CDH singularity \citep{Vo2013,Letson2017}. 

In this work, we have reported several key findings that have not been previously observed in three-timescale systems. First, although we have identified the same type of CDH singularity as documented in \citep{Letson2017} with its center manifold transverse to $\ls$ at the CDH point (see our proof in Appendix \ref{app:CDH-vc}), SAOs in our system are not guaranteed to occur in the neighborhood of a CDH singularity. Specifically, we observed that CDH singularities on the lower $\ls$ did not support SAOs in either the case of $\gsyn=4.3$ or $\gsyn=4.4$. We have explained in subsection \ref{sec:noSTOs-CDH} why neither of the canard and DHB mechanisms near the lower CDH gives rise to SAOs. Our analysis suggests that the absence of SAOs near the lower CDH might be due to the proximity of this CDH to the actual fold of the superslow manifold $\mss$. Further analytical work is still required to confirm this observation and should be considered for future work. 
Secondly, we have explored the conditions underlying the robust occurrence of MMOs in a three-timescale setting. Our analysis has revealed that the existence of CDH singularities critically determines whether or not the two MMO mechanisms (canard and DHB) can coexist and interact, which in turn greatly impacts the robustness of MMOs. In particular, we have found that MMOs occurring near a CDH singularity are much more robust than MMOs when CDH singularities are absent. 

In addition to uncovering the relationship between CDH singularities and robustness of MMOs, we have also provided a detailed investigation on how the features and mechanisms of MMOs without or with CDH singularities vary with respect to timescale variations. Table \ref{tab:vary-eps-delta} outlines a summary of the different mechanisms and robustness properties of MMOs as we vary timescales. When $\gsyn=4.3$, where no CDH was found near the small oscillations, the only mechanism for the MMOs is the DHB mechanism as we justified in Section \ref{4p3}. 
In this case, speeding up the fast variable $V$ via decreasing $\varepsilon$ led to a total of three transitions between MMOs and non-MMO solutions due to its effect on the upper DHB point (Figure \ref{fig:two-par-4p3}A). Initially, MMOs disappeared as the decrease of $\varepsilon$ moved the DHB closer to the green square where the SAO phase began, resulting in insufficient time for generating small oscillations (see Figure \ref{fig:de-c1}A). However, as the further reduction of $\varepsilon$ resulted in the cross of the DHB with the green square or a complete vanish of the DHB, we observed a recovery of MMOs originating from the saddle foci equilibria along the superslow manifold $\mss$ (see Figures \ref{fig:de-c1}B,C and \ref{fig:4p3-saddlefocus}). Eventually, SAOs disappeared entirely when $V$ became so rapid that it failed to remain in proximity of $\mss$ to generate small oscillations. 

As one would expect, MMOs when $\gsyn=4.3$ is also sensitive to increasing $\delta$, which speeds up the superslow variable and thus makes the DHB mechanism less relevant. Interestingly, however, increasing $\delta$ does not just simply eliminate MMOs. Instead, it led to a total of five transitions between MMOs and non-MMO states (see Figures~\ref{fig:c1-phi2} and \ref{fig:gsyn_4p3-phi2}B). Our analysis suggests that these complex transitions occur due to a spike-adding like mechanism. When no big spike is lost with the increase of $\delta$, a transition from MMOs to non-MMOs will take place. However, if an entire big spike is lost during this process, changes to the SAOs will be reversed and MMOs will recover again. Ultimately, MMOs will be completely lost as $\delta$ is increased to a point where the DHB mechanism is no longer relevant.  

On the other hand, MMOs at $\gsyn=4.3$ show strong robustness to increasing $\varepsilon$ or decreasing $\delta$. This is not surprising as both of these changes manifest the DHB mechanism by moving the three timescale problem closer to (3S, 1SS) separation. As a result, we observed more DHB characteristic in the MMOs as demonstrated in the case of increasing $\varepsilon$ (see Figures \ref{fig:gsyn_4p3_varyC1}A and \ref{fig:gsyn4p3-bif-largerC1}). Nonetheless, decreasing $\delta$ led to an interesting non-monotonic effect on SAOs, where the amplitude and the number of SAOs exhibit an alternation between increase and decrease as $\delta$ is reduced (see Figure \ref{fig:gsyn_4p3-phi2}A). Our analysis showed that the mechanism underlying such non-monotonic effects on SAOs over the decrease of $\delta$ is similar to the mechanism that drives multiple MMOs/non-MMOs transitions as $\delta$ is increased. 

Unlike $\gsyn=4.3$, an upper CDH occurs near the SAOs at $\gsyn=4.4$. In this case, we showed that this CDH allowed for the coexistence and interaction of canard and DHB mechanisms, resulting in MMOs with strong robustness against timescale variations. Since both MMO solutions for $\gsyn=4.3$ and $4.4$ exhibit DHB mechanisms, they show similar responses and robustness to increasing $\varepsilon$ and decreasing $\delta$, both of which lead to more DHB characteristic in the MMOs. In contrast to $\gsyn=4.3$, MMOs at $\gsyn=4.4$ are also robust against changes that speed up fast or superslow variables. This is due to the existence of the canard mechanism in addition to the DHB. Instead of eliminating the upper DHB in the case of $\gsyn=4.3$ (Figure \ref{fig:two-par-4p3}A), decreasing $\varepsilon$ at $\gsyn=4.4$ brings the DHB point closer to the CDH (Figure \ref{fig:effects-eps-hb-4p4}, right panel) and, consequently, closer to the midst of small oscillations (Figure \ref{fig:4p4_increasing_c1}). Moreover, this timescale change manifests the canard mechanism by moving the system closer to the $\varepsilon\to 0$ singular limit. Hence, MMOs at $\gsyn=4.4$ persist as $\varepsilon$ decreases  and are organized by both mechanisms. On the other hand, increasing $\delta$ diminishes the relevance of the DHB mechanism and moves the trajectory further away from the superslow manifold, resulting in MMOs with more canard-like features. 

To summarize, the coexistence of canard and DHB mechanisms for MMOs at $\gsyn=4.4$, due to the presence of a nearby upper CDH, leads to significantly enhanced robustness against timescale variations compared with $\gsyn=4.3$, where only one mechanism is present. While we have not examined the case when only a canard mechanism is present\RED{, based} on our analysis, we expect that MMOs with only canard mechanism would show more sensitivities to timescale variations that reduces the relevance of the canard mechanism such as increasing $\varepsilon$ or decreasing $\delta$. Similarly, such MMOs should show stronger robustness to decreasing $\varepsilon$ or increasing $\delta$. It would be of interest to explore such a scenario for future investigation. Furthermore, we did not notice any non-monotonic effects on the features of SAOs as $\phi_2$ is decreased when $\gsyn=4.4$, unlike what we observed in $\gsyn=4.3$. This difference is likely attributed to the presence of an additional canard mechanism which may have hindered the occurrence of complex non-monotonic behaviors. A complete analysis of this phenomenon could be investigated in future work. 

\begin{table*}[!t]
\caption{Effects of varying $\epsilon$ and $\delta$ on MMOs.}
\label{tab:vary-eps-delta}
\begin{center}
\begin{tabular}{| m{4cm} | m{4.5cm}| m{5.5cm}|}
\hline
  & $g_{\rm syn}=4.3$  & $g_{\rm syn}=4.4$  \\
 \hline 
 Mechanisms & No upper CDH; Only DHB mechanism & Upper CDH exists; Canard and DHB mechanisms coexist 
\\  \hline
Increasing $C_1$ (or $\epsilon$) & MMOs are preserved with more DHB characteristics. & MMOs are preserved with more DHB characteristics. \\
\hline
Decreasing $C_1$ (or $\epsilon$) & Two MMOs/non-MMOs transitions are observed before a complete loss of MMOs. &  MMOs are preserved and organized by both mechanisms.\\
\hline
Increasing $\phi_2$ (or $\delta$) & Four MMOs/non-MMOs transitions are observed before MMOs are entirely lost. &  MMOs are preserved with more canard-like and less DHB features.\\
\hline
Decreasing $\phi_2$ (or $\delta$) &  MMOs are preserved, but a non-monotonic effect on the small oscillations is observed. &   MMOs are preserved with both DHB and canard characteristics. \\
\hline
\end{tabular} 
\end{center}
\end{table*}

\begin{acknowledgements}
This work was supported in part by NIH/NIDA R01DA057767 to YW, as part of the NSF/NIH/DOE/ANR/BMBF/BSF/NICT/AEI/ISCIII\\Collaborative Research in Computational Neuroscience Program. 

\end{acknowledgements}

\section*{Conflict of Interest}
The authors have no conflicts to disclose.


\section*{Data Availability}
Data sharing is not applicable to this article as no new data were created or analyzed in this study.

\appendix

\section{Dimensional Analysis}\label{ap:nondim}

Here we present dimensional analysis of \eqref{eq:main} to reveal its timescales. To this end, we introduce a dimensionless timescale $t_s=t/Q_t$ with reference timescale $Q_t=\min(\tau_w(V))/\phi_1=18.86\,\rm ms = O(10) \,\rm ms$. This transforms \eqref{eq:main} to the dimensionless system \eqref{eq:slow} given in the introduction, namely,
\begin{widetext}
\begin{equation}\label{eq:slow-appendix}
    \begin{array}{rcl}
        \frac{C_1}{g_{\rm max}Q_t}\frac{dV_1}{dt_s}:=\varepsilon\frac{dV_1}{dt_s} &=& f_1(V_1,\, w_1, \, V_2), \vspace{0.1in}\\
        \frac{dw_1}{dt_s}&=&Q_t\phi_1(w_\infty(V_1)-w_1)/\tau_{w}(V_1)  := g_1(V_1, \, w_1), \vspace{0.1in} \\
        \frac{dV_2}{dt_s}& =&\frac{Q_t g_{\rm max}}{C_2} \bar{f}_2(V_2, \, w_2)  :=  f_2(V_2, \,w_2), \vspace{0.1in}\\
        \frac{dw_2}{dt_s}&=&Q_t\phi_2(w_\infty(V_2)-w_2)/\tau_{w}(V_2) := \delta g_2(V_2, \, w_2),
\end{array}
\end{equation}
where $\varepsilon=\frac{C_1}{g_{\rm max}Q_t}=0.1$,  $\delta=Q_t\phi_2/\mathrm{min}(\tau_w(V_2))=0.053$, $g_{\rm max}=8 \,\rm mS/cm^2$,
\small{\[f_1(V_1,w_1,V_2)=\frac{(I_1-g_{\rm Ca}m_\infty(V_1)(V_1-V_{\rm Ca})-g_{\rm K}w_1(V_1-V_{\rm K})-g_{\rm L}(V_1-V_{\rm L})-g_{\rm syn}S(V_2)(V_1-V_{\rm syn}))}{g_{\rm max}}\]} 
and 
\small{\[\bar{f}_2(V_2,w_2)=\frac{(I_2-g_{\rm Ca}m_\infty(V_2)(V_2-V_{\rm Ca})-g_{\rm K}w_2(V_2-V_{\rm K})-g_{\rm L}(V_2-V_{\rm L})}{g_{\rm max}}\]} 
\end{widetext}

From system \eqref{eq:slow-appendix}, we can see that the voltage $V_1$ evolves on a fast timescale ($\frac{C_1}{g_{\rm max}}=1\,\rm ms$), $(w_1,V_2)$ evolve on a slow timescale ($\frac{\mathrm{min}(\tau_w(V_1))}{\phi_1}\approx \frac{C_2}{g_{\rm max}}=O(10)\,\rm ms$), whereas $w_2$ evolves on a superslow timescale ( $\frac{\mathrm{min}(\tau_w(V_2))}{\phi_2}=O(100)\,\rm ms$).  

\section{Folded saddle-node (FSN) singularities of \eqref{desingu_slowreduced}}\label{app:cond-FSN}

To derive the conditions for FSN singularities, we note that the eigenvalues $\lambda$ of the folded singularities satisfy the following algebraic equation
\begin{equation} \label{equil-folded}
    \lambda^3 - {\rm{tr}}(J_D)\lambda^2+ \frac{1}{2} \left( ({\rm{tr}}(J_D))^2-{\rm{tr}}(J_D^2) \right)\lambda -\det(J_D)=0,
\end{equation}
where ${\rm{tr}}(J_D)$ and $\det(J_D)$ denote the trace and determinant of the Jacobian matrix $J_D$ of the desingularized system \eqref{desingu_slowreduced}, respectively.
As discussed before, $\det(J_D)\equiv 0$ along the folded singularity curve \RED{$\mathcal{M}$}. This also directly follows from the following calculations 
\begin{equation*}
\begin{array}{rcl}
\det (J_D)&=&\det \begin{pmatrix}
  F_{V_1} &  F_{V_2} & F_{w_2}\\
G_{V_1} & G_{V_2} & G_{w_2}\\
H_{V_1} & H_{ V_2} & H_{w_2}\\
\end{pmatrix} \Bigg\vert_{\RED{\mathcal{M}}} \vspace{0.1in}\\ &=& \det\begin{pmatrix}\label{det-desingu}
 F_{V_1} &  F_{V_2} & F_{w_2}\\
G_{V_1} & G_{V_2}&0\\
H_{V_1} & H_{ V_2}&0\\
\end{pmatrix} \Bigg\vert_{\RED{\mathcal{M}}} \vspace{0.1in} \\ &=& F_{w_2} \left( G_{V_1} H_{V_2}-  G _{V_2} H_{V_1} \right)\equiv 0,
\end{array}
\end{equation*}
\RED{where the entries in $J_D$ denote partial derivatives}. The last equality in the above equations holds because $G_{V_1} H_{V_2}\equiv G_{V_2} H_{V_1}=\delta f_2 g_2 F_{1V_1V_2} F_{1 V_1V_1}$ (see \eqref{desingu_slowreduced})\RED{, where $F_{1V_1V_1}=\frac{\partial^2 F_1}{\partial V_1^2}$ and $F_{1V_1V_2}=\frac{\partial^2 F_1}{\partial V_1 \partial V_2}$}. Hence the remaining two nontrivial eigenvalues $\lambda_w$ and $\lambda_s$ satisfy
\begin{equation*} 
    \lambda^2 - {\rm{tr}}(J_D)\lambda+ \frac{1}{2} \left( ({\rm{tr}}(J_D))^2-{\rm{tr}}(J_D^2) \right)=0.
\end{equation*}
It follows that an $\fsn$ bifurcation is given by 
\begin{equation*}
    \frac{1}{2}(\rm{tr}(J_D))^2-\rm{tr}(J_D^2)=F_{V_1}G_{V_2}-F_{V_2}G_{V_1}-F_{w_2}H_{V_1}=0.
\end{equation*}
Plugging in functions $F$, $G$ and $H$ from \eqref{desingu_slowreduced}, we can rewrite the above condition as the equation \eqref{FSNI} in subsection \ref{subsec:canard}, namely,
\begin{equation} \label{eq:app-FSN}
    \begin{array}{rcl}
     \fsn &:=& \{ (V_1,\,w_1,\,V_2,\,w_2) \in \RED{\mathcal{M}}: \vspace{0.1in} \\  && f_2 Q(V_1,\,V_2,\,w_2)=\delta P(V_1,\,V_2,\,w_2)\},
     \end{array}
\end{equation}
where $Q(V_1,\,V_2,\,w_2):=F_{V_1}F_{1 V_1 V_2}-F_{V_2} F_{1 V_1 V_1}$ and $P(V_1,\,V_2,\,w_2):=g_2 F_{1 V_2} f_{2 w_2} F_{1 V_1 V_1}$. 

Next we prove there are two ways an $\fsn$ can occur.
In case 1, suppose $Q\neq 0$. An $\fsn$ can occur when $f_2=\delta K_1(V_1,\,V_2,\,w_2)$, where $K_1(V_1,\,V_2,\,w_2) :=P/Q$. Note that $\fsn$ is a folded singularity point \eqref{folded}, that is, $F=g_1- F_{1 V_2}f_2=0$. This implies $g_1=F_{1 V_2}f_2=\delta K_2(V_1,\,V_2,\,w_2)$, where $K_2(V_1,\,V_2,\,w_2):=F_{1 V_2}K_1(V_1,\,V_2,\,w_2)$. Thus, the first way that an $\fsn$ can occur is described as the following.
\begin{equation*}
    \begin{array}{rcl}
    \mathrm{FSN}^1 &:=& \{ (V_1,\,w_1,\,V_2,\,w_2) \in \ls: \vspace{0.1in} \\ && f_2=\delta K_1(V_1,\,V_2,\,w_2), \ g_1= \delta K_2(V_1,\,V_2,\,w_2)\},
    \end{array}
\end{equation*}
which is the equation \eqref{FSN-1way} in subsection \ref{subsec:canard}. This implies that an $\fsn^1$ point becomes a CDH when $\delta=0$. 

In case 2, suppose $f_2\neq 0$. An $\fsn$ can occur when $Q(V_1,V_2,w_2)=\delta K_3(V_1,V_2,w_2)$
where $K_3(V_1,V_2,w_2) = P/f_2$. As a result, the second way that an $\fsn$ can occur is defined as
\begin{equation*}
    \begin{array}{rcl}
    \mathrm{FSN}^2 &:=& \{ (V_1,\,w_1,\,V_2,\,w_2) \in \RED{\mathcal{M}}: \vspace{0.1in} \\ && Q(V_1,\,V_2,\,w_2)= \delta K_3(V_1,\,V_2,\,w_2)\},
    \end{array}
\end{equation*}
which is the equation \eqref{FSN-2way} in subsection \ref{subsec:canard}. This implies that an $\fsn^2$ is $O(\delta)$ close to the intersection of $\RED{\mathcal{M}}$ and $Q(V_1,\,V_2,\,w_2)=0$. 

\section{CDH at double singular limit exhibits two linearly independent critical eigenvectors} \label{app:CDH-vc}
Recall the Jacobian matrix $J_D$ at an $\fsn^1$, which becomes a CDH at the double singular limit, has two zero eigenvalues. We prove that the center subspace associated with the two zero eigenvalues is two dimensional and is \RED{transverse} to the fold surface of the critical manifold (see subsection \ref{subsec:interaction} and Figure \ref{fig:weak-eigen-direc}).

The Jacobian matrix $J_D$ of the desingularized system at the double singular limit ($\varepsilon\to 0, \ \delta\to 0)$ is given by
\begin{equation*}
J_D=\begin{pmatrix}
 F_{V_1} &  F_{V_2} & F_{w_2}\\
G_{V_1} & G_{V_2}&0\\
0 & 0 &0\\
\end{pmatrix},
\end{equation*}
in which $G_{V_1}=-F_{1 V_1 V_1}f_2$ and $G_{V_2}=-F_{1 V_1 V_2}f_2$. At a CDH, $f_2=g_1=0$. It follows that $G_{V_1}=G_{V_2}=0$. Thus, the Jacobian matrix $J_D$ at a CDH singularity becomes
\begin{equation*}
J_D=\begin{pmatrix}
 F_{V_1} &  F_{V_2} & F_{w_2}\\
0 & 0&0\\
0 & 0&0\\
\end{pmatrix},
\end{equation*}
which has a nullity of 2, implying the existence of two linearly independent critical eigenvectors associated with zero eigenvalues. Moreover, the center subspace at the CDH, given by the plane $F_{V_1}V_1+ F_{V_2}V_2+ F_{w_2}w_2=0$ (see pink planes in Figure \ref{fig:weak-eigen-direc}C and D), is transverse to the fold surface $\ls$ of the critical manifold. This is because the normal vector of $\ls$ at the CDH, which is given by $\mathbf{n_f}=\frac{\mathbf{n_s}}{|\mathbf{n_s}|}$ where $\mathbf{n_s}=\left( F_{1 V_1 V_1}, F_{1 V_1 V_2}, 0\right)$ is not parallel to the normal vector of the center subspace given by $\left( F_{V_1}, F_{V_2}, F_{w_2}\right)$. Thus, the CDH singularity (an $\fsn^1$ at the double singular limit) considered in system \eqref{eq:main} is a saddle-node bifurcation of folded singularities, with the center manifold of the $\fsn^1$ transverse to the fold of the critical manifold. 

Finally, we verify a previous claim that we made in subsection \ref{subsec:interaction} that the eigenvector associated with the first trivial zero eigenvalue of $J_D$ (denoted by $\mathbf{v_0}$) is always tangent to $\ls$. Through direct computations, we obtain $\mathbf{v_0}=\frac{\mathbf{u_0}}{|\mathbf{u_0}|}$ with 
\begin{equation*}
\begin{array}{rcl}
\mathbf{u_0}&=& ( F_{V_2} F_{w_2} F_{1 V_1 V_2}, -F_{V_2} F_{w_2} F_{1 V_1 V_1} , \vspace{0.1in} \\ && \left( F_{V_2}\right)^2 F_{1 V_1 V_1} -  F_{V_1} F_{V_2} F_{1 V_1 V_2} ). 
\end{array}
\end{equation*}
The dot product of $\mathbf{v_0}$ and the unit normal vector of $\ls$ at the CDH is given by
\begin{equation*}
    \begin{array}{rcl}
\mathbf{v_0}\cdot \mathbf{n_f} &=& \frac{\left( F_{V_2} F_{w_2} F_{1 V_1 V_2} F_{1 V_1 V_1}-F_{V_2} F_{w_2} F_{1 V_1 V_1} F_{1 V_1 V_2} \right)}{|u_0||n_s|}  \vspace{0.1in} \\ &=& 0.
    \end{array}
\end{equation*}
It follows directly that $\mathbf{v_0}$ is tangent to $\ls$ as expected.

\section{Decreasing $\phi_2$ when $g_{\rm syn}=4.3$ causes non-monotonic effects on SAOs}\label{app:de-phi2-4p3}

\begin{figure*}[!htp]
    \begin{center}
		\begin{tabular}{@{}p{0.48\linewidth}@{\quad}p{0.48\linewidth}@{}}
			\subfigimg[width=\linewidth]{\bf{\small{(A)}}}{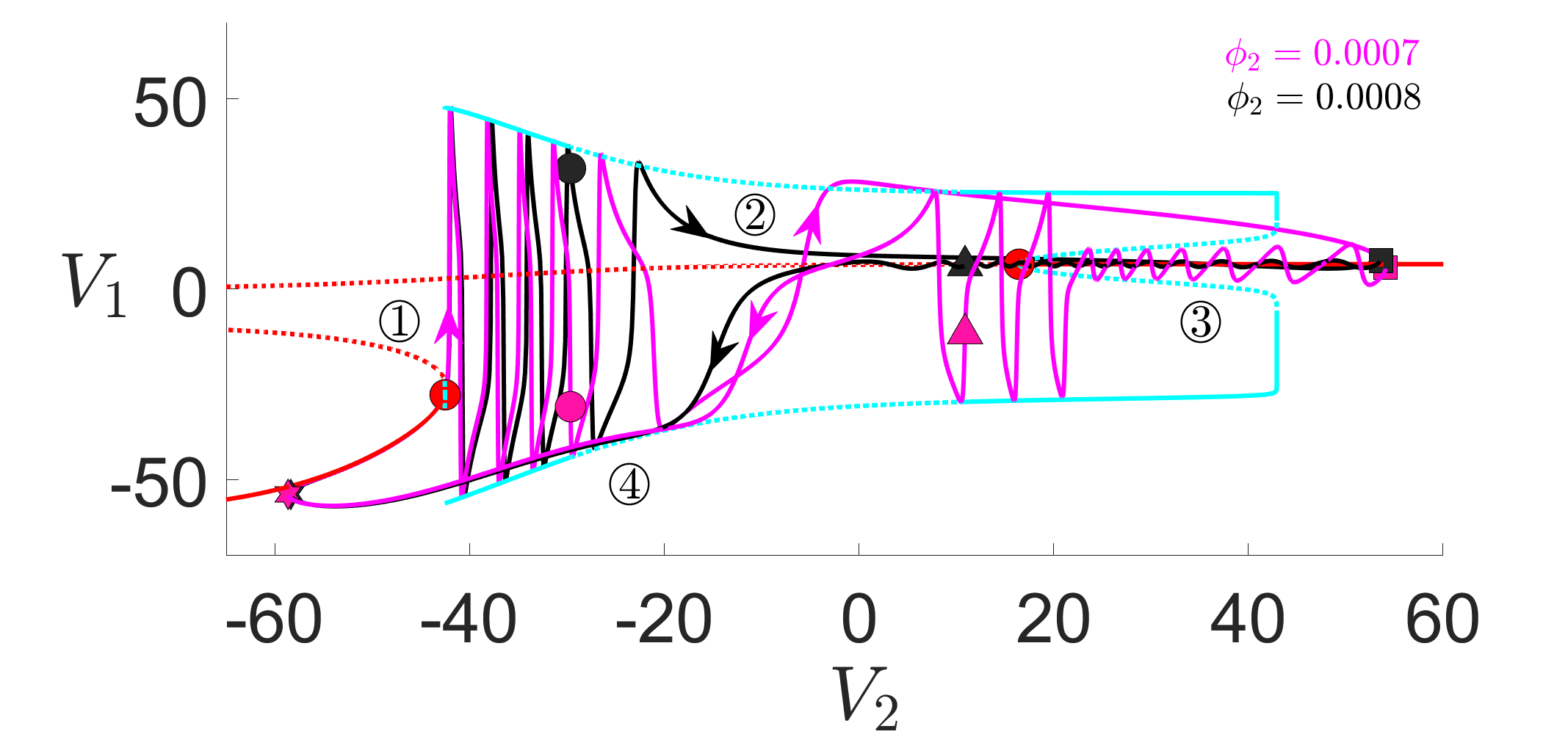}&
			\subfigimg[width=\linewidth]{\bf{\small{(B)}}}{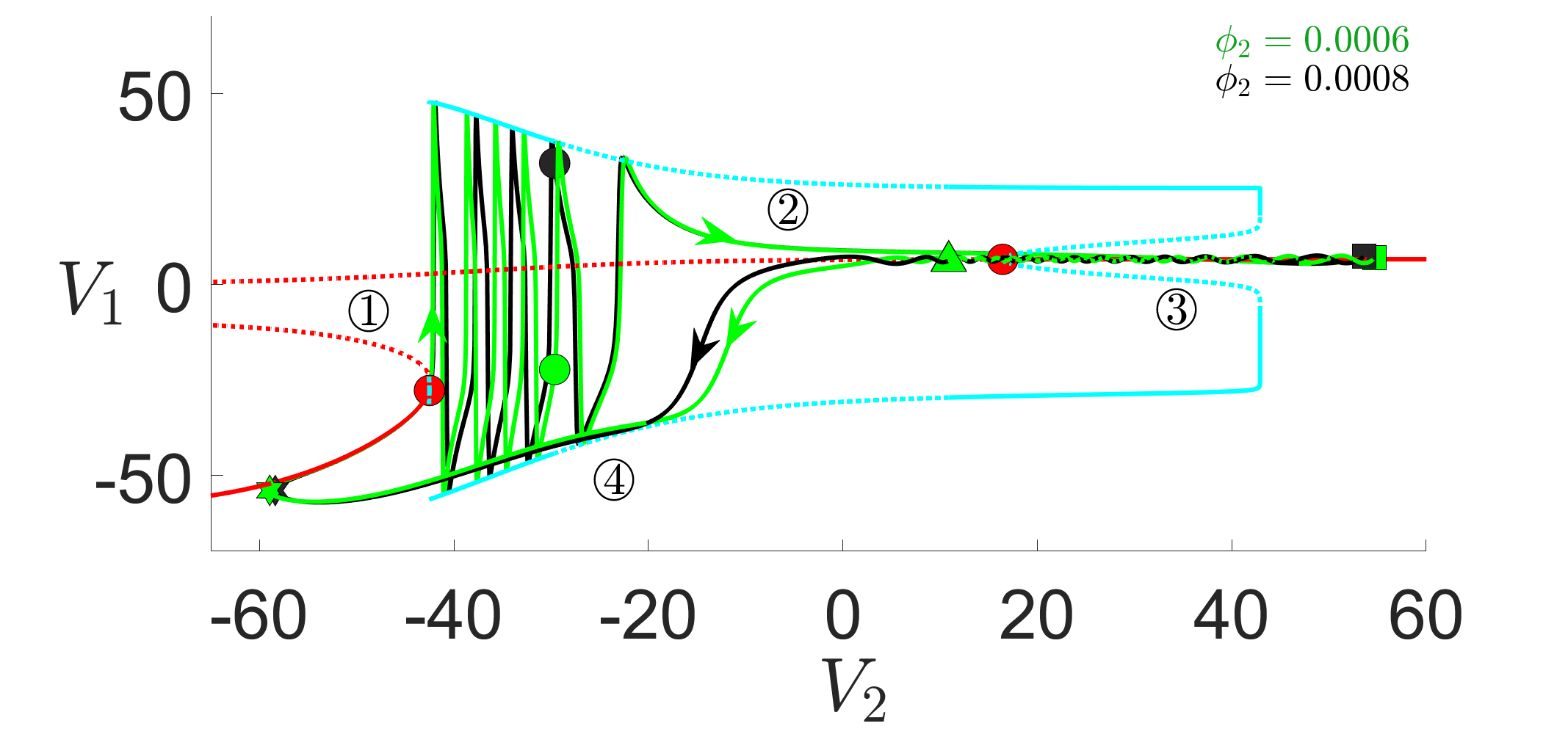}
		\end{tabular}
	\end{center}
\caption{Projections of solutions of \eqref{eq:main} with $g_{\rm syn}=4.3$, $C_1=8$, and different $\phi_2$ values onto $(V_2,V_1)$-plane.		
(A) Projections of trajectories with $\phi_2=0.0008$ (black) and $\phi_2=0.0007$ (magenta). (B) Projections of trajectories with $\phi_2=0.0008$ (black) and $\phi_2=0.0006$ (green).	The black and magenta stars, circles, squares and triangles have the same meaning as the green symbols in Figure \ref{fig:transition}. Other symbols and curves have the same meaning as in Figure \ref{fig:Ms-V1V2-gsyn4p3}.}
\label{fig:4p3-phi2-p8-p7}
\end{figure*}

In this subsection, we explain the mechanism underlying the non-monotonic effects of decreasing $\phi_2$ on SAOs (see Figure \ref{fig:gsyn_4p3-phi2}A in subsection \ref{sec:vary_eps_delta_mmo_4p3}). We claim that (1) if an additional (full) big spike is generated during the process of decreasing $\phi_2$, the number of SAOs will decrease and their amplitudes will increase (e.g., when $\phi_2$ decreases from $0.0008$ to $0.0007$); (2) if no additional big spike is generated during the process, the changes to the SAOs will be reversed (e.g., there are more SAOs with smaller amplitudes when $\phi_2$ decreases from $0.0007$ to $0.0006$). Below, we examine the two opposite effects using the projections of the corresponding solutions and bifurcation diagrams onto $(V_2, V_1)$ space (see Figure \ref{fig:4p3-phi2-p8-p7}). 

Figure \ref{fig:4p3-phi2-p8-p7}A shows that the amplitude of SAOs with $\phi_2=0.0008$ (black) is smaller than SAOs with $\phi_2=0.0007$ (magenta). There are also more black SAOs when $\phi_2=0.0008$ (see Figure \ref{fig:gsyn_4p3-phi2}A), which however is not obvious in the $(V_2, w_2, V_1)$ projection as the black SAOs are hardly visible. 
Compared with the black solution, the magenta one has a slower evolving rate of $V_2$, which is slaved to $w_2$, during phase \textcircled{1}. Thus, one extra full spike occurs for the magenta trajectory compared to the black trajectory during the jump at phase \textcircled{2}. As a result, the black trajectory, after its last big spike, approaches the black square at its maximum $V_2$ along $\mss$. In contrast, the magenta trajectory, due to the extra full spike that occurs during the jump, approaches the maximum $V_2$ from the outside of the periodic orbit branch. Consequently, the magenta trajectory that is being pushed further away from $\mss$ exhibits fewer small oscillations with larger amplitudes compared with the black trajectory, which remains close to $\mss$. 

When $\phi_2$ is reduced from $0.0007$ to $0.0006$, the solution trajectory changes to the green curve in Figure \ref{fig:4p3-phi2-p8-p7}B. Slowing down of $w_2$ even further keeps the number of full spikes between the magenta and green trajectories the same, but now all five full green spikes occur within  phase \textcircled{1}, similar to the case when $\phi_2=0.0008$. As a result, the green and black SAOs exhibit similar characteristics as shown in Figure \ref{fig:4p3-phi2-p8-p7}B, resulting in more SAOs with smaller amplitudes than the magenta solution in Figure \ref{fig:4p3-phi2-p8-p7}A. 

\section{Increasing $\phi_2$ when $g_{\rm syn}=4.3$ leads to multiple MMOs/non-MMOs transitions}\label{app:in-phi2-4p3}

\begin{figure*}[!t]
    \begin{center}
        \begin{tabular}{@{}p{0.48\linewidth}@{\quad}p{0.48\linewidth}@{}}
            \multicolumn{2}{c}{\subfigimg[width=0.5\linewidth]{\bf{\small{(A)}}}{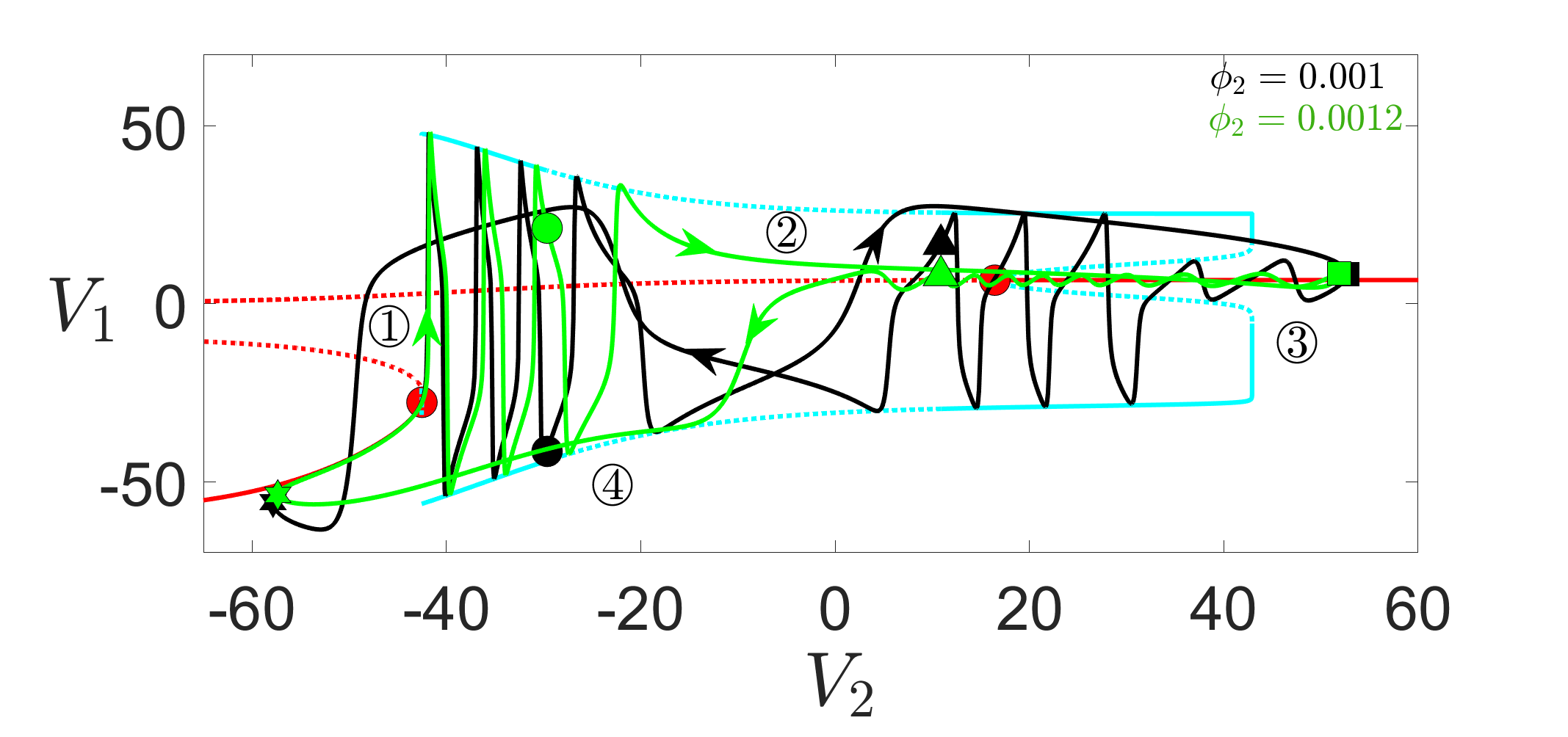}}\\
            \subfigimg[width=\linewidth]{\bf{\small{(B)}}}{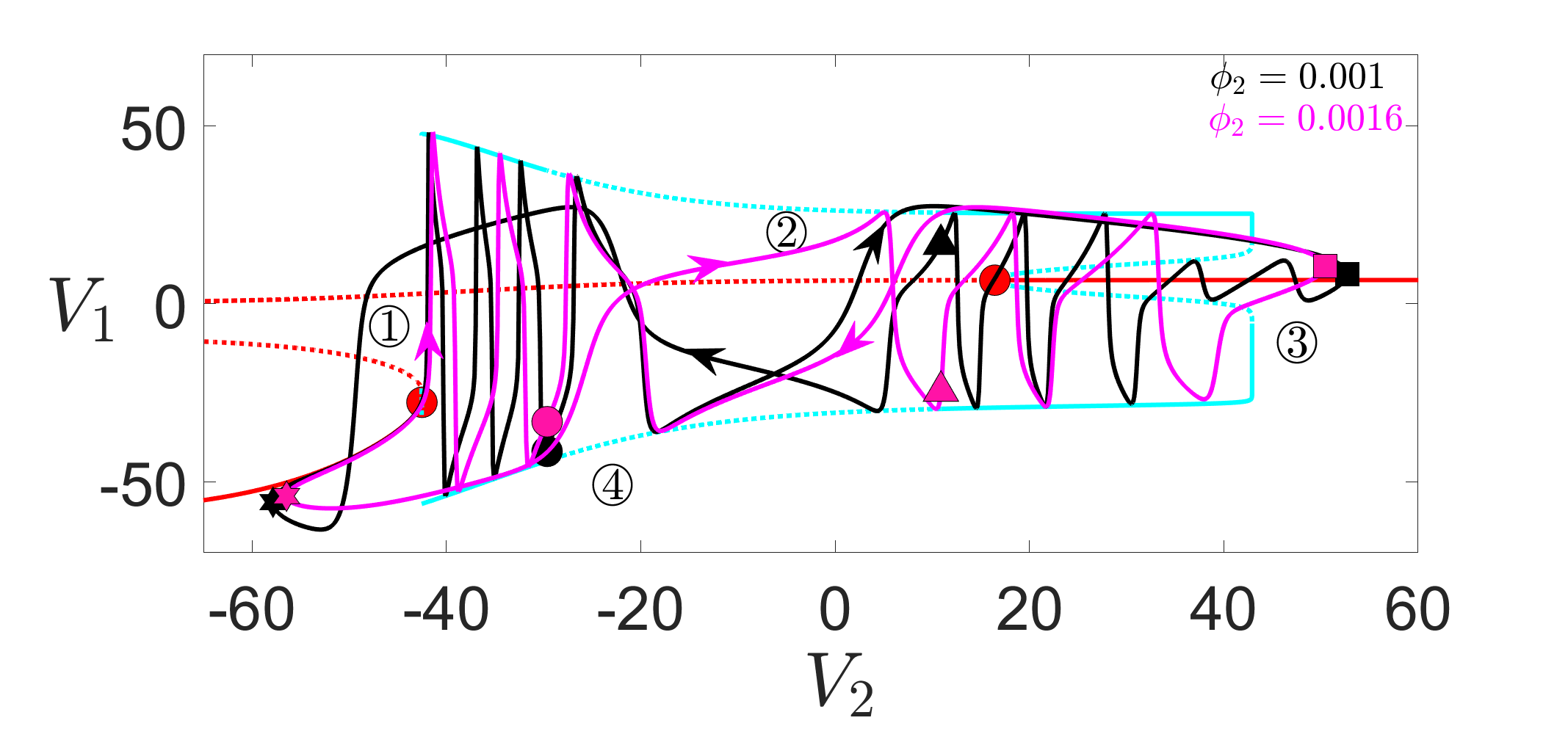}&
            \subfigimg[width=\linewidth]{\bf{\small{(C)}}}{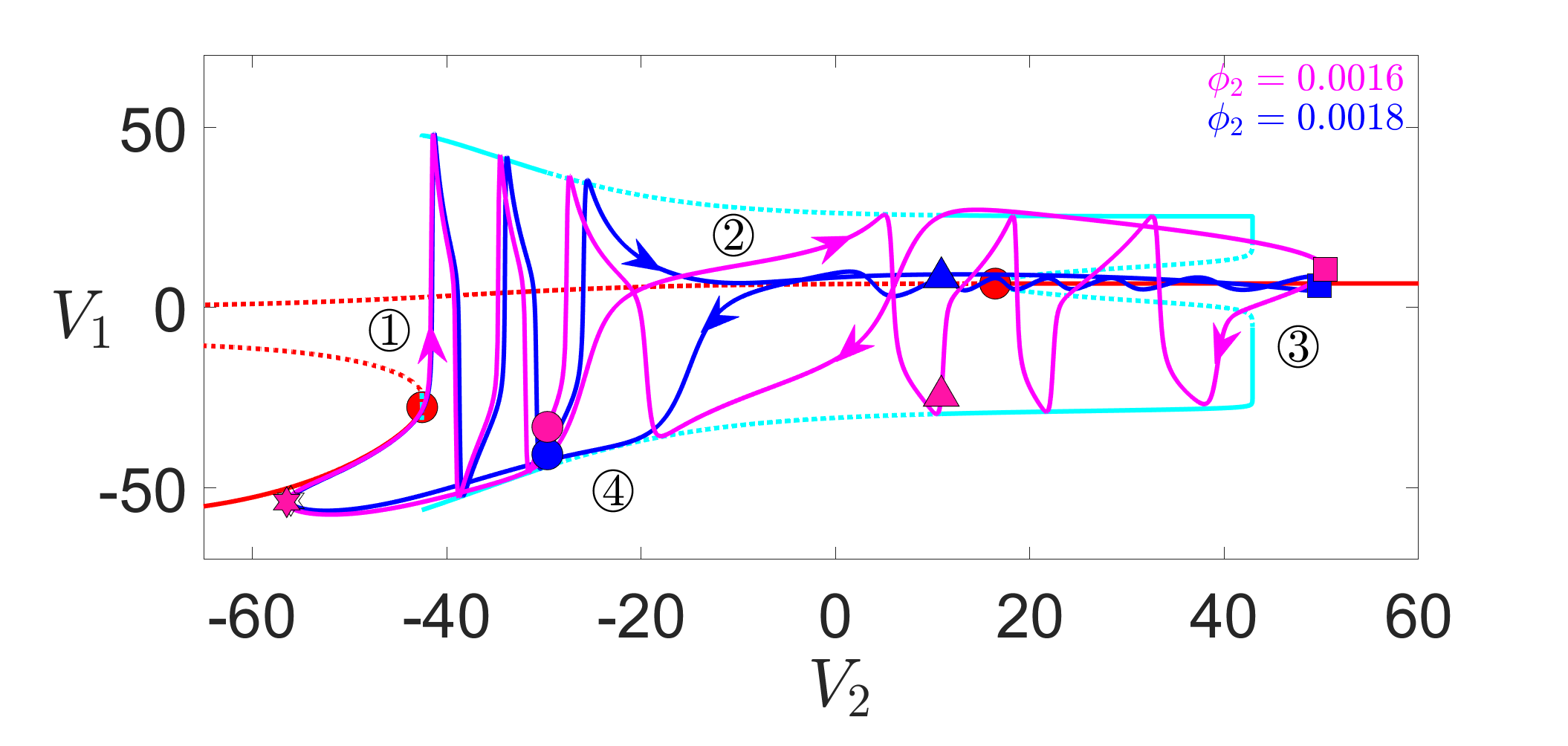}
        \end{tabular}
    \end{center}
\caption{Projections of solutions \eqref{eq:main} with $g_{\rm syn}=4.3$, $C_1=8$ ($\epsilon =0.1$) and different values of $\phi_2$ onto the $(V_2, V_1)$-space. Trajectories with $\phi_2=0.0012$, $\phi_2=0.0016$ and $\phi_2=0.0018$ are shown by green, magenta and blue curves, respectively. 
Other symbols and color codings have the same meaning as in Figure \ref{fig:Ms-V1V2-gsyn4p3}.}
\label{fig:4p3-increasing-phi2} 
\end{figure*}

In this subsection, we explain the mechanism underlying the MMOs/non-MMOs transitions induced by an increase of $\phi_2$ in the case of $g_{\rm syn}=4.3$ (see Figure \ref{fig:gsyn_4p3-phi2} in subsection \ref{sec:vary_eps_delta_mmo_4p3}), which is similar to the mechanism that causes the non-monotonic effects on SAOs when $\phi_2$ is decreased (see Appendix \ref{app:de-phi2-4p3}).

During the increase of $\phi_2$ which speeds up $w_2$, the big spikes produced during phase \textcircled{1} will gradually decrease. 
If one big (full) spike before phase \textcircled{3} is lost during this process, there will be more SAOs with smaller amplitudes (e.g., when $\phi_2$ increases from $0.001$ to $0.0012$ in Figure \ref{fig:4p3-increasing-phi2}A) or the earlier lost MMOs will reappear (e.g., when $\phi_2$ increases from 0.0016 to 0.0018 in Figure \ref{fig:4p3-increasing-phi2}C). As $\phi_2$ increases from $0.001$ to $0.0012$ in Figure \ref{fig:4p3-increasing-phi2}A, the number of full spikes during phases \textcircled{1} and \textcircled{2} decreases by one. Consequently, all three green full spikes occur within phase \textcircled{1}, while there is a large black spike during phase \textcircled{2}. As explained in the previous subsection, this causes the green trajectory to approach the maximum $V_2$ along $\mss$ and remain close to it after reaching the green square. As a result, more SAOs with smaller amplitudes are observed in the green trajectory compared to the black trajectory. Similarly, when the increase of $\phi_2$ from $0.0016$ to $0.0018$ leads to a decrease in the number of big spikes that occur before phase \textcircled{3}, more SAOs are observed and the transition from non-MMOs to MMOs occurs (see Figure \ref{fig:4p3-increasing-phi2}C).    

In contrast, increasing $\phi_2$ from $0.0012$ to $0.0016$ does not change the number of large spikes before phase \textcircled{3} (compare the green and magenta trajectories in Figure \ref{fig:4p3-increasing-phi2}A and B). However, speeding up of $\phi_2$ pushes the third big magenta spike to occur during the jump, similar to the black trajectory. This causes the trajectory to be pushed away from $\mss$, resulting in fewer SAOs with larger amplitudes. In fact, in this case, SAOs in the magenta trajectory are eliminated and therefore a transition from MMOs to non-MMOs results.


\section*{References}
\bibliographystyle{plain}
\bibliography{bib-file}
\nocite{*}

\end{document}